\documentclass[11pt]{amsart}
\usepackage{amssymb,amscd,amsmath,amsthm}

\newcommand\cA{{\mathcal{A}}}
\newcommand\cB{{\mathcal{B}}}

\newcommand\cD{{\mathcal{D}}}
\newcommand\cE{{\mathcal{E}}}
\newcommand\cF{{\mathcal{F}}}
\newcommand\cG{{\mathcal{G}}}
\newcommand\cH{{\mathcal{H}}}

\newcommand\cK{{\mathcal{K}}}
\newcommand\cL{{\mathcal{L}}}
\newcommand\cM{{\mathcal{M}}}
\newcommand\cN{{\mathcal{N}}}

\newcommand\cR{{\mathcal{R}}}
\newcommand\cS{{\mathcal{S}}}

\newcommand\cU{{\mathcal{U}}}
\newcommand\cV{{\mathcal{V}}}

\newcommand\cX{{\mathcal{X}}}

\newcommand\RR{{\mathbb{R}}}
\newcommand\CC{{\mathbb{C}}}
\newcommand\ZZ{{\mathbb{Z}}}
\newcommand\NN{{\mathbb{N}}}
\newcommand\QQ{{\mathbb{Q}}}

\newcommand{\supp}{\operatorname{supp}}
\newcommand{\codim}{\operatorname{codim}}

\newcommand{\dom}{\operatorname{dom}}

\newcommand{\ev}{\text{ev}}
\newcommand{\odd}{\text{odd}}

\newcommand{\tr}{{\rm tr}\;}
\newcommand{\trLambda}{{\rm tr}_\Lambda}

\newcommand{\Tr}{{\rm Tr}\;}

\newcommand{\res}{{\rm res}\;}

\newcommand\Ker{\mathop{\rm Ker}\nolimits}
\newcommand\Dom{\mathop{\rm Dom}\nolimits}
\newcommand{\bint}{\ensuremath{-\hspace{-2,4ex}\int}}

        \numberwithin{equation}{section}

\theoremstyle{plain}
\newtheorem{thm}{Theorem}[section]

\newtheorem{prop}[thm]{Proposition}

\newtheorem{theorem}[thm]{Theorem}

\newtheorem{proposition}[thm]{Proposition}

\theoremstyle{definition}

\newtheorem{defn}[thm]{Definition}

\newtheorem{ex}[thm]{Example}

\theoremstyle{remark}

\theoremstyle{plain}
\begin{document}
\bibliographystyle{plain}

\title[Noncommutative geometry of foliations]{Noncommutative geometry of foliations}
\author{Yuri A. Kordyukov}
\address{Institute of Mathematics, Russian Academy of Sciences, Ufa} \email{yuri@imat.rb.ru}
\thanks{Supported by Russian Foundation of Basic Research
(grant no. 04-01-00190)} \dedicatory{Dedicated to Yu. Solovyov}

\begin{abstract}
We review basic notions and methods of noncommutative geometry and
their applications to analysis and geometry on foliated manifolds.
\end{abstract}

\maketitle

\tableofcontents

\section{Introduction}
The starting point of noncommutative geometry is the passage from
geometric spaces to algebras of functions on these spaces with the
subsequent translation of basic analytic and geometric concepts
and constructions on geometric spaces to the algebraic language
and their extension to general noncommutative algebras. Such a
procedure is well-known and was applied for a long time, for
instance, in algebraic geometry, where it is related with the
study of commutative algebras. It is also well known that the
theory of $C^*$-algebras is a far-reaching generalization of the
theory of topological spaces, and the theory of von Neumann
algebras is a generalization of the classical measure and
integration theory.

The main purpose of noncommutative differential geometry, which
was initiated by Connes \cite{Co:nc} and is actively developing at
present time (cf. the recent surveys
\cite{Connes2000,ConnesLNM1831} and the books
\cite{Co,Gracia:book,Landi:book} in regard to different aspects of
noncommutative geometry), consists in the extension of the methods
described above to the analytic objects on geometric spaces and to
the noncommutative algebras. Here the main attention is focused on
that, first, a correct noncommutative generalization applied in
the classical setting, that is, to an algebra of functions on a
compact manifold should agree with its classical analogue, and,
second, it should inherit nice algebraic and analytic properties
of its classical analogue. Nevertheless, it should be noted that,
as a rule, such noncommutative generalizations are quite
nontrivial and have richer structure and essentially new features
than their commutative analogues.

Noncommutative geometry lies on the border of functional analysis
and differential geometry and is of great importance for these
areas of mathematics. On the one hand, the development of
geometric methods in the operator theory and the theory of
operator algebras allows to use fruitfully geometric intuition for
the investigation of various problems of abstract functional
analysis. On the other hand, there are many examples of singular
geometric spaces, which are badly described from the point of view
of classical measure theory and resist to the study by usual
``commutative'' methods of geometry, topology and analysis, but
one can naturally associate to them a noncommutative algebra,
which can be considered as an analogue of the algebra of
(measurable, continuous, smooth and so on) functions on the given
geometric object. Let us give some examples of such singular
objects:
\begin{enumerate}
  \item Manifolds with singularities, for instance, manifolds with isolated conic
  singular points.
  \item Discrete spaces.
  \item The dual space to a group (discrete or a Lie group).
  \item Cantor sets.
  \item The orbit space of a group action on a manifold.
  \item The leaf space of a foliation.
\end{enumerate}
Use of the notions and methods of noncommutative geometry for a
noncommutative algebra, being an analogue of an algebra of
functions on a singular geometric space, allows in many cases,
first of all, just to define some reasonable analytic and
geometric objects associated with the given space and pose
sensible problems, that, in its turn, gives possibility to apply
properties of these objects for getting an information about
geometry of this space.

This paper contains a short exposition of the methods of
noncommutative geometry as applied to the study of one of the
fundamental examples of noncommutative geometry, namely, the leaf
spaces of a foliation on a smooth manifold, or, that is the same,
to the study of the transverse structure of foliations. Our
purpose is to give a survey of the basic notions and methods of
noncommutative geometry, to show how one can associate various
objects of noncommutative geometry to foliations and how these
noncommutative analogues of classical notions are related with
classical objects on foliated manifolds, and to describe
applications of the methods of noncommutative geometry to the
study of geometry of foliations.

The author is grateful to N.I. Zhukova for useful remarks.

\section{Background information on foliation theory}\label{rev}
\subsection{Foliations}\label{s:fol}
In this Section, we recall the definition of a foliated manifold
and some notions related to foliations (concerning to different
aspects of the foliation theory see, for instance,
\cite{Camacho,Candel-Conlon1,Candel-Conlon2,Godbillon,MMrcun,Molino,M-S,Re,Tondeur}).

Let $M$ be a smooth manifold of dimension $n$. (Here and later on,
``smooth'' means ``of class $C^\infty$''. We will always assume
that all our objects under consideration are of class $C^\infty$.)

\begin{defn}
(1) An atlas $\cA=\{(U_i,\phi_i)\}$, where $\phi_i: U_i\subset
M\rightarrow \RR^n$, of the manifold $M$ is called an atlas of a
foliation of dimension $p$ and codimension $q$ ($p\leq n, p+q=n$),
if, for any $i$ and $j$ such that $U_i\cap U_j\not=\varnothing$,
the coordinate transformation $$
\phi_{ij}=\phi_i\circ\phi_j^{-1}:\phi_j(U_i\cap U_j)\subset
\RR^p\times\RR^q \to \phi_i(U_i\cap U_j)\subset \RR^p\times\RR^q$$
has the form $$ \phi_{ij}(x,y)=(\alpha_{ij}(x,y), \gamma_{ij}(y)),
\quad (x,y)\in \phi_j(U_i\cap U_j)\subset \RR^p\times\RR^q. $$

(2) Two atlases of a dimension $p$ foliation are equivalent, if
their union is again an atlas of a dimension $p$ foliation.
\medskip\par
(3) A manifold $M$ endowed with an equivalence class $\cF$ of
atlases of a dimension $p$ foliation is called a manifolds with a
dimension $p$ foliation (or a foliated manifold).
\end{defn}

An equivalence class $\cF$ of foliation atlases is also called a
complete atlas of a foliation. We will also say that $\cF$ is a
foliation on $M$.

A pair $(U,\phi)$ that belongs to the atlas of the foliation $\cF$
and also the corresponding map $\phi:U\rightarrow \RR^n$ are
called a foliated chart of the foliation $\cF$, and $U$ is called
a foliated coordinate neighborhood.

Let $\phi:U\subset M \rightarrow \RR^n$ be a foliated chart. The
connected components of the set $\phi^{-1}(\RR^p\times \{y\}),
y\in \RR^q$ are called the plaques of the foliation $\cF$.

Plaques of $\cF$ given by all possible foliated charts form a base
of a topology on $M$. This topology is called the leaf topology on
$M$. We will also denote by $\cF$ the set $M$ endowed with the
leaf topology. One can introduce the structure of a
$p$-dimensional manifold on $\cF$.

Connected components of the manifold $\cF$ are called leaves of
the foliation $\cF$. Leaves are (one-to-one) immersed
$p$-dimensional submanifolds in $M$. For any $x\in M$, there is a
unique leaf, passing through $x$. We will denote this leaf by
$L_x$.

One can give the equivalent definition of a foliation, saying,
that there is given a foliation $\cF$ of dimension $p$ on a
manifold $M$ of dimension $n$, if $M$ is represented as a disjoint
union of a family $\{L_\lambda : \lambda\in \cL\}$ of connected,
(one-to-one) immersed submanifolds of dimension $p$, and there is
an atlas $\cA=\{(U_i,\phi_i)\}$ of the manifold $M$ such that, for
any $i$ and for any $\lambda\in \cL$, the connected components of
the set $L_\lambda\cap U_i$ are given by equations of the form
$y=\text{const}$.

\begin{ex}
Let $M$ be a smooth manifold of dimension $n$, $B$ a smooth
manifold of dimension $q$ and $\pi: M\to B$ a submersion (that is,
the differential $d\pi_x:T_xM\to T_{\pi(x)}B$ is surjective for
any $x\in M$). The connected components of the pre-images of
points of $B$ under the map $\pi$ determine a codimension $q$
foliation on $M$, which is called the foliation determined by the
submersion $\pi$. If, in addition, the pre-images $\pi^{-1}(b),
b\in B,$ are connected, the foliation is called simple.
\end{ex}

\begin{ex}
Let $X$ be a nonsingular (that is, a non-vanishing) smooth vector
field on a compact manifold $M$. Then its phase curves form a
codimension one foliation.

More generally, suppose that a connected Lie group $G$ acts
smoothly on a smooth manifold $M$, and, moreover, the dimension of
the isotropy subgroup $G_x=\{x\in G: gx=x\}$ does not depend on
$x\in M$. In particular, one can assume that the action is locally
free, that means $G_x$ is discrete for any $x\in M$. Then the
orbits of the Lie group $G$ action define a foliation on $M$.
\end{ex}

\begin{ex} \emph{Linear foliation on torus.}
Consider a vector field $\tilde{X}$ on $\RR^2$ given by
\[
\tilde{X}=\alpha\frac{\partial}{\partial
x}+\beta\frac{\partial}{\partial y}
\]
with constant $\alpha$ and $\beta $. Since $\tilde{X}$ is
invariant under all translations, it determines a vector field $X$
on the two-dimensional torus ${T}^2={\RR}^2/{\ZZ}^2$. The vector
field $X$ determines a foliation $\cF$ on ${T}^2$. The leaves of
$\cF$ are the images of the parallel lines
$\tilde{L}=\{(x_0+t\alpha, y_0+t\beta): t\in\RR\}$ with the slope
$\theta=\beta/\alpha $ under the projection $\RR^2\to T^2$.

In the case when $\theta$ is rational, all leaves of $\cF$ are
closed and are circles, and the foliation $\cF$ is determined by
the fibers of a fibration $T^2\to S^1$. In the case when $\theta$
is irrational, all leaves of $\cF$ are everywhere dense in $T^2$.
\end{ex}

\begin{ex} \emph{Homogeneous foliations.}
Let $G$ be a Lie group and $H\subset G$ its connected Lie
subgroup. The family $\{ gH : g\in G\}$ of right cosets of $H$
forms a foliation $\cH$ on $G$. If $H$ is a closed subgroup, then
$G/H$ is a manifold and $\cH$ is a foliation, whose leaves are
given by the fibers of the fibration $\pi : G\to G/H$.

Moreover, suppose that $\Gamma\subset G$ is a discrete subgroup
$G$. Then the set $M=\Gamma\backslash G$ of left cosets of $\Gamma
$ is a manifold of the same dimension as $G$. If $\Gamma $ is
cocompact in $G$, $M$ is compact. In any case, because $\cH$ is
invariant under the left translations, and $\Gamma $ acts on the
left, the foliation $\cH$ is mapped by the map $G\to
M=\Gamma\backslash G$ to a well-defined foliation $\cH_\Gamma $ on
$M$, which is often denoted by $\cF(G,H,\Gamma)$ and is called a
locally homogeneous foliation. The leaf of $\cH_\Gamma $ through a
point $\Gamma g \in M$ is diffeomorphic to $H/(g\Gamma g^{-1}\cap
H)$.
\end{ex}

\begin{ex} \emph{The Reeb foliation.}
Let us describe a classical construction of a codimension 1
foliation on the three-dimensional sphere $S^3$ due to Reeb
\cite{Reeb52}. Let $D^2$ denote the disk $\{(x,y)\in \RR^2 :
x^2+y^2\leq 1\}$. We start with the foliation in the cylinder
$\{(x,y,z)\in \RR^3 : x^2+y^2\leq 1\}=D^2\times \RR$, whose leaves
are the boundary of the cylinder $x^2+y^2=1$ and surfaces
$z=c+\exp (1/(1-x^2-y^2))$ with an arbitrary constant $c$. Since
this foliation is invariant under translations in $z$, it is
mapped by the standard projection $D^2\times \RR \to D^2\times
\RR/\ZZ=D^2\times S^1$ to a foliation $\cF_R$ on the solid torus
$D^2\times S^1$. Finally, observe that the standard
three-dimensional sphere $S^3=\{x=(x_1,x_2,x_3,x_4)\in \RR^4 :
x_1^2+x_2^2+x_3^2+x^2_4=1\}$ is obtained by gluing along the
boundary of two copies of $D^2\times S^1$. More precisely,
$S^3=S^3_1\cup S^3_2$, where
\[
S^3_1=\{x\in S^3: x_1^2+x_2^2 \leq x_3^2+x^2_4\}, \quad
S^3_2=\{x\in S^3: x_1^2+x_2^2 \geq x_3^2+x^2_4\}.
\]
A diffeomorphism $S^3_1\stackrel{\cong}{\longrightarrow} D^2\times
S^1$ is given by
\[
x\in S^3_1 \mapsto \left(\frac{x_1}{\sqrt{x_3^2+x^2_4}},
\frac{x_2}{\sqrt{x_3^2+x^2_4}}, \frac{x_3}{\sqrt{x_3^2+x^2_4}},
\frac{x_4}{\sqrt{x_3^2+x^2_4}}\right) \in D^2\times S^1,
\]
where we consider $S^1$ as $\{(x,y)\in\RR^2: x^2+y^2=1\}$, and its
inverse has the form
\[
(y_1,y_2,y_3,y_4)\in D^2\times S^1 \mapsto \left(\frac{y_1}{|y|},
\frac{y_2}{|y|}, \frac{y_3}{|y|},\frac{y_4}{|y|}\right)\in S^3_1,
\]
where $|y|^2=y_1^2+y_2^2+y_3^2+y^2_4$. Similar formulas can be
written for $S^3_2$. Take the Reeb foliations $\cF_R$ on $S^3_1$
and $S^3_2$. Since the boundaries $\partial S^3_1$ and $\partial
S^3_2$ are compact leaves, there is a well-defined foliation on
the union of $S^3_1$ and $S^3_2$, that is, on the sphere $S^3$,
which is called the Reeb foliation on $S^3$. The only compact leaf
of the Reeb foliation is the common boundary $\partial
S^3_1=\partial S^3_2$ of $S^3_1$ and $S^3_2$, diffeomorphic to the
two-dimensional torus $T^2$. The other leaves are diffeomorphic to
the plane $\RR^2$.
\end{ex}

\begin{ex} \emph{Suspension.}
Let $B$ be a connected manifold and $\tilde{B}$ its universal
cover equipped with the action of the fundamental group
$\Gamma=\pi_1(B)$ by deck transformations. Suppose that there is
given a homomorphism $\phi : \Gamma \to \operatorname{Diff}(F)$ of
$\Gamma$ to the group $\operatorname{Diff}(F)$ of diffeomorphisms
of a smooth manifold $F$. Define a manifold $ M=\tilde{B}
\times_\Gamma F$ as the quotient of the manifold $\tilde{B}\times
F$ by the action of $\Gamma$ given, for any $\gamma\in \Gamma$, by
\[
\gamma (b, f)=(\gamma b, \phi(\gamma)f),\quad (b, f) \in
\tilde{B}\times F.
\]
There is a natural foliation $\cF$ on $M$, whose leaves are the
images of the sets $\tilde{B}\times \{f\}, f\in F,$ under the
projection $\tilde{B}\times F\to M$. If, for any $\gamma\in
\Gamma, \gamma\not= e$, the diffeomorphism $\phi(\gamma)$ has no
fixed points, all leaves of $\cF$ are diffeomorphic to
$\tilde{B}$.

There is defined the bundle $\pi : M\to B : [(b, f)]\mapsto b\mod
\Gamma $ such that the leaves of $\cF$ are transverse to the
fibers of $\pi$. The bundle $\pi :M\to B$ is often said to be a
flat foliated bundle.
\end{ex}

A foliation $\cF$ determines a subbundle $F=T{\mathcal F}$ of the
tangent bundle $TM$, called the tangent bundle to $\cF$. It
consists of all vectors, tangent to the leaves of $\cF$. Denote by
$\cX(M)=C^\infty(M,TM)$ the Lie algebra of all smooth vector
fields on $M$ with the Lie bracket and by $\cX(\cF)=C^\infty(M,F)$
the subspace of vector fields on $M$, tangent to the leaves of
$\cF$ at each point. The subspace $\cX(\cF)$ is a subalgebra of
the Lie algebra $\cX(M)$. Moreover, by the Frobenius theorem, a
subbundle $E$ of $TM$ is the tangent bundle to  a foliation if and
only if it is involutive, that is, the space of sections of this
bundle is a Lie subalgebra of the Lie algebra $\cX(M)$: for any
$X,Y\in C^\infty(M,E)$ we have $[X,Y]\in C^\infty(M,E)$.

Let us introduce the following objects:
\begin{itemize}
  \item $\tau=TM/T\cF$ is the normal bundle to ${\mathcal F}$;
  \item $P_\tau: TM\to \tau$ is the natural projection;
  \item $N^*{\mathcal F}=\{\nu\in T^*M:\langle\nu,X\rangle=0\ \mbox{\rm for any}\
X\in F\}$ is the conormal bundle to ${\mathcal F}$.
\end{itemize}
Usually, we will denote by $(x,y)\in I^p\times I^q$ ($I=(0,1)$ is
the open interval) the local coordinates in a foliated chart
$\phi: U\to I^p \times I^q$ and by  $(x,y,\xi,\eta) \in I^p \times
I^q \times {\RR}^p\times {\RR}^q$ the local coordinates in the
corresponding chart on $T^*M$. Then the subset $N^*{\mathcal
F}\cap \pi^{-1}(U)=U_1$ (here $\pi: T^*M\to M$ is the bundle map)
is given by the equation $\xi=0$. Therefore, $\phi$ determines
naturally a foliated chart $\phi_n: U_1\to I^p \times I^q\times
{\RR}^q$ on $N^*{\mathcal F}$ with the coordinates $(x,y,\eta)$.

Any $q$-dimensional distribution $Q\subset TM$ such that
$TM=F\oplus Q$ is called a distribution transverse to the
foliation or a connection on the foliated manifold $(M,\cF)$. Any
Riemannian metric $g$ on $M$ determines a transverse distribution
$H$, which is given by the orthogonal complement of $F$ with
respect to the given metric: $H=F^{\bot}=\{X\in TM : g(X,Y)=0\
\mbox{\rm for any}\ Y\in F\}$.

\begin{defn}
A vector field $V$ on a foliated manifold $(M,\cF)$ is called an
infinitesimal transformation of $\cF$ if $[V,X]\in \cX(\cF)$ for
any $X\in \cX(\cF)$.
\end{defn}

The set of all infinitesimal transformations of $\cF$ is denoted
by $\cX(M/\cF)$. If $V\in \cX(M/\cF)$ and $T_t:M\to M, t\in \RR$
is the flow of the vector field $V$, then the diffeomorphisms
$T_t$ are automorphisms of the foliated manifold $(M,\cF)$, that
is, they take each leaf of $\cF$ to, possibly, another leaf.

\begin{defn}
A vector field $V$ on a foliated manifold $(M,\cF)$ is called
projectable, if its normal component $P_\tau(V)$ is locally the
lift of a vector field on the local base.
\end{defn}

In other words, a vector field $V$ on $M$ is projectable, if in
any foliated chart with local coordinates $(x,y), x\in \RR^p, y\in
\RR^q,$ it has the form $$
V=\sum_{i=1}^pf^i(x,y)\frac{\partial}{\partial x^i} + \sum_{j=1}^q
g^j(y)\frac{\partial}{\partial y^j}. $$

There is a natural action of the Lie algebra $\cX(\cF)$ on
$C^\infty(M,\tau)$. The action of $X\in \cX(\cF)$ on $N\in
C^\infty(M,\tau)$ is given by
\[
\theta(X) N=P_\tau[X,\widetilde{N}],
\]
where $\widetilde{N}\in \cX(M)$ is any vector field on $M$ such
that $P_\tau(\widetilde{N})= N$. A vector field $N\in \cX(M)$ is
projectable if and only if its transverse component $P_\tau(N)\in
C^\infty(M,\tau)$ is invariant under the $\cX(\cF)$-action
$\theta$. From here, one can easily see that a vector field on a
foliated manifold is called projectable if and only if it is an
infinitesimal transformation of the foliation.

\subsection{Holonomy} Let $(M, \cF)$ be a foliated manifold. The holonomy map
is a generalization to the case of foliations of the first return
map (or the Poincar\'e map) for flows.

\begin{defn}
A smooth transversal is a compact $q$-dimensional manifold $T$,
possibly with boundary, and an embedding $i:T\rightarrow M$, whose
image is everywhere transverse to the leaves of ${\mathcal F}$:
$T_{i(t)}i(T)\oplus T_{i(t)}\cF=T_{i(t)}M$ for any $t\in T$.
\end{defn}

We will identify a transversal $T$ with the image $i(T)\subset M$.

\begin{defn}
A transversal is complete, if it meets every leaf of the
foliation.
\end{defn}

Take an arbitrary continuous leafwise path $\gamma$ with endpoints
$\gamma(0)=x$ and $\gamma(1)=y$. (We will call a path $\gamma
:[0,1]\to M$ leafwise, if its image $\gamma([0,1])$ is contained
entirely in one leaf of the foliation). Let $T_0$ and $T_1$ are
smooth transversals such that $x\in T_0$ and $y\in T_1$.

Choose a partition $t_0=0<t_1<\ldots < t_k=1$ of the segment
$[0.1]$ such that for any $i=1,\ldots, k$ the curve
$\gamma([t_{i-1},t_i])$ is contained in some foliated chart $U_i$.
Shrinking, if necessary, the neighborhoods $U_1$ and $U_2$, one
can assume that, for any plaque $P_1$ of $U_1$, there is a unique
plaque $P_2$ of $U_2$, which meets $P_1$. Shrinking, if necessary,
the neighborhoods $U_1$, $U_2$ and $U_3$, one can assume that, for
any plaque $P_2$ of $U_2$, there is a unique plaque $P_3$ of
$U_3$, which meets $P_2$ and so on. After all, we get a family
$\{U_1, U_2, \ldots, U_k\}$ of foliated coordinate neighborhoods,
which covers the curve $\gamma([0,1])$, such that, for any
$i=1,\ldots, k$ and for any plaque $P_{i-1}$ of $U_{i-1}$, there
is a unique plaque $P_i$ of $U_i$, which meets $P_{i-1}$. In
particular, we get a one-to-one correspondence between the plaques
of $U_1$ and the plaques of $U_k$.

The smooth transversal $T_0$ determines a parametrization of the
plaques of $U_1$ near $x$. Accordingly, a smooth transversal $T_1$
determines a parametrization of the plaques of $U_k$ near $y$.
Taking into account one-to-one correspondence between the plaques
of $U_1$ and $U_k$ constructed above, we get a diffeomorphism
$H_{T_0T_1}(\gamma)$ of some neighborhood of $x$ in $T_0$ to some
neighborhood of $y$ in $T_1$, which is called the holonomy map
along the path $\gamma$.

One can easily see that the germ of $H_{T_0T_1}(\gamma)$ at $x$
does not depend on the choice of a partition $t_0=0<t_1<\ldots <
t_k=1$ of the segment $[0.1]$ and of a family of foliated
coordinate neighborhoods $\{U_1, U_2, \ldots, U_k\}$. Moreover,
the germ of $H_{T_0T_1}(\gamma)$ at $x$ is not changed, if we
replace $\gamma$ by any other continuous leafwise path $\gamma_1$
with the initial point $x$ and the final point $y$, which is
homotopic to $\gamma$ in the class of continuous leafwise paths
with the initial point $x$ and the final point $y$.

There is a slightly different definition of holonomy
\cite{Haefliger62}. First, let us introduce some notions.

\begin{defn}
A map $f:V\subset M\rightarrow \RR^q$ is called distinguished, if
in a neighborhood of any point in $V$ there exists a foliated
chart $(U,\phi)$ such that the restriction of $f$ to $U$ has the
form $pr_{nq}\circ \phi$, where $pr_{nq}:\RR^n =\RR^p\times\RR^q
\rightarrow \RR^q$ is the natural projection.
\end{defn}

Consider the set $D$ of germs of distinguished maps at various
points of $M$ (or, briefly speaking, the set of distinguished
germs). Let $\sigma:D\rightarrow M$ be the map, which associates
to a distinguished germ at $x\in M$ the point $x$.

Let $(U,\phi)$ be a foliated chart, $P$ its plaque. Consider the
subset $\widetilde{P}$ of $D$, which consists of all germs of the
corresponding distinguished map $pr_{nq}\circ \phi : U\rightarrow
\RR^q$ in different points of $P$. Sets of the form
$\widetilde{P}$, determined by the plaques $P$ of all possible
foliated charts, form a base of a topology on $D$, which endows
$D$ with the structure of a $p$-dimensional manifold. This
topology is called the leaf topology on $D$. It can be easily
checked \cite{Haefliger62} that $\sigma$ is a covering map from
the manifold $D$ to the manifold $\cF$.

Now consider a continuous leafwise path $\gamma$ with the initial
point $x$ and the final point $y$. Let $\pi\in \sigma^{-1}(x)$,
and let $\tilde{\gamma}$ be the lift of $\gamma$ to $D$ with the
initial point $\pi$ via the covering map $\sigma$. The holonomy
map associated with $\gamma$ is the map $$ h_\gamma
:\sigma^{-1}(x)\rightarrow\sigma^{-1}(y), $$ which takes $\pi\in
\sigma^{-1}(x)$ to the final point of the path $\tilde{\gamma}$.

The connection between the holonomy maps defined above is
established as follows. Let $T_0$ and $T_1$ be smooth transversals
such that $x\in T_0$ and $y\in T_1$. Let $\pi\in\sigma^{-1}(x)$,
and let $f : U\subset M\to \RR^q$ be a distinguished map defined
in a neighborhood of $x$. A choice of a  distinguished map $f$
defined in a neighborhood of $x$ is equivalent to a choice of a
local coordinate system $f\circ i$ on $T_0$ defined in a
neighborhood of $x$. The diffeomorphism $H_{T_0T_1}(\gamma)$
allows to define a local coordinate system on $T_1$ defined in
some neighborhood of $y$, that, in its turn, gives a distinguished
map $f_1: V\subset M\to \RR^q$ defined in a neighborhood of $y$.
The germ of the distinguished map $f_1$ at $y$ coincides with
$h_\gamma(\pi)\in \sigma^{-1}(y)$.

If $\gamma $ is a closed leafwise path with the initial and final
points $x$, and $T$ is a smooth transversal such that $x\in T$,
then $H_{TT}(\gamma)$ is a local diffeomorphism of $T$, which
leaves $x$ fixed. The correspondence $\gamma \to H_{TT}(\gamma)$
defines a group homomorphism
\[
H_T:\pi_1(L_x,x)\to \operatorname{Diff}_x(T)
\]
from the fundamental group $\pi_1(L_x,x)$ of the leaf $L_x$ to the
group $\operatorname{Diff}_x(T)$ of germs at $x$ of local
diffeomorphisms of $T$, which leave $x$ fixed. The image of the
homomorphism $H_T$ is called the holonomy group of the leaf $L_x$
at $x$. The holonomy group of a leaf $L$ at a point $x\in L$ is
independent modulo isomorphism of the choice of a transversal $T$
and $x$. A leaf is said to have trivial holonomy, if its holonomy
group is trivial.

\begin{ex}
Let $X$ be a complete nonsingular vector field on a manifold $M$
of dimension $n$, $x_0$ a (for simplicity, isolated) periodic
point of the flow $X_t$ of the given vector field and $C$ the
corresponding closed phase curve. Let $T$ be an
$(n-1)$-dimensional submanifold of $M$, passing through $x_0$
transverse to the vector $X(x_0)$: $T_xM=T_xT\oplus \RR X(x_0)$.
For all $x\in T$, closed enough to $x_0$, there is the least
$t(x)>0$ such that the corresponding positive semi-trajectory of
the flow $\{X_t(x):t>0\}$ meets $T$: $X_{t(x)}(x)\in T$. Thus, we
get a local diffeomorphism $\phi_T : x\mapsto X_{t(x)}(x)$ of $T$,
defined in a neighborhood of $x_0$ and taking $x_0$ to itself.
This diffeomorphism is called the first return map (or the
Poincar\'e map) along the curve $C$.

If $\cF$ is the foliation on $M$ given by the trajectories of $X$,
then the holonomy group of the leaf $C$ coincides with $\ZZ$, and
the germ of $\phi_T$ at $x_0$ is a generator of this group.
\end{ex}

\begin{ex}
For the Reeb foliation of the three-dimensional sphere $S^3$ all
noncompact leaves have trivial holonomy. The holonomy group of the
compact leaf is isomorphic to $\ZZ^2$.
\end{ex}

For any smooth transversal $T$ and for any $x\in T$, there is a
natural isomorphism of the tangent space $T_xT$ with the normal
space $\tau_x$ to $\cF$. Thus, the normal bundle $\tau$ plays a
role of the tangent bundle to the (germs of) transversals to
$\cF$. For any continuous leafwise path $\gamma$ with the initial
point $x$ and the final point $y$ and for any smooth transversals
$T_0$ and $T_1$ such that $x\in T_0$ and $y\in T_1$, the
differential of the holonomy map $H_{T_0T_1}(\gamma)$ at $x$
defines a linear map $dH_{T_0T_1}(\gamma)_x: \tau_x\to\tau_y$. It
is easy to check that this map is independent of the choice of
transversals $T_0$ and $T_1$. It is called the linear holonomy map
and is denoted by $dh_{\gamma}:\tau_x\to\tau_y$. Taking the
adjoint of $dh_{\gamma}$, one gets a linear map
$dh_{\gamma}^*:N^*\cF_y\to N^*\cF_x$.

Now we turn to another notion related with the holonomy, the
notion of holonomy pseudogroup. First, recall the general
definition of a pseudogroup.

\begin{defn}
A family $\Gamma$, consisting of diffeomorphisms between open
subsets of a manifold $X$ (or, in other words, of local
diffeomorphisms of $X$) is called a pseudogroup on $X$, if the
following conditions hold:
\begin{enumerate}
  \item if $\Phi\in\Gamma$, then $\Phi^{-1}\in\Gamma$;
  \item if $\Phi_1:U\to U_1$ and $\Phi_2:U_1\to U_2$ belong to $\Gamma$,
  then $\Phi_2 \circ \Phi_1:U\to U_2$ belongs to $\Gamma$;
  \item if $\Phi :U\to U_1$ belongs to $\Gamma$, then its restriction to any open subset
  $V\subset U$ belongs to $\Gamma$;
  \item if a diffeomorphism $\Phi:U\to U_1$ coincides on some neighborhood of each point in $U$
  with an element of $\Gamma$, then $\Phi\in\Gamma$;
  \item the identity diffeomorphism belongs to $\Gamma$.
\end{enumerate}
\end{defn}

\begin{ex}
The set of all local diffeomorphisms of a manifold $X$ form a
pseudogroup on $X$. One can also consider pseudogroups, consisting
of local diffeomorphisms of a manifold $X$, which preserve a
geometric structure, for instance, the pseudogroup of local
isometries and so on.
\end{ex}

\begin{defn}
Let $(M,\cF)$ be a smooth foliated manifold and $X$ the disjoint
union of all smooth transversals to $\cF$. The holonomy
pseudogroup of the foliation $\cF$ is the pseudogroup $\Gamma$,
which consists of all local diffeomorphisms of $X$, whose germ at
any point coincides with the germ of the holonomy map along a
leafwise path.
\end{defn}

\begin{defn}
Let $(M,\cF)$ be a smooth foliated manifold and $T$ be a smooth
transversal. The holonomy pseudogroup induced by the foliation
$\cF$ on $T$ is the pseudogroup $\Gamma_T$, which consists of all
local diffeomorphisms of $T$, whose germ at any point coincides
with the germ of the holonomy map along a leafwise path.
\end{defn}

There is a special class of smooth transversals given by good
covers of the manifold $M$.

\begin{defn}
A foliated chart $\phi: U\subset M \to \RR^p\times \RR^q$ is
called regular, if it admits an extension to a foliated chart
$\overline{\phi}: V\to \RR^p\times \RR^q$ such that $\overline{U}
\subset V$.
\end{defn}

\begin{defn}
A cover of a manifold $M$ by foliated neighborhoods $\{U_i\}$ is
called good, if:
\begin{enumerate}
\item Any chart $(U_i,\phi_i)$ is a regular foliated chart;
\item If $\overline{U}_i\cap \overline{U}_j\not=\varnothing$, then
${U}_i\cap {U}_j\not=\varnothing$ and the set ${U}_i\cap {U}_j$ is
connected. The same is true for the corresponding foliated
neighborhoods $V_i$;
\item Each plaque of $V_i$ meets at most one plaque of $V_j$. A plaque
of $U_i$ meets a plaque of $U_j$ if and only if the intersection
of the corresponding plaques of $V_i$ and $V_j$ is nonempty.
\end{enumerate}
\end{defn}

Good covers always exist.

Let $\cU=\{U_i\}$ be a good cover for the foliation ${\mathcal
F}$, $\phi_i:U_i\overset{\cong}{\to} I^p \times I^q$. For any $i$,
put $$ T_i=\phi^{-1}_i(\{0\}\times I^q). $$ Then $T_i$ is a
transversal, and $T=\cup T_i$ is a complete transversal. For $y\in
T_i$ denote by $P_i(y)$ the plaque of $U_i$, passing through $y$.
For any pair of indices $i$ and $j$ such that ${U}_i\cap
{U}_j\not=\varnothing$, define $$ T_{ij}=\{y\in T_i: P_i(y)\cap
U_j\not=\varnothing\}. $$ There is defined a transition function
$f_{ij}:T_{ij}\rightarrow T_{ji}$ given for $y\in T_{ij}$ by the
formula $f_{ij}(y)=y_1,$ where $y_1\in T_{ji}$ corresponds to the
unique plaque $P_j(y_1)$, for which $P_i(y)\cap
P_j(y_1)\not=\varnothing.$ The holonomy pseudogroup $\Gamma_T$
induced by $\cF$ on $T$ coincides with the pseudogroup generated
by the maps $f_{ij}$.

\begin{defn}[cf., for instance, \cite{Haefliger84}]
A transverse structure on a foliation $\cF$ is a structure on a
complete transversal $T$, invariant under the action of the
holonomy pseudogroup $\Gamma_T$.
\end{defn}

Using the notion of transverse structure, one can single out
classes of foliations with specific transverse properties. For
instance, if a complete transversal $T$ is equipped with a
Riemannian metric, and the holonomy pseudogroup $\Gamma_T$
consists of all local isometries of this Riemannian metric, we get
a class of Riemannian foliations (see Section~\ref{riemann}).
Similarly, if a complete transversal $T$ is equipped with a
symplectic structure, and the holonomy pseudogroup $\Gamma_T$
consists of all local diffeomorphisms, preserving this symplectic
structure, we get a class of symplectic foliations (see
Section~\ref{sympl}). One can also consider Kaehler foliations,
measurable foliations and so on.

\subsection{Transverse measures}
Before we turn to the discussion of an analogue of the notion of
measure on the leaf space of a foliation, we recall some basic
facts, concerning to densities and integration of densities (cf.,
for instance, \cite{Cannas-Wein,geom:asymp}).

\begin{defn}
Let $L$ be an $n$-dimensional linear space and $\cB(L)$ the set of
bases in $L$. An $\alpha$-density on $L$ ($\alpha\in \RR$) is a
function $\rho : \cB(L)\to \CC$ such that, for any $A=(A_{ij})\in
\mathop{GL}(n,\CC)$ and $e=(e_1, e_2,\ldots, e_n)\in \cB(L)$,
\[
\rho(e\cdot A)=|\det A|^\alpha \rho(e),
\]
where $(e\cdot A)_i=\sum_{j=1}^ne_jA_{ji}, i=1,2,\ldots,n$.
\end{defn}

We will denote by $|L|^\alpha$ the space of all $\alpha$-densities
on $L$. For any vector bundle $V$ on $M$, denote by $|V|^{\alpha}$
the associated bundle of $\alpha$-densities, $|V|=|V|^1$.

For any smooth, compactly supported density $\rho$ on a smooth
manifold $M$ there is a well-defined integral $\int_M\rho$,
independent of the fact if $M$ is orientable or not. This fact
allows to define a Hilbert space $L^2(M)$, canonically associated
with $M$, which consists of square integrable half-densities on
$M$. The diffeomorphism group of $M$ acts on the space $L^2(M)$ by
unitary transformations.

\begin{defn}
A (Borel) transversal to a foliation ${\mathcal F}$ is a Borel
subset of $M$, which intersects each leaf of the foliation in an
at most countable set.
\end{defn}

\begin{defn}
A transverse measure $\Lambda$ is a countably additive Radon
measure, defined on the set of all transversals to the foliation.
\end{defn}

\begin{defn}
A transverse measure $\Lambda$ is called holonomy invariant, if
for any transversals $B_1$ and $B_2$ and any bijective Borel map
$\phi :B_1\to B_2$ such that, for any $x\in B_1$, the point
$\phi(x)$ belongs to the leaf through the point $x$, we have:
$\Lambda(B_1)=\Lambda(B_2)$.
\end{defn}

\begin{ex}
Let us call by a transverse density any section of the bundle
$|\tau|$. Since, for any smooth transversal $T$, there is a
canonical isomorphism $T_xT\cong \tau_x$, a continuous positive
density $\rho\in C(M,|\tau |)$ determines a continuous positive
density on $T$, that determines a transverse measure. This
transverse measure is holonomy invariant if and only if $\rho$ is
invariant under the linear holonomy action.
\end{ex}

\begin{ex}
Any compact leaf $L$ of the foliation $\cF$ defines a holonomy
invariant transverse measure $\Lambda$. For any transversal $T$
and for any set $A\subset T$, its measure $\Lambda(A)$ equals the
number of elements in $A\cap L$.
\end{ex}

Let $\alpha\in C^\infty(M,|T\cF|)$ be a smooth positive leafwise
density on $M$. Starting from a transverse measure $\Lambda$ and
the density $\alpha$, one can construct a Borel measure $\mu$ on
$M$ in the following way. Take a good cover $\{U_i\}$ of $M$ by
foliated coordinate neighborhoods with the corresponding
coordinate maps $\phi_i: U_i\to I^p \times I^q$ and a partition of
unity $\{\psi_i\}$ subordinate to this cover. Consider the
corresponding complete transversal $T=\bigcup_i T_i$, where
$T_i=\phi^{-1}_i(\{0\}\times I^q).$ In any foliated chart
$(U_i,\phi_i)$, the transverse measure $\Lambda$ defines a measure
$\Lambda_i$ on $T_i$, and the smooth positive leafwise density
$\alpha$ defines a family $\{\alpha_{i,y} : y\in T_i\}$, where
$\{\alpha_{i,y}\}$ is a smooth positive density on the plaque
$P_i(y)$. Observe that $\Lambda$ is holonomy invariant if and only
if for any pair of indices $i$ and $j$ such that ${U}_i\cap
{U}_j\not=\varnothing$, we have: $f_{ij}(\Lambda_i)=\Lambda_j$.

For any $u\in C^\infty_c(M)$, put $$ \int_M u(m)d\mu(m)=\sum_i
\int_{T_i}\int_{P_i(y)} \psi_i(x,y)u(x,y) d\alpha_{i,y}(x)
d\Lambda_i(y). $$ One can show that this formula defines a measure
$\mu$ on $M$, which is independent of the choice of a cover
$\{U_i\}$ and a partition of unity $\{\psi_i\}$.

A measure $\mu$ on $M$ will be called holonomy invariant, if it is
obtained by means of the above construction from a holonomy
invariant transverse measure $\Lambda$ with some choice of a
smooth positive leafwise density $\alpha$.

If we take in the above construction instead of the leafwise
density $u\cdot\alpha$ the restrictions to the leaves of an
arbitrary differential $p$-form $\omega$ on $M$, we obtain a
well-defined functional $C$ on $C^\infty_c(M,\bigwedge^pT^*M)$,
called the Ruelle-Sullivan current \cite{Ruelle-S}, corresponding
to $\Lambda$:
\[
\langle C,\omega\rangle =\sum_i \int_{T_i}\int_{P_i(y)}
\psi_i(x,y)\omega_{i,y}(x) d\Lambda_i(y), \quad \omega\in
C^\infty_c(M,\bigwedge^pT^*M),
\]
where $\omega_{i,y}$ is the restriction of $\omega$ to the plaque
$P_i(y), y\in T_i$.

A transverse measure $\Lambda$ is holonomy invariant, if and only
if the corresponding Ruelle-Sullivan current $C$ is closed:
\[
\langle C,d\sigma \rangle =0, \quad \sigma \in
C^\infty_c(M,\bigwedge^{p-1}T^*M).
\]

\begin{ex}
Suppose that a transverse measure $\Lambda$ is given by a smooth
positive transverse density $\rho\in C^\infty(M,|\tau |)$. Take a
positive leafwise density $\alpha\in C^\infty(M,|T\cF|)$. Then the
corresponding measure $\mu$ on $M$ is given by the smooth positive
density  $\alpha\otimes \rho\in C^\infty(M,|TM|)$, which
corresponds to $\alpha$ and $\rho$ under the canonical isomorphism
$|TM|\cong |T\cF|\otimes |\tau|$, defined by the short exact
sequence $0\to T\cF\to TM\to\tau\to 0$.
\end{ex}

\begin{ex}
Suppose that a holonomy invariant transverse measure $\Lambda$ is
given by a compact leaf $L$ of the foliation $\cF$, and $\alpha\in
C^\infty(M,|T\cF|)$ is a smooth positive leafwise density on $M$.
Then the corresponding measure $\mu$ on $M$ is the
$\delta$-function along $L$:
\[
\int_M f(x)\,d\mu(x)=\int_Lf(x)\,\alpha(x), \quad f\in C_c(M).
\]
\end{ex}

\begin{ex}
Suppose that the foliation $\cF$ is given by the orbits of a
locally free action of a Lie group $H$ on the compact manifold $M$
and a smooth leafwise density $\alpha $ is given by a fixed Haar
measure $dh$ on $H$. Then the corresponding measure $\mu $ on $M$
is holonomy invariant if and only if it is invariant under the
action of $H$.
\end{ex}

\subsection{Connections}\label{s:conn}
The infinitesimal expression of the holonomy on a foliated
manifold is the canonical flat connection
\[
\stackrel{\circ}{\nabla} : \cX(\cF)\times C^\infty(M,\tau)\to
C^\infty(M,\tau)
\]
in the normal bundle $\tau$, defined along the leaves of $\cF$
(the Bott connection) \cite{Bott71,Bott72}. It is given by
\begin{equation}\label{e:Bott}
{\stackrel{\circ}\nabla}_X N=\theta(X)N=
P_\tau[X,\widetilde{N}],\quad X\in \cX(\cF),\quad N\in
C^\infty(M,\tau),
\end{equation}
where $\widetilde{N}\in C^\infty(M, TM)$ is any vector field on
$M$ such that $P_\tau(\widetilde{N})= N$. Thus, the restriction of
$\tau$ to any leaf of $\cF$ is a flat vector bundle. The parallel
transport in $\tau$ along any leafwise path $\gamma:x\to y$
defined by ${\stackrel{\circ}\nabla}$ coincides with the linear
holonomy map $dh_{\gamma}:\tau_x\to\tau_y$.

\begin{defn}
A connection $\nabla : \cX(M)\times C^\infty(M,\tau)\to
C^\infty(M,\tau)$ in the normal bundle $\tau$ is called adapted,
if its restriction to $\cX(\cF)$ coincides with the Bott
connection ${\stackrel{\circ}\nabla}$.
\end{defn}

One can construct an adapted connection, starting with an
arbitrary Riemannian metric $g_M$ on $M$. Denote by $\nabla^g$ the
Levi-Civita connection, defined by $g_M$. An adapted connection
$\nabla$ is given by
\begin{equation}\label{e:adapt}
\begin{aligned}
\nabla_XN &=P_\tau[X,\widetilde{N}],\quad X\in \cX(\cF),\quad N\in
C^\infty(M,\tau)\\ \nabla_XN&=P_\tau\nabla^g_X\widetilde{N},\quad
X\in C^\infty(M,F^{\bot}),\quad N\in C^\infty(M,\tau),
\end{aligned}
\end{equation}
where $\widetilde{N}\in C^\infty(M,TM)$ is any vector field such
that $P_\tau(\widetilde{N})=N$. One can show that the adapted
connection $\nabla$ described above has zero torsion.

The Bott connection ${\stackrel{\circ}\nabla}$ on $\tau$
determines a connection $({\stackrel{\circ}\nabla})^*$ on
$\tau^*\cong N^*\cF$ by the formula
\begin{equation}\label{e:dualBott}
({\stackrel{\circ}\nabla})^*_X\omega(N)=X[\omega(N)]-\omega(\nabla_XN),
\end{equation}
for any $X\in TM, \omega\in C^\infty(M,\tau^*), N\in
C^\infty(M,\tau)$.

\begin{defn}
An adapted connection $\nabla$ in the normal bundle $\tau$ is
called holonomy invariant, if, for any $X\in \cX(\cF)$, $Y\in
\cX(M)$ and $N\in C^\infty(M,\tau)$, we have
\[
(\theta(X)\nabla)_YN=0,
\]
where, by definition,
\[
(\theta(X)\nabla)_YN=\theta(X)[\nabla_YN]-\nabla_{\theta(X)Y}N-
\nabla_Y[\theta_YN].
\]

A holonomy invariant adapted connection in $\tau$ is called a
basic (or projectable) connection.
\end{defn}

A fundamental property of basic connections is the fact that their
curvature $R_\nabla$ is a basic form, i.e. $i_XR_\nabla=0,
\theta(X)R_\nabla=0$ for any $X\in \cX(\cF)$. There are
topological obstructions for the existence of basic connections
for an arbitrary foliations found independently by Kamber and
Tondeur and Molino  (cf., for instance,
\cite{KamberTond74,Molino71}).

\subsection{Riemannian foliations}\label{riemann}
\begin{defn}
A foliation $(M,\cF)$ is called Riemannian, if it has a transverse
Riemannian structure. In other words, a foliation $(M,\cF)$ is
called Riemannian, if there is a cover $\{U_i\}$ of $M$ by
foliated coordinate charts, $\phi_i:U_i\to I^p\times I^q$, and
Riemannian metrics $g^{(i)}(y)= \sum_{\alpha\beta}
g^{(i)}_{\alpha\beta}(y) dy^\alpha dy^\beta$, defined on the local
bases $I^q$ of $\cF$ such that, for any coordinate transformation
$$ \phi_{ij}(x,y)=(\alpha_{ij}(x,y), \gamma_{ij}(y)), \quad
(x,y)\in \phi_j(U_i\cap U_j), $$ the map $\gamma_{ij}$ preserves
the metric on $I^q$, $\gamma_{ij}^*(g^{(j)})=g^{(i)}$.
\end{defn}

The class of Riemannian foliations was introduced in the papers of
Reinhart \cite{Re1,Re2}. There are several equivalent
characterizations of Riemannian foliations. Before we formulate
the corresponding result (cf., for instance, \cite[Chapter
IV]{Re}), we introduce some auxiliary notions.

\begin{defn}
A distribution on a Riemannian manifold is called totally
geodesic, if every geodesic, which is tangent to the given
distribution at some point, is tangent to it along the whole its
length.
\end{defn}

\begin{defn}
The second fundamental form of a distribution $H$ on a Riemannian
manifold $(M, g)$ is the tensor $\cS$ (which takes any vector
$X\in T_xM$ at a point $x\in M$ to a linear map $\cS_X:T_xM\to
T_xM$) given by
\begin{gather*}
g(\cS_MN,X)=\frac{1}{2}g(\nabla^g_MN+\nabla^g_NM,X),\\
g(\cS_MN,L)=0,\\ g(\cS_MX,Y)=0,\\ g(\cS_MX,N)=g(\cS_MN,X),\\
\cS_X=0,
\end{gather*}
where $X,Y\in \cX(\cF)$ and $L,M,N\in C^\infty(M,F^{\bot})$,
$\nabla^g$ denote the Levi-Civita connection determined by $g$.
\end{defn}

\begin{defn}
A map $f:M\to B$ of Riemannian manifolds $M$ and $B$ is called a
Riemannian submersion, if the tangent map $df_m:T_mM\to T_{f(m)}B$
at any point $m\in M$ is surjective and induces an isometric map
between the normal space $T_mM/T_mf^{-1}(f(m))$ to the level set
$f^{-1}(f(m))$ of $f$ at $m$ and the tangent space $T_{f(m)}B$ to
$B$ at $f(m)$.
\end{defn}

\begin{theorem}\label{rev:blike}
A foliation $(M,{\mathcal F})$ is Riemannian if and only if there
is a Riemannian metric $g$ on $M$, satisfying any of the following
equivalent conditions:
\begin{enumerate}
\item The distribution $H=F^{\bot}$ is totally geodesic.
\item The second fundamental form of $H$ vanishes.
\item The induced metric $g_\tau$ on the normal bundle $\tau$ is
holonomy invariant: ${\stackrel{\circ}\nabla}_X g_\tau(M,N)=0$ for
any $X\in \cX(\cF)$ and for any $M, N\in C^\infty(M,\tau)$, where,
by definition,
\[
{\stackrel{\circ}\nabla}_X g_\tau(M,N)=X[g_\tau(M,N)]-
g_\tau({\stackrel{\circ}\nabla}_XM, N) -g_\tau(M,
{\stackrel{\circ}\nabla}_XN).
\]
\item For any vector fields $M$ and $N$, which are defined on an open set, are
orthogonal to the leaves and are infinitesimal transformations of
the foliation, and for any $X\in \cX(\cF)$, we have $ X
[g(M,N)]=0.$

\item In any foliated chart $\phi:U\to I^p\times I^q$ with the
local coordinates $(x,y)$, the restriction $g_H$ of $g$ to $H$ is
written in the form  $$ g_H=\sum_{\alpha,\beta=1}^q
g_{\alpha\beta}(y)\theta^{\alpha} \theta^{\beta}, $$ where
$\theta^{\alpha}\in H^*$ is the 1-form, corresponding to the form
$dy^\alpha$ under the isomorphism $H^*\stackrel{\cong}{\to}
T^*{\RR}^q$, and $g_{\alpha\beta}(y)$ depend only on the
transverse variables $y\in \RR^q$.

\item $(M,{\mathcal F})$ locally has the structure of a Riemannian submersion,
i. e., for any foliated chart $\phi:U\to I^p\times I^q$, there
exists a Riemannian metric on $I^q$ such that the corresponding
distinguished map $pr_{nq}\circ\phi:U\to I^q$ is a Riemannian
submersion.
\item The adapted connection $\nabla$ on the normal bundle $\tau$
given by (\ref{e:adapt}) is a (torsion-free) Riemannian
connection: for any $Y\in \cX(M)$ and $M, N\in C^\infty(M,\tau)$
\[
Y[g_\tau(M,N)]=g_\tau(\nabla_YM,N)+g_\tau(M,\nabla_YN).
\]
\item The holonomy group of the adapted connection $\nabla$ on the normal bundle $\tau$
given by (\ref{e:adapt}) at any point preserves the metric.
\end{enumerate}
\end{theorem}

\begin{defn}
Any Riemannian metric on $M$, satisfying the equivalent conditions
of Theorem~\ref{rev:blike}, is called bundle-like.
\end{defn}

One can prove that a (torsion-free) Riemannian connection on the
normal bundle $\tau$ to a Riemannian foliation $\cF$ is unique. It
is uniquely determined by the transverse metric $g_\tau$ and is
called the transverse Levi-Civita connection for $\cF$. Thus, the
transverse Levi-Civita connection is an adapted connection.
Moreover, the transverse Levi-Civita connection turns out to be a
holonomy invariant and, therefore, a basic connection. In
particular, this shows the existence of a basic connection for any
Riemannian foliation.

The existence of a bundle-like metric on a foliated manifold
imposes strong restrictions on the geometry of the foliation.
There are structure theorems for Riemannian foliations obtained by
Molino (cf. \cite{Molino,Molino82}). Using these structure
theorems, many questions, concerning to Riemannian foliations, can
be reduced to the case of Lie foliations, that is, of foliations,
whose transverse structure is modelled by a finite-dimensional Lie
group  (cf. Example~\ref{ex:lie}).

\begin{ex}
Any foliation defined by a submersion $\pi:M\to B$ is Riemannian.
\end{ex}

\begin{ex}
The orbits of a locally free isometric action of a Lie group on a
Riemannian manifold define a Riemannian foliation. On the other
hand, flows, whose orbits form a Riemannian foliation, are called
Riemannian flows. There are examples of Riemannian flows, which
are not isometric. Concerning to Riemannian flows, see, for
instance, \cite{Carriere84}, and also \cite[Appendix A]{M-S}.
\end{ex}

\begin{ex}\label{ex:lie}
Let $M$ be a smooth manifold, $\mathfrak g$ a real
finite-dimensional Lie algebra and $\omega$ an 1-form on $M$ with
values in $\mathfrak g$,  satisfying the conditions:
\begin{enumerate}
  \item the map $\omega_x:T_xM\to \mathfrak g$ is surjective for any $x\in M$;
  \item $d\omega+\frac{1}{2}[\omega,\omega]=0$.
\end{enumerate}
The distribution $F_x=\ker \omega_x$ is integrable and, therefore,
defines a codimension $q=\dim \mathfrak g$ foliation on $M$. Such
a foliation is called a Lie $\mathfrak g$-foliation. The class of
Lie foliations was introduced in \cite{Fedida}. Any Lie foliation
is Riemannian.

In the case $\mathfrak g=\RR$, a Lie foliation is precisely a
codimension one foliation given by a non-vanishing closed 1-form.
Actually, it can be easily seen that a codimension one foliation
is Riemannian if and only if it is given by a non-vanishing closed
1-form.
\end{ex}

\begin{ex}
A foliation $\cF$ on a manifold $M=\tilde{B}\times_\Gamma F$,
obtained by the suspension from a manifold $B$ and a homomorphism
$\phi : \Gamma \to \operatorname{Diff}(F)$ of the fundamental
group $\Gamma=\pi_1(B)$ is Riemannian if and only if for any
$\gamma\in\Gamma $ the diffeomorphism $\phi(\gamma)$ preserves a
Riemannian metric on $F$.
\end{ex}

Let $\cF$ be a transversely oriented Riemannian foliation and $g$
a bundle-like Riemannian metric. The induced metric on the normal
bundle $\tau$ yields the transverse volume form $v_\tau\in
C^\infty(M, \wedge^q \tau^*)=C^\infty(M, \wedge^q N^*\cF)$, which
is holonomy invariant and defines, therefore, a holonomy invariant
transverse measure on ${\mathcal F}$.

\subsection{Symplectic foliations}\label{sympl}

\begin{defn}
A foliation $(M,\cF)$ is called symplectic, if it has a transverse
symplectic structure. In other words, a foliation $(M,\cF)$ is
called symplectic, if there is a cover $\{U_i\}$ of $M$ by
foliated coordinate charts, $\phi_i:U_i\to I^p\times I^q$, and
symplectic forms $\omega_i$, defined on the local bases $I^q$ of
$\cF$ such that, for any coordinate transformation $$
\phi_{ij}(x,y)=(\alpha_{ij}(x,y), \gamma_{ij}(y)), \quad (x,y)\in
\phi_j(U_i\cap U_j), $$ the map $\gamma_{ij}$ preserves the
symplectic form, $\gamma_{ij}^*\omega_j=\omega_i$.
\end{defn}

\begin{defn}
A presymplectic manifold is a manifold equipped with a closed
$2$-form of constant rank.
\end{defn}

A transverse symplectic structure on $(M,\cF)$ uniquely determines
a presymplectic structure $\omega$ on $M$ such that $T\cF$
coincides with the kernel of $\omega$. On the other hand, if
$(M,\omega)$ is a presymplectic manifold, then the kernel of
$\omega$ determines an integrable distribution on $M$ and $\omega
$ induces a transverse symplectic structure on the corresponding
foliation $(M,\cF)$ (cf. for instance,
\cite{Bernshtein-R73,Block-Ge}, in \cite{Bernshtein-R73}
symplectic foliations are called Hamiltonian).

If $(M,\cF)$ is a symplectic foliation and $\omega$ is the
corresponding presymplectic structure on $M$, then the $q$-form
$\wedge^q\omega$ defines the holonomy invariant transverse density
$|\wedge^q\omega |\in |\wedge^q\tau^*|$, and, therefore, a
holonomy invariant transverse measure, which can be naturally
called the transverse Liouville measure.

\begin{ex}
A foliation $\cF$ on a manifold $M=\tilde{B}\times_\Gamma F$,
obtained by the suspension from a manifold $B$ and a homomorphism
$\phi : \Gamma \to \operatorname{Diff}(F)$ of the fundamental
group $\Gamma=\pi_1(B)$ is symplectic if and only if for any
$\gamma\in\Gamma $ the diffeomorphism $\phi(\gamma)$ preserves a
symplectic structure on $F$.
\end{ex}

\begin{ex}\label{ex:sympl}
Recall (cf., for instance, \cite{LM87}) that a submanifold
$\Sigma$ of a symplectic manifold $X$ is called coisotropic, if,
for any $\sigma\in \Sigma$, the skew-orthogonal complement
$(T_{\sigma}\Sigma)^{\bot}$ of $T_{\sigma}\Sigma$ is contained in
$T_{\sigma}\Sigma$. If $\Sigma$ is a coisotropic submanifold, then
the distribution $(T_{\sigma}\Sigma)^{\bot}$ is integrable, and
the corresponding foliation ${\mathcal F}_{\Sigma}$ is called the
characteristic foliation of the coisotropic submanifold $\Sigma$.
It is well-known that there is a canonical symplectic structure on
$T_{\sigma}\Sigma/(T_{\sigma}\Sigma)^{\bot}$, therefore, the
foliation ${\mathcal F}_{\Sigma}$ is symplectic. Moreover, if
${\mathcal F}_{\Sigma}$ is simple, then the set $\Gamma\subset
T^*M\times T^*M$, which consists of all $(\nu,\nu')\in T^*M\times
T^*M$ such that $\nu$ and $\nu'$ lie on the same leaf of
${\mathcal F}_{\Sigma}$ is a canonical relation.

As shown in \cite{Block-Ge} (see also \cite{Gotay}), any
presymplectic manifold can be obtained, using this construction,
that is, as a coisotropic submanifold of a symplectic manifold.

A particular example of the construction described above is the
following one. For a foliated manifold $(M,\cF)$, consider
$T^{*}M$ as a symplectic manifold with the standard symplectic
structure. Then $N^*{\mathcal F}$ is a coisotropic submanifold in
$T^*M$. The corresponding characteristic foliation ${\mathcal
F}_N$ is the natural lift of $\cF$ to the conormal bundle and is
called a horizontal (or linearized) foliation. Thus, the
linearized foliation ${\mathcal F}_N$ is symplectic. The
coordinate chart $\phi_n: N^*{\mathcal F}\to I^p \times I^q\times
{\RR}^q$ determined by a foliated coordinate chart $\phi$ on $M$
(cf. Section~\ref{s:fol}) is a foliated chart for ${\mathcal F}_N$
with plaques given by the level sets $y={\rm const}, \eta={\rm
const}$. Informally speaking, the leaf space $N^*{\mathcal
F}/{\mathcal F}_N$ of ${\mathcal F}_N$ can be considered as the
cotangent bundle to the leaf space $M/\cF$ of $\cF$.
\end{ex}

\subsection{Differential operators}\label{s:do}
Let $(M,{\mathcal F})$ be a compact foliated manifold and $E$ a
smooth vector bundle on $M$ (unless otherwise is stated, we will
assume that vector bundles under consideration are smooth and
complex). We start with some general definitions, concerning to
differential operators on $M$.

\begin{defn}
A linear differential operator $A$ of order $\mu$, acting in
$C^{\infty}(M,E)$, is called a tangential differential operator,
if, in any foliated chart $\phi : U\subset M\to I^p\times I^q$ and
any trivialization of $E$ over it, $A$ is of the form
\begin{equation}
A =\sum _{|\alpha|\leq \mu}a_{\alpha}(x,y)D^{\alpha}_{x}, \quad
(x,y)\in I^{p} \times I^{q}, \label{rev:(1.14)}
\end{equation}
where $a_{\alpha}$ are matrix-valued function on $I^p\times I^q$,
$D_x=\frac{1}{i}\frac{\partial}{\partial x}$.
\end{defn}

\begin{defn}
For a tangential differential operator $A$ given by
(\ref{rev:(1.14)}) in some foliated chart $\phi : U\subset M\to
I^p\times I^q$ and a trivialization of $E$ over it, define its
tangential (complete) symbol $$ \sigma (x,y,\xi) = \sum_{| \alpha
| \leq \mu} a_{\alpha}(x,y) \xi ^{\alpha},\quad (x,y)\in I^p\times
I^q,\quad \xi\in\RR^p, $$ and its tangential principal symbol $$
\sigma _{\mu}(x,y,\xi) = \sum _{ | \alpha | =\mu}
a_{\alpha}(x,y)\xi ^{\alpha},\quad (x,y)\in I^p\times I^q,\quad
\xi\in\RR^p. $$
\end{defn}

The tangential principal symbol is invariantly defined as a
section of the bundle $\cL(\pi_F^*E)$ on $T^*\cF$ (where $\pi_F :
T^*\cF\to M$ is the natural projection).

\begin{defn}
A tangential differential operator $A$ is called tangentially
elliptic, if its tangential principal symbol $\sigma_{\mu}$ is
invertible for $\xi \neq 0$.
\end{defn}

Let $D^m(M,E)$ denote the set of all differential operators of
order $m$ and $D^\mu({\mathcal F},E)$ denote the set of all
tangential differential operators of order $\mu$, acting in
$C^{\infty}(M,E)$.

Introduce classes $D^{m,\mu}(M,{\mathcal F},E)$, which are
linearly generated by arbitrary compositions of tangential
differential operators of order $\mu$ and differential operators
of order $m$ on $M$. In other words, an operator $A\in
D^{m,\mu}(M,{\mathcal F},E)$ is of the form $ A=\sum_{\alpha}
B_{\alpha}C_{\alpha}$, where $B_{\alpha}\in D^{m}(M,E)$,
$C_{\alpha}\in D^{\mu}({\mathcal F},E)$. If $A_1\in
D^{m_1,\mu_1}(M,{\mathcal F},E)$, $A_2\in
D^{m_2,\mu_2}(M,{\mathcal F},E)$, then $A_1\circ A_2\in
D^{m_1+m_2,\mu_1+\mu_2}(M,{\mathcal F},E)$ and, if $A\in
D^{m,\mu}(M,{\mathcal F},E)$, then $A^{*}\in D^{m,\mu}(M,{\mathcal
F},E)$. The classes $D^{m,\mu}(M,{\mathcal F},E)$ can be extended
to classes $\Psi ^{m,\mu}(M,{\mathcal F},E)$, which contain, for
instance, parametrices for elliptic operators of class
$D^{m,\mu}(M,{\mathcal F},E)$ (cf. \cite{noncom}).

We will use the standard classes of pseudodifferential operators
$\Psi^k(M,E)$ (for the theory of pseudodifferential operators see,
for instance, \cite{H3,Taylor,Treves1,Shubin:pdo}). Recall that a
pseudodifferential operator on $M$ is a linear operator $P:
C^{\infty}(M)\to {\cD}'(M)$, which can be represented in a
coordinate patch $X\subset \RR^n$ as
\[
Pu(x)=\int e^{(x-y)\xi}p(x,\xi)\,u(y)\,dy\,d\xi, \quad x\in X,
\]
where $u\in C^\infty_c(X)$, $p(x,\xi)\in S^m(X\times\RR^n)$ is the
complete symbol of $P$. Usually, we will assume that the complete
symbol $p$ can be represented as an asymptotic sum $p\sim
p_m+p_{m-1}+\ldots$, where $p_k$ is homogeneous  of degree $k$ in
$\xi$ for $|\xi |>1$. The principal symbol $p_m$ of $P$ is
well-defined as a section of the bundle $\cL(\pi^*E)$ on
$\tilde{T}^*M=T^*M\setminus \{0\}$, where $\pi : T^*M\to M$ is the
natural projection.

\begin{defn}
The transversal principal symbol $\sigma _{P}$ of an operator
$P\in \Psi ^{m}(M,E)$ is the restriction of its principal symbol
$p_m$ to $\widetilde{N}^{\ast}{\mathcal F}=N^{\ast}{\mathcal
F}\setminus \{0\}$.
\end{defn}

\begin{defn}
A operator $P\in \Psi^m(M,E)$ is called transversally elliptic, if
its transversal principal symbol $\sigma _{P}(\nu)$ is invertible
for any $\nu\in \widetilde{N}^*{\mathcal F}$.
\end{defn}

Now suppose that a closed foliated manifold $(M,{\mathcal F})$ is
equipped with a Riemannian metric $g_M$. Let $H$ be the orthogonal
complement of $F=T\cF$. Thus, there is a decomposition of $TM$
into the direct sum $TM=F\oplus H$ and the corresponding bigrading
of the exterior power bundle $\Lambda^*T^*M$: $$
\Lambda^kT^*M=\oplus_{i+j=k}\Lambda^{i,j}T^*M, \quad
\Lambda^{i,j}T^*M=\Lambda^iH^*\otimes \Lambda^jF^*. $$ There is
(cf., for instance, \cite[Proposition 10.1]{BGV},\cite{Tondeur})
the corresponding decomposition of the de Rham differential $d$
into the sum of bigraded components of the form
\begin{equation}\label{e:d}
d=d_F+d_H+\theta.
\end{equation}
Here
\begin{enumerate}
\item  $ d_F=d_{0,1}: C^{\infty}(M,\Lambda^{i,j}T^*M)\to
C^{\infty}(M,\Lambda^{i,j+1}T^*M)$ is the tangential de Rham
differential, which is a first order tangentially elliptic
operator, independent of the choice of $g$;
\item $d_H=d_{1,0}: C^{\infty}(M,\Lambda^{i,j}T^*M)\to
C^{\infty}(M,\Lambda^{i+1,j}T^*M)$ is the transversal de Rham
differential, which is a first order transversally elliptic
operator;
\item $\theta=d_{2,-1}: C^{\infty}(M,\Lambda^{i,j}T^*M)\to
C^{\infty}(M,\Lambda^{i+2,j-1}T^*M)$ is a zero order differential
operator, which is the contraction operator by the 2-form $\theta$
on $M$ with values in $F$, $\theta \in C^{\infty}(M,F\otimes
\Lambda^2\tau^*)$, given by $$ \theta(X,Y)=p_F([X,Y]),\quad X,Y\in
C^{\infty}(M,H), $$ where $P_F:TM\to F$ is the natural projection.
In particular, $\theta$ vanishes if and only if $H$ is integrable.
\end{enumerate}

There is a similar decomposition for the adjoint:
\begin{equation}\label{e:delta}
 \delta = \delta_F + \delta_H +
\theta^{*},
\end{equation} where $\delta_F$, $ \delta_H$ and $
\theta^{*}$ are the adjoints of $d_F$, $d_H$ and $\theta$ in the
Hilbert space $L^2(M,\Lambda T^{*}M)$ accordingly.

The Laplace operator $ \Delta_g=d\delta + \delta d $ of the metric
$g_M$ can be written in the form
\begin{equation}\label{e:L}
\Delta_g =\Delta_F + \Delta_H +\Delta_{-1,2}+ K_1+K_2 +K_3,
\end{equation}
where
\begin{itemize}
\item $ \Delta_F=d_F\delta_F+\delta_Fd_F\in D^{0,2}(M,{\mathcal F},
\Lambda T^{*}M )$ is the tangential Laplacian.
\item $ \Delta_H=d_H \delta_H+ \delta_Hd_H\in D^{2,0}(M,{\mathcal
F},\Lambda T^{*}M )$ is the transversal Laplacian.
\item $\Delta_{-1,2}=\theta\theta^{*}+ \theta^{*}\theta
\in D^{0,0}(M,{\mathcal F}, \Lambda T^{*}M )$.
\item $K_1 =  d_F \delta_H+ \delta_H
d_F + \delta_Fd_H+d_H \delta_F\in D^{1,0}(M,{\mathcal F}, \Lambda
T^{*}M )$.
\item $K_2 =  d_F \theta^{*}+ \theta^{*}
d_F + \delta_F\theta+\theta \delta_F\in D^{0,0}(M,{\mathcal F},
\Lambda T^{*}M).$
\item $K_3 = d_H \theta^{*}+ \theta^{*}
d_H +  \delta_H\theta+\theta
 \delta_H\in D^{1,0}(M,{\mathcal F}, \Lambda T^{*}M )$.
\end{itemize}

We also introduce the first order differential operator $D_H=d_H +
d^*_H$ in $C^{\infty}(M,\bigwedge H^{*})$, which is called the
transverse signature operator.

A basic property of geometric operators on manifolds equipped with
Riemannian foliation is that, if $\mathcal F$ is a Riemannian
foliation and $g_{M}$ is a bundle-like metric, then the operators
$d_F \delta_H+ \delta_Hd_F$ and $\delta_Fd_H+d_H\delta_F$ belong
to $D^{0,1}(M,{\mathcal F}, \Lambda T^{*}M )$. In particular,
$K_1\in D^{0,1}(M,{\mathcal F}, \Lambda T^{*}M )$.

\section{Operator algebras of foliations}
In this Section, we will describe the noncommutative algebras
associated with the leaf space of a foliation. First, we will
define an algebra, consisting of very nice functions, on which all
basic operators of analysis are defined, then, depending on a
problem in question, we will complete this algebra and obtain an
analogue of the algebra of measurable, continuous or smooth
functions. The role of a ``nice'' algebra is played by the algebra
$C^\infty_c(G)$ of smooth compactly supported functions on the
holonomy groupoid $G$ of the foliation. Therefore, we start with
the notion of holonomy groupoid of a foliation.

\subsection{Holonomy groupoid}\label{s:groupoid}
A foliation ${\cF}$ defines an equivalence relation ${\mathcal
R}\subset M\times M$ on $M$: $(x, y)\in {\mathcal R}$ if and only
if $x$ and $y$ lie on the same leaf of the foliation ${\cF}$.
Generally, ${\mathcal R}$ is not a smooth manifold, but one can
resolve its singularity, constructing a smooth manifold $G$,
called the holonomy groupoid or the graph of the foliation, which
``almost everywhere'' coincides with $\cR$ and which can be used
in many cases as a substitution for $\cR$. The idea of the
holonomy groupoid appeared in the papers of Ehresmann, Reeb and
Thom (cf. \cite{Ehresmann63,Thom64}) and was completely realized
by Winkelnkemper \cite{Winkeln}. First of all, we give the general
definition of a groupoid (see
\cite{Renault80,M-S,Cannas-Wein,Mackenzie,Paterson99} for
groupoids and related subjects).

\begin{defn}
We say that a set $G$ has the structure of a groupoid with the set
of units $G^{(0)}$, if there are defined maps
\begin{itemize}
\item $\Delta : G^{(0)}\rightarrow G$ (the diagonal map or the unit map);
\item an involution $i:G\rightarrow G$ called the inversion and
written as $i(\gamma)=\gamma^{-1}$;
\item a range map $r:G\rightarrow G^{(0)}$ and a source map $s:G\rightarrow G^{(0)}$;
\item an associative multiplication $m: (\gamma,\gamma')\rightarrow
\gamma\gamma'$ defined on the set
\[
G^{(2)}=\{(\gamma,\gamma')\in G\times G : r(\gamma')=s(\gamma)\},
\]
\end{itemize}
satisfying the conditions
\begin{itemize}
\item $r(\Delta(x))=s(\Delta(x))=x$ and $\gamma\Delta(s(\gamma))=\gamma$,
$\Delta(r(\gamma))\gamma=\gamma$;
\item $r(\gamma^{-1})=s(\gamma)$ and $\gamma\gamma^{-1}=\Delta(r(\gamma))$.
\end{itemize}
\end{defn}

Alternatively, one can define a groupoid as a small category,
where each morphism is an isomorphism.

It is convenient to think of an element $\gamma\in G$ as an arrow
$\gamma :x \to y$, going from $x=s(\gamma)$ to $y=r(\gamma)$.

We will use the standard notation (for $x,y\in G^{(0)}$):
\begin{itemize}
\item $G^x=\{\gamma\in G:r(\gamma)=x\} =r^{-1}(x)$,
\item $G_x=\{\gamma\in G:s(\gamma)=x\} =s^{-1}(x)$,
\item $G^x_y=\{\gamma\in G : s(\gamma)=x, r(\gamma)=y\}$.
\end{itemize}

\begin{defn}
A groupoid $G$ is called smooth (or a Lie groupoid), if $G^{(0)}$,
$G$ and $G^{(2)}$ are smooth manifolds, $r$, $s$, $i$ and $m$ are
smooth maps, $r$ and $s$ are submersions, and $\Delta$ is an
immersion.
\end{defn}

\begin{ex}
\emph{Lie groups.} A Lie group $H$ defines a smooth groupoid as
follows: $G=H$, $G^{(0)}$ consists of a single point, the maps $i$
and $m$ are given by the group operations in $H$.
\end{ex}

\begin{ex}
\emph{Trivial groupoid.} Let $X$ be an arbitrary set. Put $G=X$,
$G^{(0)}=X$, the maps $s$ and $r$ are the identity maps (that is,
in other words, each element $x\in G^{(0)}= X$ is identified with
an unique element $\gamma : x\to x$).
\end{ex}

\begin{ex}
\emph{Equivalence relations.} Any equivalence relation $R\subset
X\times X$ defines a groupoid, if we put $G^{(0)}=X$, $G=R$, the
maps $s:R\to X$ and $r:R\to X$ are given by $s(x,y)=y$,
$r(x,y)=x$. Thus, pairs $(x_1,y_1)$ and $(x_2,y_2)$ can be
multiplied if and only if $y_1=x_2$, and in this case
$(x_1,y_1)(x_2,y_2) = (x_1,y_2)$.

In the particular case of $R=X\times X$, we get a so called
principal or pair groupoid.
\end{ex}

\begin{ex}
\emph{Group actions.} Let a Lie group $H$ act smoothly from the
left on a smooth manifold $X$. The crossed product groupoid
$X\rtimes H$ is defined as follows: $G^{(0)}=X$, $G=X\times H$.
The maps $s: X\times H\to X$ and $r:X\times H\to X$ have the form
$s(x,h)=h^{-1}x$, $r(x,h)=x$. Thus, pairs $(x_1,h_1)$ and
$(x_2,h_2)$ can be multiplied if and only if $x_2=h^{-1}_1x_1$,
and in this case $(x_1,h_1)(x_2,h_2) = (x_1,h_1h_2)$.
\end{ex}

\begin{ex}
\emph{Fundamental groupoid} (cf. for instance, \cite{Spanier}).
Let $X$ be a topological space, $G=\Pi(X)$ the set of homotopy
classes of paths in $X$ with all possible endpoints. More
precisely, if $\gamma : [0,1] \to X$ is a path from $x=\gamma(0)$
to $y=\gamma(1)$, then we denote by $[\gamma]$ the homotopy class
of $\gamma$ with fixed $x$ and $y$. Define the groupoid $\Pi(X)$
as the set of triples $(x, [\gamma], y)$, where $x,y\in X$,
$\gamma $ is a path with the initial point $x=\gamma(0)$ and the
final point $y=\gamma(1)$, where the multiplication is given by
the product of paths. The groupoid $\Pi(X)$ is called the
fundamental groupoid of $X$.
\end{ex}

\begin{ex}
\emph{The Haefliger groupoid $\Gamma_n$}
\cite{Bott72,Haefliger72,Haefliger:CIME}. Let $M$ be a smooth
manifold. A groupoid $\Gamma_M$ consists of the germs of local
diffeomorphisms of $M$ at arbitrary points of $M$.
$(\Gamma_M)^{(0)}=M$. If $\gamma\in \Gamma_M$ is the germ at $x\in
M$ of a diffeomorphism $f$ from some neighborhood $U$ of $x$ on an
open set $f(U)$, then $s(\gamma)=x$, $r(\gamma)=f(x)$. The
multiplication in $\Gamma_M$ is given by the composition of maps.
If $M=\RR^n$, then the groupoid $\Gamma_M$ is denoted by
$\Gamma_n$.
\end{ex}

The holonomy groupoid $G=G(M,{\mathcal F})$ of a foliated manifold
$(M,{\mathcal F})$ is defined in the following way. Let $\sim_h$
be an equivalence relation on the set of continuous leafwise paths
$\gamma:[0,1]\rightarrow M$, setting $\gamma_1\sim_h \gamma_2$, if
$\gamma_1$ and $\gamma_2$ have the same initial and final points
and the same holonomy maps: $h_{\gamma_1} = h_{\gamma_2}$. The
holonomy groupoid $G$ is the set of $\sim_h$-equivalence classes
of leafwise paths. The set of units $G^{(0)}$ is a manifold $M$.
The multiplication in $G$ is given by the product of paths. The
corresponding range and source maps $s,r:G\rightarrow M$ are given
by $s(\gamma)=\gamma(0)$ and $r(\gamma)=\gamma(1)$. Finally, the
diagonal map $\Delta:M\rightarrow G$ takes any $x\in M$ to the
element in $G$ given by the constant path $\gamma(t)=x, t\in
[0,1]$. To simplify the notation, we will identify $x\in M$ with
$\Delta(x)\in G$.

For any $x\in M$ the map $s$ maps $G^x$ on the leaf $L_x$ through
$x$. The group $G^x_x$ coincides with the holonomy group of $L_x$.
The map $s:G^x\rightarrow L_x$ is the covering map associated with
the group $G^x_x$, called the holonomy covering.

One can also introduce the holonomy groupoid $G(L)$ of a leaf $L$
of $\cF$ as the set of $\sim_h$-equivalence classes of piecewise
smooth paths in $L$.

The holonomy groupoid $G$ has the structure of a smooth (in
general, non-Hausdorff and non-paracompact) manifold of dimension
$2p+q$. Recall the construction of an atlas on $G$ \cite{Co79}.

Let $\phi: U\to I^p\times I^q, \phi': U'\to I^p\times I^q$ be two
foliated charts, $\pi=pr_{nq}\circ\phi: U\to \RR^q$,
$\pi'=pr_{nq}\circ\phi': U'\to \RR^q$ the corresponding
distinguished maps. The foliated charts $\phi$, $\phi'$ are called
compatible, if, for any $m\in U$ and $m'\in U'$ with
$\pi(m)=\pi'(m')$, there is a leafwise path $\gamma$ from $m$ to
$m'$ such that the corresponding holonomy map $h_{\gamma}$ takes
the germ $\pi_m$ of $\pi$ at $m$ to the germ $\pi'_{m'}$ of $\pi'$
at $m'$.

For any pair of compatible foliated charts $\phi$ and $\phi'$
denote by $W(\phi,\phi')$ the subset in $G$, consisting of all
$\gamma\in G$ from $s(\gamma)=m=\phi^{-1}(x,y)\in U$ to
$r(\gamma)=m'={\phi'}^{-1}(x',y)\in U'$ such that the
corresponding holonomy map $h_{\gamma}$ takes the germ $\pi_m$ of
the map $\pi=pr_{nq}\circ\phi$ at $m$ to the germ $\pi'_{m'}$ of
the map $\pi'=pr_{nq}\circ\phi'$ at $m'$. There is a coordinate
map
\begin{equation}
\label{rev:wchart} \Gamma:W(\phi,\phi')\to I^p\times I^p\times
I^q,
\end{equation}
which takes each element $\gamma\in W(\phi,\phi')$ such that
$s(\gamma)=m=\varkappa^{-1}(x,y)$,
$r(\gamma)=m'={\phi'}^{-1}(x',y)$ and $h_{\gamma}\pi_m=\pi'_{m'}$
to the triple $(x,x',y)\in I^p\times I^p\times I^q$.

As shown in \cite{Co79}, the coordinate neighborhoods
$W(\phi,\phi')$ form an atlas of a $(2p+q)$-dimensional manifold
(in general, non-Hausdorff and non-paracompact) on $G$. Moreover,
the groupoid $G$ is a smooth groupoid.

Non-Hausdorffness of the holonomy groupoid is related with the
phenomenon of one-sided holonomy. The simplest example of a
foliation with the non-Hausdorff holonomy groupoid is given by the
trajectories of a nonsingular vector field on the plane, having a
one-sided limit cycle. As shown in \cite{Winkeln}, the holonomy
groupoid is Hausdorff if and only if the holonomy maps
$H_{T_0T_1}(\gamma_1)$ and $H_{T_0T_1}(\gamma_2)$ along any
leafwise paths $\gamma_1$ and $\gamma_2$ with the initial point
$x$ and the final point $y$, given by smooth transversals $T_0$
and $T_1$, passing through $x$ and $y$ accordingly, coincide, if
they coincide on some open subset $U\subset T_0$ such that $x\in
\bar{U}$. In particular, the holonomy groupoid is Hausdorff, if
the holonomy is trivial, or real analytic. Moreover, the holonomy
groupoid of a Riemannian foliation is Hausdorff. In the following,
we will always assume that $G$ is a Hausdorff manifold.

\begin{ex}
If a simple foliation $\cF$ is defined by a submersion $\pi :M\to
B$, then its holonomy groupoid $G$ consists of all $(x,y)\in
M\times M$ such that $\pi(x)=\pi(y)$, and, moreover, $G^{(0)}=M$,
the maps $s:G\to M$ and $r:G\to M$ are given by $s(x,y)=y$,
$r(x,y)=x$.
\end{ex}

\begin{ex}
If a foliation $\cF$ is given by the orbits of a free smooth
action of a connected Lie group $H$ on a manifold $M$, then its
holonomy groupoid coincides with the crossed product groupoid
$M\rtimes H$.
\end{ex}

Besides the holonomy groupoid, there are another groupoids, which
can be associated with the foliation. First of all, it is the
groupoid given by the equivalence relation on $M$, setting points
$x$ and $y$ to be equivalent, if they lie on the same leaf of the
foliation (the coarse groupoid). As noted above, this groupoid is
not smooth. One can also consider the fundamental groupoid of the
foliation $\Pi(M,\cF)$, which also consists of equivalence classes
of leafwise paths, where two leafwise paths are called equivalent,
if they are homotopic in the class of leafwise paths with fixed
endpoints. The fundamental groupoid of the foliation $\Pi(M,\cF)$
is a smooth groupoid (cf., for instance, \cite{phillips87}).

There is a foliation ${\mathcal G}$ of dimension $2p$ on the
holonomy groupoid $G$. In any coordinate chart $W(\phi,\phi')$
given by a pair of compatible foliated charts $\phi$ and $\phi'$,
the leaves of ${\mathcal G}$ are given by equations of the form
$y={\rm const}$. The leaf of $\cG$ through $\gamma\in G$ consists
of all $\gamma'\in G$ such that $r(\gamma)$ and $r(\gamma')$ lie
on the same leaf of $\cF$ and coincides with the holonomy groupoid
of this leaf. The holonomy group of a leaf of $\cG$ coincides with
the holonomy group of the corresponding leaf of $\cF$. (The last
statement corrects an erroneous one made in \cite{Winkeln}. This
fact was noted, for instance, by Molino in his review of
\cite{Winkeln} in {\em Mathematical Reviews} (see MR 85j:57043).)

The differential of the map $(r,s):G\to M\times M$ maps
isomorphically the tangent bundle $T{\mathcal G}$ to $\cG$ to the
bundle $F\boxtimes F$ on $M\times M$, therefore, there is a
canonical isomorphism $T{\mathcal G}\cong r^*F\oplus s^*F$.

A distribution $H$ on $M$ transverse to $\cF$ determines a
distribution $HG$ on $G$ transverse to $\cG$. For any $X\in H_y$,
there is a unique vector $\widehat{X}\in T_{\gamma}G$ such that $d
s(\widehat{X})= dh^{-1}_{\gamma}(X)$ and $d r(\widehat{X})=X$,
where $dh_{\gamma} : H_x\to H_y$ is the linear holonomy map
associated with $\gamma$. The space $H_\gamma G$ consists of all
vectors of the form $\widehat{X}\in T_{\gamma}G$ for different
$X\in H_y$. In any coordinate chart $W(\phi,\phi')$ on $G$, the
tangent space $T_\gamma {\mathcal G}$ to $\cG$ at some $\gamma $
with the coordinates $(x,x',y)$ consists of vectors of the form
$X\frac{\partial}{\partial x}+X'\frac{\partial}{\partial x'}$ and
the distribution $H_\gamma G$ consists of vectors
$X\frac{\partial}{\partial x}+X'\frac{\partial}{\partial
x'}+Y\frac{\partial}{\partial y}$ such that
$X\frac{\partial}{\partial x}+Y\frac{\partial}{\partial y}\in
H_{(x,y)}$ and $X'\frac{\partial}{\partial
x'}+Y\frac{\partial}{\partial y}\in H_{(x',y)}$.

Let $g_M$ be a Riemannian metric on $M$ and $H=F^{\bot}$. Then a
Riemannian metric $g_G$ on $G$ is defined as follows. All the
components in $T_{\gamma}G=F_y\oplus F_x\oplus H_\gamma G$ are
mutually orthogonal, and, by definition, $g_G$ coincides with
$g_M$ on $F_y\oplus H_\gamma G\cong F_y\oplus H_y=T_yM$ and with
$g_F$ on $F_x$.

If $\cF$ is Riemannian and $g_M$ is a bundle-like Riemannian
metric, then $g_G$ is bundle-like, and, therefore, $\cG$ is
Riemannian. Moreover, in this case the maps $s:G\to M$ and $r:G\to
M$ are Riemannian submersions and locally trivial fibrations. In
particular, the holonomy coverings $G^x$ of leaves of $\cF$ are
diffeomorphic \cite{Winkeln}.

\subsection{The $C^*$-algebra of a foliation and noncommutative topology}\label{s:calgebra}
In this Section, we will describe the construction of the
$C^*$-algebra associated with an arbitrary foliation, which is an
analogue of the algebra of continuous functions on the leaf space
of the foliation.

\begin{defn} (cf., for instance, \cite{Dixmier,Murphy,Ped,Tak})
A $C^*$-algebra is an involutive Banach algebra $A$ such that
\[
\|a^*a\|=\|a\|^2, \quad a\in A.
\]
\end{defn}

\begin{ex}
The simplest example of a $C^*$-algebra is given by the algebra
$C_0(X)$ of continuous functions on a locally compact Hausdorff
topological space $X$, vanishing at the infinity, which is endowed
with operations of the pointwise addition and the multiplication,
with the standard involution and with the uniform norm
\[
\|f\|=\sup_{x\in X}|f(x)|, \quad f\in C_0(X).
\]

The Gelfand-Naimark theorem allows to reconstruct uniquely from a
commutative $C^*$-algebra $A$ the locally compact Hausdorff
topological space $X$ such that $A\cong C_0(X)$. More precisely,
$X$ coincides with the set $\widehat{A}$ of all characters of the
algebra $A$, i. e., of all continuous homomorphisms $A\to \CC$,
endowed with the topology of pointwise convergence.
\end{ex}

The previous example permits to consider an arbitrary
$C^*$-algebra as the algebra of continuous functions on some
virtual space. By this reason, the theory of $C^*$-algebras is
often called as noncommutative topology.

\begin{ex}
The algebra $\cL(H)$ of bounded operators in a Hilbert space $H$
equipped with the involution given by taking the adjoints and with
the operator norm is a $C^*$-algebra.

By the second Gelfand-Naimark theorem, any $C^*$-algebra is
isometrically $\ast$-isomorphic to some norm closed
$\ast$-subalgebra of the algebra $\cL(H)$ for some Hilbert space
$H$.
\end{ex}

There are two ways to define the $C^*$-algebras associated with a
foliation. The first way makes use of the auxiliary choices of a
smooth Haar system, the second one requires no auxiliary choices
and uses the language of half-densities.

\subsubsection{Definitions, using a Haar system.}
In this Section we give the definition of the $C^*$-algebras
associated with an arbitrary smooth groupoid $G$. In fact, the
assumption of smoothness of a groupoid is not essential here, and
all the definitions can be generalized to the case of topological
groupoids (cf., for instance, \cite{Renault80}).

We will only consider Hausdorff groupoids. For the definition of
the $C^*$-algebra  of a foliation in the case when the holonomy
groupoid is non-Hausdorff, cf. \cite{Co:survey} (and also
\cite{crainic:cyclic,CrainicM00,CrainicM04}).

\begin{defn}
A smooth Haar system on a smooth groupoid $G$ is a family of
positive Radon measures $\{\nu^x:x\in G^{(0)}\}$ on $G$,
satisfying the following conditions:
\begin{enumerate}
  \item The support of the measure $\nu^x$ coincides with $G^x$, and
  $\nu^x$ is a smooth measure on $G^x$.
  \item The family $\{\nu^x:x\in G^{(0)}\}$ is left-invariant,
  that is, for any continuous function $f\in
  C(G^x), f\geq 0,$ and for any $\gamma\in G, s(\gamma)=x,
  r(\gamma)=y,$ we have
\[
\int_{G^y}f(\gamma_1)d\nu^y(\gamma_1)=\int_{G^x}
f(\gamma\gamma_1)\,d\nu^x(\gamma_1).
\]
  \item The family $\{\nu^x:x\in G^{(0)}\}$ is smooth, that is,
  for any $\phi\in C^\infty_c(G)$ the function
\[
x\in G^{(0)}\mapsto \int_{G^x}\phi(\gamma)d\nu^x(\gamma)
\]
is a smooth function on $G^{(0)}$.
\end{enumerate}
\end{defn}

For a compact foliated manifold $(M,\cF)$, there is a
distinguished class of smooth Haar systems $\{\nu^x:x\in
G^{(0)}\}$ on the holonomy groupoid $G$ of $\cF$ given by smooth
positive leafwise densities $\alpha \in C^\infty(M,|T{\mathcal
F}|)$. The positive Radon measure $\nu ^{x}$ on $G^x, x\in M,$ is
given as the lift of the density $\alpha $ via the holonomy
covering $s:G^x\to M$. In the following, we will assume that a
Haar system on $G$ constructed in such a way is fixed.

Let $G$ be a smooth groupoid, $G^{(0)}=M$ and $\{\nu^x:x\in M\}$ a
smooth Haar system. Introduce the structure of an involutive
algebra on $C^{\infty}_c(G)$ by
\begin{align*}
k_1\ast k_2(\gamma)&=\int_{G^x} k_1(\gamma_1)
k_2(\gamma^{-1}_1\gamma)\,d\nu^x(\gamma_1),\quad \gamma\in G^x,\\
k^*(\gamma)&=\overline{k(\gamma^{-1})}, \quad \gamma\in G.
\end{align*}
For any $x\in M$, there is a natural representation of
$C^{\infty}_c(G)$ in the Hilbert space $L^2(G^x,\nu^x)$ given, for
$k\in C^{\infty}_c(G)$ and $\zeta \in L^2(G^x,\nu^x)$, by
\[
R_x(k)\zeta(\gamma)=\int_{G^x}k(\gamma^{-1}\gamma_1)
\zeta(\gamma_1) d\nu^x(\gamma_1),\quad r(\gamma)=x.
\]
The completion of the involutive algebra $C^{\infty}_c(G)$ in the
norm
\[
\|k\|=\sup_x\|R_x(k)\|
\]
is called the reduced $C^*$-algebra of the groupoid $G$ and
denoted by $C^{\ast}_r(G)$. There is also defined the full
$C^{\ast}$-algebra of the groupoid $C^{\ast}(G)$, which is the
completion of $C^{\infty}_c(G)$ in the norm
\[
\|k\|_{\text{max}}=\sup \|\pi(k)\|,
\]
where supremum is taken over the set of all $\ast$-representations
$\pi$ of the algebra $C^{\infty}_c(G)$ in Hilbert spaces.

Any $k\in C^{\infty}_c(G)$ defines a bounded operator $R(k)$ in
$C^{\infty}_c(M)$. For any $u\in C^{\infty}_c(M)$, we have
\[
R(k)u(x)=\int_{G^x}k(\gamma)u(s(\gamma))d\nu^x(\gamma), \quad x\in
M.
\]
The correspondence $k\mapsto R(k)$ defines a representation of the
algebra $C^{\infty}_c(G)$ in $C^{\infty}_c(M)$.

\begin{ex}
\emph{Lie groups.} If a groupoid $G$ is given by a Lie group $H$,
then a smooth Haar system on $G$ is given by a left-invariant Haar
measure $dh$ on $H$, and the multiplication in $C^\infty_c(G)$ is
the convolution operation defined for any $u, v\in C^\infty_c(H)$
by
\[
(u\ast v)(g)=\int_Hu(h)v(h^{-1}g)\,dh, \quad g\in H,
\]
the involution is given by
\[
u^*(g)=\overline{u(g^{-1})}, \quad u\in C^\infty_c(H),
\]
and the operator algebras $C^*_r(G)$ and $C^*(G)$ are the group
$C^*$-algebras $C^*_r(H)$ and $C^*(H)$ (cf., for instance,
\cite{Ped}).
\end{ex}

\begin{ex}
\emph{Trivial groupoid.} Let $X$ be a smooth manifold. Put $G=X$,
$G^{(0)}=X$, the maps $s$ and $r$ are the identity maps. The
operator algebras $C^*_r(G)$ and $C^*(G)$ coincide with the
commutative $C^*$-algebra $C_0(X)$.
\end{ex}

\begin{ex}
\emph{Principal groupoid.} Let $X$ be a smooth manifold. Put
$G^{(0)}=X$, $G=X\times X$. The maps $s:X\times X\to X$ and
$r:X\times X\to X$ are given by $s(x,y)=y$, $r(x,y)=x$. If we
choose a Haar system on $G$, taking as $\nu^x$ a fixed smooth
measure $\mu$ in each $G^x=X$, then the operations in
$C^\infty_c(G)$ are given by
\begin{align*}
(k_1\ast k_2)(x,y)&=\int_Xk_1(x,z)k_2(z,y)d\mu(z), \quad (x,y)\in
X\times X, \\ k^*(x,y)&=\overline{k(y,x)}, \quad (x,y)\in X\times
X,
\end{align*}
where $k, k_1, k_2\in C^\infty_c(G)$. Thus, elements of
$C^\infty_c(G)$ can be considered as the kernels of integral
operators in $C^\infty(X)$. The representation $k\to R_x(k)$ takes
each $k\in C^\infty_c(G)\subset C^\infty(X\times X)$ to the
integral operator in $L^2(G^x,\nu^x)\cong L^2(X,\mu)$ with the
integral kernel $k$:
\[
R_x(k)u(y)=\int_Xk(y,z)u(z)d\mu(z), \quad u\in L^2(X,\mu).
\]
Finally, the operator algebras $C^*_r(G)$ and $C^*(G)$ coincide
with the algebra $C_0(X)\otimes \cK(L^2(X,\mu))$.
\end{ex}

\begin{ex}\label{ex:actions}
\emph{Group actions.} Let a Lie group $H$ act smoothly from the
left on a smooth manifold $X$. Consider the corresponding crossed
product groupoid $G=X\rtimes H$. Then $G^x=H$ for any $x\in X$ and
a smooth Haar system on $G=X\rtimes H$ is given by a
left-invariant Haar measure $dh$ on $H$. The multiplication in
$C^\infty_c(G)$ is defined for any $u, v\in C^\infty_c(X\times H)$
by
\[
(u\ast v)(x,g)=\int_Hu(x,h)v(h^{-1}x, h^{-1}g)\,dh, \quad (x,g)\in
X\times H,
\]
the involution is defined for $u\in C^\infty_c(X\times H)$ by
\[
u^*(x,g)=\overline{u(g^{-1}x, g^{-1})}, \quad \quad (x,g)\in
X\times H.
\]
The operator algebras $C^*_r(G)$ and $C^*(G)$ associated with the
crossed product groupoid $G=X\rtimes H$ coincide with the crossed
products $C_0(X)\rtimes_r H$ and $C_0(X)\rtimes H$ of the algebra
$C_0(X)$ by the group $H$ with respect to the induced action of
$H$ on $C_0(X)$ (cf., for instance, \cite{Ped}).

If the group $H$ is discrete, elements of $C^\infty_c(G)$ are
families $\{a_\gamma \in C^\infty(X) : \gamma\in H\}$ such that
$a_\gamma\not=0$ for finitely many elements $\gamma$. It is
convenient to write them as $a=\sum_{\gamma\in H} a_\gamma
U_\gamma $. The multiplication in $C^\infty_c(G)$ is written as
\[
(a_{\gamma_1} U_{\gamma_1} )(b_{\gamma_2} U_{\gamma_2}) =
(a_{\gamma_1} T_{\gamma_1}(b_{\gamma_2}))U_{\gamma_1\gamma_2},
\]
where $T_\gamma$ denotes the operator in $C_0(X)$ induced by the
action of $\gamma\in H$:
\[
T_\gamma f(x)=f(\gamma^{-1}x),\quad x\in X,\quad f\in C_0(X).
\]
The involution in the algebra $C^\infty_c(G)$ is given by
\[
(a_\gamma U_\gamma )^*=\bar{a}_\gamma U_{\gamma^{-1}}.
\]
\end{ex}

Let $G$ be the holonomy groupoid of a foliation $\cF$ on a compact
manifold $M$. Elements of the algebra $C^\infty_c(G)$ can be
considered as families of the kernels of integral operators along
the leaves of the foliation (more precisely, on the holonomy
coverings $G^x$). Namely, each $k\in C^\infty_c(G)$ corresponds to
the family $\{R_x(k):x\in M\}$, where $R_x(k)$ is the integral
operator in $L^2(G^x,\nu^x)$ given by the integral kernel
\[
K(\gamma_1,\gamma_2)=k(\gamma_1^{-1}\gamma_2), \quad \gamma_1,
\gamma_2\in G^x.
\]
The product of elements $k_1$ and $k_2$ from $C^\infty_c(G)$
corresponds to the composition of integral operators
$\{R_x(k_1)R_x(k_2):x\in M\}$. Finally, the representation $R$
means the natural action of such families of leafwise integral
operators in $L^2(M)$. In this case, the algebra $C^*_r(G)$ will
be often called the reduced $C^*$-algebra of the foliation and
denoted by $C^*(M,\cF)$.

\begin{ex}
Consider the simplest example of a foliation, given by the linear
foliation on the torus. Thus, suppose that $M=T^2=\RR^2/\ZZ^2$ is
the two-dimensional torus and a foliation $\cF_\theta$ is given by
the trajectories of the vector field $X=\frac{\partial}{\partial
x}+ \theta\frac{\partial}{\partial y},$ where $\theta\in \RR$ is a
fixed irrational number.

Since this foliation is given by the orbits of a free group action
of $\RR$ on $T^2$, its holonomy groupoid coincides with the
crossed product groupoid  $T^2\rtimes \RR$. Thus, $G=T^2\times
\RR$, $G^{(0)}=T^2$, $s(x,y,t)=(x-t,y-\theta t)$,
$r(x,y,t)=(x,y)$, $(x,y)\in T^2, t\in\RR,$ and the multiplication
is given by
\[
(x_1,y_1,t_1)(x_2,y_2,t_2)=(x_1,y_1,t_1+t_2),
\]
if $x_2=x_1-t_1, y_2=y_1-\theta t_1$.

The $C^*$-algebra $C^*_r(G)$ of the linear foliation on $T^2$
coincides with the crossed product $C(T^2)\rtimes_r \RR$.
Therefore, the product $k_1\ast k_2$ of $k_1, k_2\in
C^\infty_c(T^2\times \RR)\subset C(T^2)\rtimes_r \RR$ is given by
\begin{multline*}
(k_1\ast k_2)(x,y,t)=\int_{-\infty}^\infty
k_1(x_1,y_1,t_1)k_2(x-t_1, y-\theta t_1, t-t_1)\,dt_1, \\ (x,y)\in
T^2,\quad t\in\RR,
\end{multline*}
and, for any $k\in C^\infty_c(T^2\times \RR)$,
\[
k^*(x,y,t)=\overline{k(x-t,y-\theta t,-t)}, \quad (x,y)\in
T^2,\quad t\in\RR.
\]
For any $k\in C^\infty_c(T^2\times \RR)$ and for any $(x,y)\in
T^2$, the operator $R_{(x,y)}(k)$ in $L^2(G^{(x,y)})\cong
L^2(\RR)$ has the form: for any $u\in L^2(\RR)$
\[
R_{(x,y)}(k)u(t)=\int_{-\infty}^\infty k(x-t_1, y-\theta t_1,
t-t_1)u(t_1)\,dt_1, \quad t\in\RR.
\]
Finally, for any $k\in C^\infty_c(T^2\times \RR)$, the
corresponding operator $R(k)$ in $L^2(T^2)$ is given by
\[
R(k)u(x,y)=\int_{-\infty}^\infty k(x,y,t) u(x-t,y-\theta t)\,dt,
\quad (x,y)\in T^2.
\]

If $\theta$ is rational, then the linear foliation on $T^2$ is
given by the orbits of a free group action of $S^1$ on $T^2$, and
its holonomy groupoid coincides with the crossed product groupoid
$G=T^2\rtimes S^1$.
\end{ex}

We will also need the twisted version of the notion of the
$C^*$-algebra of a foliation $(M,\cF)$, with coefficients in a
Hermitian vector bundle $E$ on $M$. Denote by
$C^{\infty}_c(G,{\mathcal L}(E))$ the space of smooth, compactly
supported sections of the vector bundle $(s^*E)^*\otimes r^*E$ on
$G$. In other words, the value of $k\in C^{\infty}_c (G,{\mathcal
L}(E))$ at any $\gamma \in G$ is a linear map $k(\gamma) :
E_{s(\gamma)}\rightarrow E_{r(\gamma)}.$ The structure of an
involutive algebra on $C^{\infty}_c(G,{\mathcal L}(E))$ is defined
by analogous formulas
\begin{align*}
k_1\ast k_2(\gamma)&=\int_{G^x} k_1(\gamma_1)
k_2(\gamma^{-1}_1\gamma)\,d\nu^x(\gamma_1),\quad \gamma\in G,
\quad r(\gamma)=x, \\ k^*(\gamma) & = k(\gamma^{-1})^*, \quad
\gamma\in G.
\end{align*}
with the only difference that the product $k_1(\gamma_1)
k_2(\gamma^{-1}_1\gamma)$ means here the composition of the linear
maps $k_2(\gamma^{-1}_1\gamma): E_{s(\gamma)}\rightarrow
E_{s(\gamma_1)}$ and $k_1(\gamma_1) : E_{s(\gamma_1)}\rightarrow
E_{r(\gamma_1)}$ and $k(\gamma^{-1})^*$ means the adjoint of
$k(\gamma^{-1}): E_{r(\gamma)} \rightarrow E_{s(\gamma)}.$

Consider the vector bundle $\widetilde{E}=s^*(E)$ on $G$. Let
$\widetilde{E}^x$ be the restriction of the bundle $\widetilde{E}$
to $G^x$, and $L^2(G^x,\widetilde{E}^x)$ the Hilbert space of
$L^2$-sections of the bundle $\widetilde{E}^x$, determined by the
fixed Hermitian structure on $\widetilde{E}^x$ and the measure
$\nu^x$. For any $x\in M$, we introduce a representation $R_x$ of
the algebra $C^{\infty}_c(G,{\mathcal L}(E))$ in
$L^2(G^x,\widetilde{E}^x)$, given, for $k\in
C^{\infty}_c(G,{\mathcal L}(E))$, by
\begin{equation}\label{e:Rx}
R_x(k)\zeta(\gamma)=\int_{G^x}k(\gamma^{-1}\gamma_1)\zeta(\gamma_1)
d\nu^x(\gamma_1),\quad \zeta \in L^2(G^x,\widetilde{E}^x).
\end{equation}

The completion of $C^{\infty}_c(G,{\mathcal L}(E))$ in the norm
$\|k\|=\sup_x\|R_x(k)\|$ is called the reduced (twisted)
$C^*$-algebra of the foliation with coefficients in $E$ and
denoted by $C^{\ast}_r(G,E)$ or $C^{\ast}(M,{\mathcal F},E)$.
Denote by $C^{\ast}(G,E)$ the full $C^{\ast}$-algebra of the
foliation, defined in the same way as in the case of trivial
coefficients.

Any $k\in C^{\infty}_c(G,{\mathcal L}(E))$ defines an operator
$R_E(k)$ in $C^{\infty}(M,E)$. For any $u\in C^{\infty}(M,E)$, we
have
\begin{equation}
\label{rev:tang}
R_E(k)u(x)=\int_{G^x}k(\gamma)u(s(\gamma))d\nu^x(\gamma), \quad
x\in M.
\end{equation}
The correspondence $k\mapsto R_E(k)$ determines a representation
of the algebra $C^{\infty}_c(G,{\mathcal L}(E))$ in $L^2(M,E)$.

Let us give some facts, which relate the structure of the reduced
$C^*$-algebra of a foliation $C^*(M,\cF)$ with the topology of
$\cF$ (for more details cf. \cite{F-Sk,Hil-Skan83}).

\begin{thm}[\cite{F-Sk}]
Let $(M,\cF)$ be a foliated manifold.
\begin{enumerate}
  \item The $C^*$-algebra $C^*(M,\cF)$ is simple if and only if
  $\cF$ is minimal, i.e. every its leaf is dense in $M$.
  \item The $C^*$-algebra $C^*(M,\cF)$ is primitive if and only if $\cF$ is
   (topologically) transitive, i.e. it has a leaf, dense in $M$.
  \item If $\cF$ is amenable in the sense that
  $C^*(M,\cF)=C^*(G)$, the $C^*$-algebra $C^*(M,\cF)$ has a
  representation, consisting of compact operators if and only if $\cF$ has
  a compact leaf.
\end{enumerate}
\end{thm}

In the work \cite{F-Sk}, a description of the space of primitive
ideals of the $C^*$-algebra $C^*(M,\cF)$ is also given.

\subsubsection{Definition, using half-densities}
In this Section, we will give the definitions of the operator
algebras associated with a foliated manifold, which make no choice
of a Haar system. For this, we will use the language of
half-densities.

Let $(M,\cF)$ be a compact foliated manifold. Consider the vector
bundle of leafwise half-densities $|T{\mathcal F}|^{1/2}$ on $M$.
Pull back $|T{\mathcal F}|^{1/2}$ to the vector bundles
$s^*(|T{\mathcal F}|^{1/2})$ and $r^*(|T{\mathcal F}|^{1/2})$ on
the holonomy groupoid $G$, using the source map $s$ and the range
map $r$. Define a vector bundle $|T\cG|^{1/2}$ on $G$ as $$
|T\cG|^{1/2}=r^*(|T{\mathcal F}|^{1/2})\otimes s^*(|T{\mathcal
F}|^{1/2}). $$ The bundle $|T\cG|^{1/2}$ is naturally identified
with the bundle of leafwise half-densities on the foliated
manifold $(G,{\mathcal G})$.

The structure of an involutive algebra on
$C^{\infty}_c(G,|T\cG|^{1/2})$ is defined as
\begin{align*}
\sigma_1\ast \sigma_2(\gamma)&=\int_{\gamma_1\gamma_2=\gamma}
\sigma_1(\gamma_1)\sigma_2(\gamma_2),\quad \gamma\in G,\\
\sigma^*(\gamma)&=\overline{\sigma(\gamma^{-1})},\quad \gamma\in
G,
\end{align*}
where $\sigma, \sigma_1, \sigma_2\in C^{\infty}_c (G,
|T\cG|^{1/2})$. The formula for $\sigma_1\ast \sigma_2$ should be
interpreted in the following way. If we write $\gamma : x\to y,
\gamma_1: z\to y $ and $\gamma_2: x\to z$, then
\begin{align*}
\sigma_1(\gamma_1)\sigma_2(\gamma_2) \in & |T_y{\mathcal
F}|^{1/2}\otimes |T_z{\mathcal F}|^{1/2}\otimes |T_z{\mathcal
F}|^{1/2}\otimes |T_x{\mathcal F}|^{1/2} \\
 & \cong |T_y{\mathcal
F}|^{1/2}\otimes |T_z{\mathcal F}|^{1}\otimes |T_x{\mathcal
F}|^{1/2},
\end{align*}
and, integrating the $|T_z{\mathcal F}|^{1}$-component
$\sigma_1(\gamma_1)\sigma_2(\gamma_2)$ with respect to $z\in M$,
we get a well-defined section of the bundle $r^*(|T{\mathcal
F}|^{1/2})\otimes s^*(|T{\mathcal F}|^{1/2})=|T\cG|^{1/2}. $

In the case of nontrivial coefficients, taking values in an
Hermitian vector bundle $E$ on $M$, it is necessary to consider
the space $C^{\infty}_c(G,{\mathcal L}(E)\otimes |T\cG|^{1/2})$,
where the structure of an involutive algebra is defined by the
same formulas.

The formula (\ref{rev:tang}) can be rewritten in the language of
half-densities as follows. For any $\sigma\in
C^{\infty}_c(G,{\mathcal L}(E)\otimes |T\cG|^{1/2})$ and $u\in
C^{\infty}(M, E\otimes |TM|^{1/2})$, the element $R_E(\sigma)u$ of
$C^{\infty}(M, E\otimes |TM|^{1/2})$ is given by $$
R_E(\sigma)u(x)=\int_{G^x} \sigma(\gamma)s^*u(\gamma),\quad x\in
M. $$ This formula should be interpreted as follows. Using a
canonical isomorphism $|T M|^{1/2}\cong |T{\mathcal
F}|^{1/2}\otimes |TM/T{\mathcal F}|^{1/2}$ of vector bundles on
$M$ and a canonical isomorphism $s^*(|TM/T{\mathcal
F}|^{1/2})\cong r^*(|TM/T{\mathcal F}|^{1/2})$ of vector bundles
on $G$ given by the holonomy, we get
\[
s^*u\in C^{\infty}_c(G, s^*(E\otimes |T{\mathcal F}|^{1/2}\otimes
|TM/T{\mathcal F}|^{1/2})),
\]
and
\[
\sigma\cdot s^*u \in C^{\infty}_c(G, r^*E \otimes r^*(|T
M|^{1/2})\otimes s^*(|T{\mathcal F}|)).
\]
Integration of the component in $s^*(|T{\mathcal F}|)$ over $G^x$,
i.e. with a fixed $r(\gamma)=x\in M$, gives a well-defined section
$R_E(\sigma)u$ of $E\otimes |T M|^{1/2}$ on $M$.

\subsection{Von Neumann algebras and noncommutative measure theory}
The initial datum of noncommutative measure theory is a pair
$(\cM,\phi)$, consisting of a von Neumann algebra $\cM$ and a
weight $\phi$ on $\cM$.

\begin{defn}
A von Neumann algebra is an involutive subalgebra of the algebra
$\cL(H)$ of bounded operators in a Hilbert space $H$, closed in
the weak operator topology.
\end{defn}

\begin{defn}
A weight on a von Neumann algebra $\cM$ is a function $\phi$,
defined on the set $\cM_+$ of positive elements of $\cM$, with
values in $\bar{\RR}_+=[0.+\infty]$, satisfying the conditions
\begin{gather*}
\phi(a+b)=\phi(a)+\phi(b), \quad a,b \in \cM_+,\\ \phi(\alpha
a)=\alpha \phi(a), \quad \alpha \in \RR_+,\quad a\in \cM_+.
\end{gather*}

A weight on a von Neumann algebra $\cM$ is called a trace, if
\[
\phi(a^*a)=\phi(aa^*), \quad a\in \cM_+.
\]
\end{defn}

\begin{defn}
A weight $\phi$ on a von Neumann algebra $\cM$ is called
\begin{enumerate}
  \item faithful, if, for any $a\in \cM_+$, the identity $\phi(a)=0$
  implies $a=0$;
  \item normal, if, for any bounded increasing net
  $\{a_\alpha\}$ of elements from $\cM_+$ with the least upper bound
  $a$, the following identity holds:
  \[
\phi(a)=\sup_\alpha\phi(a_\alpha).
  \]
  \item semifinite, if the linear span of the set $\{x\in \cM_+:\phi(x)<\infty\}$
  is $\sigma$-weakly dense in $\cM$.
\end{enumerate}
\end{defn}

Every von Neumann algebra has a faithful, normal, semifinite
weight.

\begin{ex}
The trace functional $\tr $ on the von Neumann algebra $\cL(\cH)$
of bounded linear operators in a Hilbert space $\cH$ is a
faithful, normal, semifinite trace.
\end{ex}

\begin{ex}
If a measure space $X$ is endowed with a $\sigma$-finite measure
$\mu$, then elements of $L^\infty(X,\mu)$ considered as
multiplication operators in the Hilbert space $L^2(X,\mu)$ form a
von Neumann algebra. Moreover, the identity
\[
\phi(f)=\int_X f(x)\,d\mu(x), \quad f\in L^\infty(X,\mu),
\]
defines a faithful, normal, semifinite trace on $L^\infty(X,\mu)$.
\end{ex}

The foundations of the noncommutative integration theory for
foliations were laid by Connes in \cite{Co79}. Let $(M,\cF)$ be a
compact foliated manifold. Let $\rho$ be a strictly positive
continuous transverse density, $\alpha$ a strictly positive smooth
leafwise density, and $\nu=s^*\alpha$ the corresponding smooth
Haar system. The measures $\rho$ and $\alpha$ can be combined to
construct a measure $\mu$ on $M$. Finally, the measure $\mu$ and
the Haar system $\nu$ define a measure $m$ on $G$:
\[
\int_Gf(\gamma)\,dm(\gamma)=\int_M\left(\int_{G^x}f(\gamma)\,d\nu^x(\gamma)
\right)\,d\mu(x), \quad f\in C_c(G).
\]
In the previous Section we defined the representation $R_x$ of the
involutive algebra $C^{\infty}_c(G)$ in the Hilbert space
$L^2(G^x,\nu^x)$ for any $x\in M$. Consider a representation $R$
of the algebra $C^{\infty}_c(G)$ in $L^2(G,m)=\int^\oplus_M
L^2(G^x,\nu^x)\,d\mu(x)$ defined as the direct integral of the
representations $R_x$:
\[
R=\int^\oplus_M R_x\,d\mu(x).
\]

\begin{defn}
The von Neumann algebra $W^{\ast}(M,\mathcal F)$ of the foliation
$\mathcal F$ is defined as the closure of the image of
$C^{\infty}_{c}(G)$ under the representation $R$ in the weak
operator topology of ${\mathcal L} (L^2(G,m))$.
\end{defn}

Since the definition of the von Neumann algebra
$W^{\ast}(M,\mathcal F)$ depends only on the class of the measure
$m$ (i.e. the family of all sets of $m$-measure zero), it is easy
to see that this definition does not depend on the choice of
$\rho$ and $\alpha$.

If we assume that the union of all leaves with nontrivial holonomy
has measure zero, elements of $W^{\ast}(M,\mathcal F)$ can be
considered as measurable families $\{T_l: l\in M/\cF\}$, where
$T_l$ is a bounded operator in $L^2(l)$ for each leaf $l$ of $\cF$
(see more details in \cite{Co79,CoLNP80,Co}).

A holonomy invariant measure $\mu$ on $M$ defines a normal
semi-finite trace ${\rm tr}_{\mu}$ on the von Neumann algebra
$W^{\ast}(M,{\mathcal F})$. For any bounded measurable function
$k$ on $G$, the value ${\rm tr}_{\mu}(k)$ is finite and is given
by
\begin{equation*}
{\rm tr}_{\mu}\;(k) = \int_{M}k(x)\,d\mu(x).
\end{equation*}
The paper \cite{Co79} gives a description of weights on the von
Neumann algebra $W^{\ast}(M,{\mathcal F})$. As explained in
\cite{CoLNP80}, the construction of \cite{Co79} can be interpreted
as a correspondence between weights on $W^{\ast}(M,{\mathcal F})$
and operator-valued densities on the leaf space $M/\cF$.

\begin{ex}
Consider the linear foliation $\cF_\theta$ on the torus $T^2$,
$\theta\in \RR$ is a fixed irrational number. The Lebesgue measure
$\mu=dx\,dy$ on $T^2$ is a holonomy invariant measure. The value
of the corresponding faithful normal semi-finite trace ${\rm
tr}_{\mu}$ on $k\in C^\infty_c(T^2\times \RR)\subset
W^{\ast}(T^2,{\mathcal F}_\theta)$ is given by
\[
{\rm tr}_{\mu}\;(k)=\int_{T^2}k(x,y,0)\,dx\,dy.
\]
\end{ex}

A von Neumann algebra is called a factor if its center consists of
operators of the form $\lambda\cdot I, \lambda\in\CC$.

\begin{thm}
Let $(M,\cF)$ be a foliated manifold. A von Neumann algebra
$W^{\ast}(M,{\mathcal F})$ is a factor if and only if the
foliation is ergodic, that is, any bounded measurable function,
constant along the leaves of the foliation $\cF$, is constant on
$M$.
\end{thm}

It is known that von Neumann algebras are classified in three
classes: type I, II and III. Any von Neumann algebra $\cM$ is
canonically represented as a direct sum $\cM_I\oplus
\cM_{II}\oplus \cM_{III}$ of von Neumann algebras, where $\cM_I$,
$\cM_{II}$ and $\cM_{III}$ are von Neumann algebras of type I, II
and III accordingly.

\begin{thm}
Let $(M,\cF)$ be a foliated manifold. The von Neumann algebra
$W^{\ast}(M,{\mathcal F})$ is of:
\begin{enumerate}
  \item type I if and only if the leaf space is isomorphic to
  the standard Borel space.
  \item type II if and only if there is a holonomy invariant transverse measure and
  the algebra is not of type I.
  \item type III if and only if there is no holonomy invariant transverse measure.
\end{enumerate}
\end{thm}

\subsection{$C^*$-modules and vector bundles}
A noncommutative generalization of the notion of vector bundle (as
a topological object) is the notion of Hilbert $C^*$-module.
Hilbert $C^*$-modules are also natural generalizations of Hilbert
spaces, which arise if we replace the field of scalars $\CC$ to an
arbitrary $C^*$-algebra. The theory of Hilbert $C^*$-modules
appeared in the papers \cite{Paschke,Rieffel74} and has found many
applications in the theory of operator algebras and its
applications (see more detailed expositions of basic facts of this
theory in \cite{Lance,Manuilov-Tr,Jensen-T,Wegge-Olsen}).

\begin{defn}\cite{Paschke,Rieffel74}
Let $B$ be a $C^*$-algebra. A pre-Hilbert $B$-module is a right
$B$-module $X$ equipped with a sesquilinear map (linear in the
second argument) $\langle \cdot, \cdot\rangle_X : X\times X\to B$,
satisfying the following conditions:
\begin{enumerate}
  \item $\langle x,x\rangle_X \geq 0$ for any $x\in X$;
  \item $\langle x,x\rangle_X = 0$ if and only if $x=0$;
  \item $\langle y,x\rangle_X = \langle x,y\rangle^*_X$ for any $x,y\in X$;
  \item $\langle x,yb\rangle_X=\langle x,y\rangle_Xb$ for any $x,y\in
X, b\in B$.
\end{enumerate}
The map $\langle \cdot, \cdot\rangle_X$ is called a $B$-valued
inner product.
\end{defn}

Let $X$ be a pre-Hilbert $B$-module. It can be shown that the
formula $\|x\|_X=\|\langle x,x\rangle_X\|^{1/2}$ defines a norm on
$X$. If $X$ is complete in the norm $\| \cdot \|_X$, then $X$ is
called a Hilbert $C^*$-module. In a general case, the action of
$B$ and the inner product on $X$ are extended to its completion
$\tilde{X}$, making $\tilde{X}$ into a Hilbert $C^*$-module.

\begin{ex}
If $E$ is a Hermitian vector bundle on a compact manifold $X$,
then the space of its continuous sections $C(X,E)$ is a Hilbert
module over the algebra $C(X)$ of continuous functions on $X$. The
action of $C(X)$ on $C(X,E)$ is given by
\[
(a\cdot s)(x)=a(x)s(x),\quad x\in X,
\]
and the inner product
\[
\langle s_1, s_2\rangle(x)= \langle s_1 (x), s_2(x)\rangle_{E_x},
\quad x\in X.
\]
\end{ex}

\begin{defn}
A vector bundle $E$ on a foliated manifold $(M,\cF)$ is called
holonomy equivariant, if there is given a representation $T$ of
the holonomy groupoid $G$ of the foliation $\cF$ in the fibers of
$E$, that is, for any $\gamma\in G, \gamma:x\rightarrow y$, there
is defined a linear operator $T(\gamma):E_x\rightarrow E_y$ such
that $T(\gamma_1\gamma_2)=T(\gamma_1)T(\gamma_2)$ for any
$\gamma_1,\gamma_2 \in G$ with $r(\gamma_2)=s(\gamma_1)$.

A Hermitian (resp. Euclidean) vector bundle $E$ on a foliated
manifold $(M,\cF)$ is called holonomy equivariant, if it is a
holonomy equivariant vector bundle and the representation $T$ is
unitary (resp. orthogonal): $T(\gamma^{-1})=T(\gamma)^*$ for any
$\gamma\in G$.
\end{defn}

For any holonomy equivariant vector bundle $E\rightarrow M$, the
action of the groupoid $G$ on $E$ defines a horizontal foliation
$\cF_E$ on $E$ of the same dimension as the foliation $\cF$. The
leaf of $\cF_E$ through a point $v\in E$ consists of all points of
the form $T(\gamma)^{-1}(v)$ with $\gamma\in G$,
$r(\gamma)=\pi(v)$. Thus, any holonomy equivariant vector bundle
is foliated in the sense of \cite{KamberTondLNM}.

\begin{defn}
A vector bundle $p:P\to M$ is called foliated, if there is a
foliation $\bar{\cF}$ on $P$ of the same dimension as $\cF$ such
that its leaves are transversal to the fibers of $p$ and are
mapped by $p$ to the leaves of $\cF$.
\end{defn}

Equivalently, one can say that a foliated vector bundle is a
vector bundle $P$ on $M$, which has a flat connection in the space
$C^\infty(M,P)$ defined along the leaves of $\cF$, that is, an
operator
\[
\nabla : \cX(\cF)\times C^\infty(M,P)\to C^\infty(M,P),
\]
satisfying, for any $f \in C^\infty(M), X\in \cX(\cF), s \in
C^\infty(M,P)$, the standard conditions
\[
\nabla_{fX}=f\nabla_{X}, \quad  \nabla_X(fs)=(X f)s +f \nabla_Xs,
\]
and also the flatness condition
\[
[\nabla_X, \nabla_Y]=\nabla_{[X,Y]}, \quad X, Y\in \cX(\cF).
\]
The parallel transport along leafwise paths associated with the
connection $\nabla $ defines an action of the fundamental groupoid
of the foliation $\Pi(M,\cF)$ in the fibers of the foliated vector
bundle $P$. In general, the parallel transport may depend on the
holonomy of the corresponding path, therefore, this action does
not necessarily pull down to an action of the holonomy groupoid in
the fibers of $P$, and, therefore, a foliated vector bundle is not
necessarily holonomy equivariant.

Let $E$ be a holonomy equivariant vector bundle. The holonomy
groupoid $G_E$ of the horizontal foliation $\cF_E$ on $E$ is
described as follows (see, for instance, \cite{Hil-Skan}):
$$ G_E=r^*(E)=\{(\gamma,v)\in G\times E : r(\gamma)=\pi(v)\}, $$
$(G_E)^{(0)}=E$, the source map $s:G_{E}\rightarrow E$ is given by
$s(\gamma,v)= T(\gamma)^{-1}(v)$, the range map
$r:G_{E}\rightarrow E$ by $r(\gamma,v)=v$ and the composition has
the form $(\gamma,v)(\gamma',v')=(\gamma\gamma',v)$, where
$v'=T(\gamma)^{-1}(v)$.

\begin{ex}
The normal bundle $\tau_x=T_xM/T_x{\mathcal F}, x\in M,$ is a
holonomy equivariant vector bundle, if it is equipped with the
action of the holonomy groupoid $G$ by the linear holonomy map
$dh_{\gamma}:\tau_x\rightarrow \tau_y, \gamma:x\rightarrow y$. The
corresponding partial flat connection defined along the leaves of
$\cF$ is the Bott connection ${\stackrel{\circ}\nabla}$ (cf.
(\ref{e:Bott})). The normal bundle $\tau$ is a holonomy
equivariant Euclidean vector bundle, if $\cF$ is a Riemannian
foliation.
\end{ex}

\begin{ex}
The dual example to the previous one is given by the conormal
bundle $N^*{\mathcal F}$ equipped with the action of the holonomy
groupoid $G$ by the linear holonomy map $(dh^*_{\gamma})^{-1}:
N^*_x{\mathcal F}\rightarrow N^*_y{\mathcal F}$ for
$\gamma:x\rightarrow y$. The corresponding flat connection defined
along the leaves of $\cF$ is the connection
$({\stackrel{\circ}\nabla})^*$, dual to the Bott connection
(cf.~(\ref{e:dualBott})).

The horizontal foliation on ${N}^*{\mathcal F}$ coincides with the
linearized foliation ${\mathcal F}_N$ (see the symplectic
description of this foliation in Example~\ref{ex:sympl}). The leaf
$\widetilde{L}_{\eta}$ of ${\mathcal F}_N$ through a point
$\eta\in N^*{\mathcal F}$ consists of all points of the form
$dh_{\gamma}^{*}(\eta)$ with $\gamma\in G$ such that
$r(\gamma)=\pi(\eta)$ (cf. also an invariant definition in
\cite{Molino}).

The holonomy groupoid of the linearized foliation ${\mathcal
F}_N$, denoted by $G_{{\mathcal F}_N}$, is described as follows:
$$ G_{{\mathcal F}_N}= \{(\gamma,\eta)\in G\times N^*{\mathcal F}
: r(\gamma)=\pi(\eta)\}
$$ with the source map $s:G_{{\mathcal F}_N}\rightarrow
N^{*}{\mathcal F}, s(\gamma,\eta)=dh_{\gamma}^{*}(\eta)$, the
range map $r:G_{{\mathcal F}_N}\rightarrow N^{*}{\mathcal F},
r(\gamma,\eta)=\eta$ and the composition
$(\gamma,\eta)(\gamma',\eta') = (\gamma\gamma',\eta)$ defined in
the case when $\eta'=dh_{\gamma}^{*}(\eta)$. The projection $\pi:
N^{*}{\mathcal F}\rightarrow M$ induces a map $\pi_G:G_{{\mathcal
F}_N}\rightarrow G$ by the formula $$
\pi_G(\gamma,\eta)=\gamma,\quad (\gamma,\eta)\in G_{{\mathcal
F}_N}. $$

Let us also note that the canonical relation $\Gamma$ in $T^*M$
associated with $N^*\cF$ (cf. Example~\ref{ex:sympl}) is given by
the (one-to-one) immersion $(r,s): G_{{\mathcal F}_N}\rightarrow
T^*M\times T^*M.$
\end{ex}

For any holonomy equivariant vector bundle $E$ on a foliated
manifold $(M,\cF)$, there is defined a pre-Hilbert
$C^\infty_c(G)$-module $\mathcal E_\infty$. As a linear space,
$\mathcal E_\infty$ coincides with $C^\infty_c(G, r^*E)$. The
module structure on $\cE_\infty$ is introduced in the following
way: the action of $f\in C^\infty_c(G)$ on $s\in \cE_\infty$ is
given by
\[
(s\ast f)(\gamma)=\int_{G^y}
s(\gamma')f({\gamma'}^{-1}\gamma)d\nu^y(\gamma'), \quad \gamma\in
G^y,
\]
and the inner product on $\cE_\infty $ with values in
$C^\infty_c(G)$ is given by
\[
\langle s_1, s_2\rangle (\gamma)=\int_{G^y}\langle
s_1({\gamma'}^{-1}), s_2({\gamma'}^{-1}\gamma )
\rangle_{E_{s(\gamma')}} d\nu^y(\gamma'), \quad s_1, s_2\in
\cE_\infty.
\]
The completion of $\cE_\infty $ in the norm $\|s\|=\|R(\langle s,
s\rangle)\|^{1/2}$ defines a $C^*$-Hilbert module over $C^*_r(G)$,
that we denote by $\cE$. It is equipped with a
$C^*_r(G)$-sesquilinear form $\langle \cdot, \cdot\rangle $ with
values in $C^*_r(G)$, which is the extension by continuity of the
sesquilinear form on $\cE_\infty $. The $C^*$-Hilbert module $\cE$
can be considered as a noncommutative analogue of the algebra of
continuous sections of the bundle $E$ considered as a bundle on
the leaf space $M/\cF$.

There is also a left action of $C^\infty_c(G)$ on $\cE_\infty$
given by
\begin{equation}\label{e:lambdaf}
(f\ast s)(\gamma)=\int_{G^y}
f(\gamma')T(\gamma')[s({\gamma'}^{-1}\gamma)]d\nu^y(\gamma'),
\quad \gamma\in G^y,
\end{equation}
where $f\in C^\infty_c(G)$ and $s\in \cE_\infty$. Unlike the right
action, the left action does not extend to an action of $C^*_r(G)$
by bounded endomorphisms of the $C^*$-Hilbert module $\cE$ over
$C^*_r(G)$. As shown in \cite{Co} (see also \cite{Co86}), for any
$f\in C^\infty_c(G)$ the formula (\ref{e:lambdaf}) defines an
endomorphism $\lambda(f)$ of the $C^*$-Hilbert module $\cE$ with
the adjoint, given by
\[
(\lambda(f)^* s)(\gamma)=\int_{G^y}
f^\sharp(\gamma')T(\gamma')[s({\gamma'}^{-1}\gamma)]d\nu^y(\gamma'),
\quad \gamma\in G^y,
\]
where $f^\sharp(\gamma)=\bar{f}(\gamma^{-1})\Delta(\gamma)$ and
$\Delta(\gamma) : E_{r(\gamma)}\to E_{r(\gamma)}$ is the linear
operator given by
\[
\Delta(\gamma)=(T(\gamma)^{-1})^* T(\gamma)^{-1}.
\]
Thus, if the Hermitian structure on $E$ is not holonomy invariant,
the representation $\lambda$ is unbounded with respect to the
$C^*$-norm on $C^*_r(G)$. Nevertheless, it can be shown that the
homomorphism $\lambda$ is a densely defined closable homomorphism
of $C^*$-algebras. Hence, the domain of the closure
$\bar{\lambda}$ of the homomorphism $\lambda$ endowed with the
graph norm $\|x\|_\lambda=\|x\|+\|\lambda(x)\|$ is a Banach
algebra $B$, dense in the $C^*$-algebra $C^*_r(G)$. Thus, the
$C^*$-Hilbert module $\cE$ is a $B$-$C^*_r(G)$-bimodule.

In particular, if we take $E$ to be the holonomy equivariant
bundle $\bigwedge^jN^*\cF$, we get the $C^*$-Hilbert module of
transverse differential $j$-forms $\Omega^j$ over the
$C^*$-algebra $C^*_r(G)$ as the completion of the pre-Hilbert
$C^\infty_c(G, |T\cG|^{1/2})$-module
$\Omega^j_\infty=C^\infty_c(G, r^*\bigwedge^jN^*\cF\otimes
|T\cG|^{1/2})$. There is a product $\Omega^j_\infty \times
\Omega^k_\infty \stackrel{\wedge}{\to} \Omega^{j+k}_\infty$,
compatible with the bimodule structure:
\begin{equation}\label{e:product}
(\omega\wedge\omega_1)(\gamma)=\int_{G^y}\omega(\gamma_1)\wedge
\gamma_1\omega_1(\gamma_1^{-1}\gamma), \quad \omega\in
\Omega^j_\infty, \quad \omega_1\in \Omega^k_\infty.
\end{equation}

For any holonomy equivariant vector bundle $E$ on a foliated
manifold $(M,\cF)$, there is defined a natural representation of
$C^\infty_c(G)$ in the space $C^\infty(M,E)$ of its smooth
sections over $M$. It is obtained as the restriction of the
representation $R_E$ of the algebra $C^\infty_c(G,\cL(E))$ to
$C^\infty_c(G)$, which is embedded in $C^\infty_c(G,\cL(E))$ by
means of the homomorphism $k\in C^\infty_c(G) \mapsto k T\in
C^\infty_c(G,\cL(E)).$ If $\{\nu^x: x\in G^{(0)}\}$ is a smooth
Haar system on $G$, then this representation, which we will also
denote by $R_E$, is defined as follows. For any $u\in
C^{\infty}(M,E)$, the section $R_E(k)u\in C^{\infty}(M,E)$ is
given by
\begin{equation}
\label{rev:tang1}
R_E(k)u(x)=\int_{G^x}k(\gamma)\,T(\gamma)[u(s(\gamma))]\,d\nu^x(\gamma),
\quad x\in M.
\end{equation}

If $E$ is a holonomy equivariant Hermitian vector bundle on $M$,
then $R_E$ is a $\ast$-representation. In this case, we define a
$C^{\ast}$-algebra $C^{\ast}_{E}(G)$ as the closure of
$R_{E}(C^{\infty}_{c}(G))$ in the uniform topology of ${\mathcal
L}(L^{2}(M,E))$.

By Theorem 2.1 in \cite{F-Sk}, there is an estimate
\begin{equation}
\label{tg:inject} \|k\|_{C^{\ast}_{r}(G,E)} \leq \|R_{E}(k)\|,
\quad k \in C^{\infty}_{c}(G,{\mathcal L}(E)),
\end{equation}
therefore, the reduced $C^{\ast}$-algebra $C^{\ast}_{r}(G,E)$ is
the quotient of the algebra $C^{\ast}_{E}(G)$, and there is a
natural projection $\pi : C^{\ast}_{E}(G) \rightarrow
C^{\ast}_{r}(G,E)$.

\section{Noncommutative topology}
\subsection{$K$-theory and $K$-homology}
One of the main tools for the investigation of topological spaces
is the topological $K$-theory (see
\cite{At:Kth,Karoubi:book,Mishch:book}). Recall that the group
$K^0(X)$ associated with a compact topological space $X$ is
generated by stable equivalence classes of locally trivial
finite-dimensional complex vector bundles on $X$. Equivalently,
elements of $K^0(X)$ are described as stable isomorphism classes
of finitely generated projective modules over the algebra $C(X)$
or as equivalence classes of projections in the matrix algebra
over $C(X)$ (see below).

In this Section, we briefly recall the definition of the
$K$-theory for $C^*$-algebras, the noncommutative analogue of the
topological $K$-theory, and also the definition of the
$K$-homology groups (see, for instance,
\cite{Blackadar,Co,Higson:Roe,Karoubi:book,Rordam,Wegge-Olsen} for
further information).

Let $A$ be a unital $C^*$-algebra. Denote by $M_n(A)$ the algebra
of $n\times n$-matrices with elements from $A$. We will assume
that $M_n(A)$ is embedded into $M_{n+1}(A)$ by the map $X\to
\begin{pmatrix} X & 0\\ 0 & 0 \end{pmatrix}$. Let
$M_\infty(A)=\underset{\rightarrow}{\lim}\, M_n(A)$.

The group $K_0(A)$ is defined as the set of homotopy equivalence
classes of projections ($p^2=p=p^*$) in $M_\infty(A)$ equipped
with the direct sum operation
\[
p_1\oplus p_2=\begin{pmatrix} p_1 & 0\\ 0 & p_2 \end{pmatrix}.
\]

Denote by $GL_n(A)$ the group of invertible $n\times n$-matrices
with elements from $A$. We will assume that $GL_n(A)$ is embedded
into $GL_{n+1}(A)$ by the map $X\to
\begin{pmatrix} X & 0\\ 0 & 1 \end{pmatrix}$. Let
$GL_\infty(A)=\underset{\rightarrow}{\lim}\, GL_n(A)$. The group
$K_1(A)$ is defined as the set of homotopy equivalence classes of
unitary matrices ($u^*u=uu^*=1$) in $GL_\infty(A)$, equipped with
the direct sum operation.

If $A$ has no unit and $A^+$ is the algebra, obtained by adding
the unit to the algebra $A$, then we have the homomorphism
$i:\CC\to A^+ : \lambda\mapsto \lambda\cdot 1$, inducing a
homomorphism $i_*: K_0(\CC)\to K_0(A^+)$, and $K_0(A)$ is defined
as the kernel of this homomorphism. Moreover, by definition,
$K_1(A)=K_1(A^+)$.

For an arbitrary algebra $\cA$ over $\CC$, there are defined the
groups $K_0(A)$ and $K_1(A)$ of the algebraic $K$-theory (see, for
instance, \cite{Bass,Milnor}). The group $K_0(A)$ is defined
similarly to the group $K_0(A)$ of the topological $K$-theory,
using idempotents ($e^2=e$) in $M_\infty(A)$ instead of
projections. The group $K_1(A)$ is defined as the quotient of the
group $GL_\infty(A)$ by the commutant $[GL_\infty(A),
GL_\infty(A)]$.

The construction of the $K$-homology for noncommutative algebras
is based on the notion of Fredholm module. This notion is a
functional-analytic abstraction of the notion of elliptic
pseudodifferential operator on a compact manifold. The idea that
abstract elliptic operators can be naturally considered as
elements of the $K$-homology groups, suggested by Atiyah
\cite{Atiyah-global}, was realized by Kasparov (see
\cite{Kasparov75}).

Recall that a Hilbert space $H$ is called $\ZZ_2$-graded, if there
is given its decomposition into a direct sum of Hilbert subspaces
$H=H_0\oplus H_1$. Equivalently, a $\ZZ_2$-grading on $H$ is
defined by a self-adjoint operator $\gamma\in \cL(H)$ such that
$\gamma^2=1$. Given the decomposition $H=H^+\oplus H^-$, the
operator $\gamma$ has the matrix $\begin{pmatrix} 1 & 0\\ 0 & -1
\end{pmatrix}.$

\begin{defn}
A Fredholm module (or a $K$-cycle) over a $C^*$-algebra $A$ is a
pair $(H,F)$, where
\begin{enumerate}
  \item $H$ is a Hilbert space, equipped with a
$\ast$-representation $\rho$ of the algebra $A$;
  \item $F$ is a bounded operator in $H$ such that, for any
  $a\in A$, the operators $(F^2-1)\rho(a)$, $(F-F^*)\rho(a)$ and
  $[F,\rho(a)]$ are compact in $H$.
\end{enumerate}
A Fredholm module $(H,F)$ is called even, if the Hilbert space $H$
is endowed with a $\ZZ_2$-grading $\gamma$, the operators
$\rho(a)$ are even, $\gamma\rho(a)=\rho(a)\gamma$, and the
operator $F$ is odd, $\gamma F=-F\gamma$. In the opposite case, it
is called odd.
\end{defn}

Recall that, for any $p\geq 1$, the Schatten class
$\cL^p({\mathcal H})$ consists of all compact operators $T$ in a
Hilbert space $\cH$ such that the operator $|T|^p$ is a trace
class operator. Let $\mu_1(T)\geq \mu_2(T)\geq \ldots$ be the
characteristic numbers ($s$-numbers) of a compact operator $T$ in
$\cH$, that is, the eigenvalues of the operator $|T|=\sqrt{T^*T}$,
taken with multiplicities. Then
\[
T\in \cL^p({\mathcal H}) \Leftrightarrow \tr |T|^p =
\sum_{n=1}^\infty |\mu_n(T)|^p <\infty.
\]

\begin{defn}
A Fredholm module $(H,F)$ is called $p$-summable, if
$(F^2-1)\rho(a)$, $(F-F^*)\rho(a)$ and $[F,\rho(a)]$ belong to
$\cL^p(H)$.
\end{defn}

The homology groups $K^0(A)$ ($K^1(A)$) are defined as the sets of
homotopy equivalence classes of even (resp. odd) Fredholm modules
over $A$. The direct sum operation defines the structure of an
abelian group on $K^0(A)$ and $K^1(A)$.

\begin{ex}[\cite{Co:nc}]
Let $M$ be a compact manifold, $E$ a Hermitian vector bundle on
$M$. Then the pair $(H,F)$, where:
\begin{itemize}
  \item $H=L^2(M,E)$;
  \item $F\in \Psi^0(M,E)$ is an elliptic operator with the principal symbol
$\sigma_F$ such that $\sigma_F^2=1, \sigma^*_F=\sigma_F$ (for
instance, one can take an operator of the form $D(1+D^2)^{-1/2}$
with some self-adjoint elliptic operator $D\in\Psi^m(M,E), m>0$)
\end{itemize}
is a Fredholm module over $C^\infty(M)$.
\end{ex}

\begin{ex}[\cite{Co:nc}]
Let $(M,\cF)$ be a compact foliated manifold, $E$ a holonomy
equivariant Hermitian vector bundle on $M$. Then the pair $(H,F)$,
where:
\begin{itemize}
  \item $H=L^2(M,E)$;
  \item $F\in \Psi^0(M,E)$ is a transversally elliptic operator with the holonomy
  invariant (see Definition~\ref{d:hinv}) transversal principal
  symbol $\sigma_F$ such that $\sigma_F^2=1, \sigma^*_F=\sigma_F$.
\end{itemize}
is a Fredholm module over $C^\infty_c(G)$.
\end{ex}

A Fredholm module  $(H,F)$ defines an index map
$\operatorname{ind}:K_{*}(A) \rightarrow \ZZ$.

In the even case, with respect to the decomposition $H=H^{+}\oplus
H^-$ given by the $\ZZ_{2}$-grading of $H$, the operator $F$ takes
the form
\begin{equation}
     F= \left(
        \begin{array}{cc}
            0& F^-  \\
            F^+& 0
        \end{array}
        \right),    \quad F_{\pm} : H^{\mp} \rightarrow H^{\pm}.
\end{equation}
For any projection $e \in M_{q}(A)$, the operator $e(F^{+}\otimes
1)e$, acting from $e(H^{+}\otimes \CC^q)$ to $e(H^{-}\otimes
\CC^q)$, is Fredholm, its index depends only on the homotopy class
of $e$. Therefore, a map $\operatorname{ind}:K_{0}(A) \rightarrow
\ZZ$ is defined as
    \begin{equation}
        \operatorname{ind}[e]= \operatorname{ind} e(F^{+}\otimes 1)e.
        \label{eq:LIFNCG.index-even}
    \end{equation}

In the odd case, for a unitary matrix $U\in GL_{q}(A)$, the
operator $(P\otimes 1) U (P\otimes 1)$, where $P=\frac{1+F}2$, is
a Fredholm operator. Moreover, the index of the operator
$(P\otimes 1) U (P\otimes 1)$ depends only on the homotopy class
of $U$. Therefore, the map $\operatorname{ind}:K_{1}(A)
\rightarrow \ZZ$ is given by
\begin{equation}
     \operatorname{ind}[U]=  \operatorname{ind} (P\otimes 1) U (P\otimes 1).
        \label{eq:LIFNCG.index-odd}
\end{equation}

In both cases, both in the even one, and in the odd one, the map
$\operatorname{ind}$ depends only on the class defined by the
Fredholm module $(H,F)$ in the $K$-homology group.

Computation of the $K$-theory for the $C^*$-algebra of a foliation
is not an easy problem. In the case of the linear foliation on
$T^2$, the computation of the $K$-theory of the $C^*$-algebra of
the foliation can be reduced to the computation of the $K$-theory
of the $C^*$-algebra $A_\theta$ (see Example~\ref{ex:Atheta}) and
was given in \cite{Pimsner-V80}. For arbitrary foliations on $T^2$
and $S^3$, the computation of the $K$-theory of the $C^*$-algebra
of the foliation was done in \cite{Torpe-Reeb}.

In \cite{Baum-Connes} (see also \cite{Co-skandal}) a geometrical
construction of elements of the group $K(C^*(M,\cF))$ was
suggested. For any topological groupoid $C$, one can define the
classifying space $BC$, which is constructed, using a modification
of the classical Milnor's construction of the classifying space of
a group \cite{Haefliger72,Mi56}. In particular, the classifying
space $BG$ of the holonomy groupoid of a foliation is defined (see
\cite[Chapter 9]{Co:survey}). The normal bundle $\tau$ on $M$
defines a vector bundle on $BG$, which is also denoted by $\tau$.
The geometric group $K_{*,\tau}(BG)$ is defined as the relative
$K$-theory group
\[
K_{*,\tau}(BG)=K_*(B(\tau),S(\tau)),
\]
where $B(\tau)$ is the unit ball subbundle of the bundle $\tau$
and $S(\tau)$ is the unit sphere bundle.

In \cite{Baum-Connes} (cf. also  \cite{Co-skandal}), a map (the
topological index)
\[
\mu : K_{*,\tau}(BG) \to K(C^*(M,\cF))
\]
was constructed by means of topological constructions. The
Baum-Connes conjecture claims that $\mu$ is an isomorphism. We
refer the reader to the papers
\cite{Hector94,Macho89,Macho00,Natsume1,Natsume2,Takai87,Takai89,Takai90,Torpe-Reeb,Tu-hyper,Tu-moyen}
for various aspects of the Baum-Connes conjecture for foliations
and related computations of the $K$-theory for $C^*$-algebras of
foliations.

\subsection{Strong Morita equivalence}\label{s:morita}
The notion of isometric $\ast$-iso\-morp\-hism between
$C^*$-algebras is a natural analogue of the notion of
homeomorphism of topological spaces. However, there is a wider
equivalence relation for $C^*$-algebras, strong Morita
equivalence, introduced by Rieffel \cite{Rieffel74,Rieffel82},
which preserves many invariants of $C^*$-algebras, for instance,
the K-theory, the space of irreducible representations, the cyclic
cohomology, and coincides with the isomorphism on the class of
commutative $C^*$-algebras. The isometric $\ast$-isomorphism of
$C^*$-algebras of foliations is related with the isomorphism of
the corresponding holonomy groupoids, while their strong Morita
equivalence is related with the isomorphism of the corresponding
leaf spaces. Therefore, in many respects the notion of strong
Morita equivalence is a more adequate notion for the study of
transverse geometric structures on foliated manifolds.

\begin{defn}
Let $A$ and $B$ be $C^*$-algebras. An $A$-$B$-equivalence bimodule
is an $A$-$B$-bimodule $X$, endowed with $A$-valued and $B$-valued
inner products $\langle \cdot, \cdot\rangle_A$ and $\langle \cdot,
\cdot\rangle_B$ accordingly, such that $X$ is a right Hilbert
$B$-module and a left Hilbert $A$-module with respect to these
inner products, and, moreover,
\begin{enumerate}
  \item $\langle x,y\rangle_A z=x\langle y, z\rangle_B$  for any
  $x, y, z\in X$;
  \item The set $\langle X,X\rangle_A$ generates a dense subset in
  $A$, and the set $\langle X,X\rangle_B$ generates a dense subset in
  $B$.
\end{enumerate}

We call algebras $A$ and $B$ strongly Morita equivalent, if there
is an $A$-$B$-equivalence bimodule.
\end{defn}

It is not difficult to show that the strong Morita equivalence is
an equivalence relation.

For any linear space $L$, denote by $\tilde{L}$ the conjugate
complex linear space, which coincides with $L$ as a set and has
the same addition operation, but the multiplication by scalars is
given by the formula $\lambda\tilde{x}=(\bar{\lambda} x)^\sim $.
If $X$ is an $A$-$B$-equivalence bimodule, then $\tilde{X}$ is
endowed with the structure of a $B$-$A$-equivalence bimodule. For
instance, $b\tilde{x}a=(a^*xb^*)^\sim .$

\begin{thm}\cite{Rieffel74}\label{t:morita}
Let $X$ be an $A$-$B$-equivalence bimodule. Then the map $E \to
X\otimes_BE$ defines an equivalence of the category of Hermitian
$B$-modules and the category of Hermitian $A$-module with the
inverse, given by the map $F \to F\otimes_B\tilde{X}$.
\end{thm}

In particular, Theorem~\ref{t:morita} implies that two commutative
$C^*$-algebras are strongly Morita equivalent if and only if they
are isomorphic.

The following theorem relates the notion of strong Morita
equivalence with the notion of stable equivalence.

\begin{thm}\cite{BroGreRie}
Let $A$ and $B$ are $C^*$-algebras with countable approximate
units. Then these algebras are strongly Morita equivalent if and
only if they are stably equivalent, i.e. $A\otimes \cK \cong
A\otimes \cK$, where $\cK$ denotes the algebra of compact
operators in a separable Hilbert space.
\end{thm}

\begin{ex}
If $\cF$ is a simple foliation given by a bundle $M\to B$, then
the $C^*$-algebra $C^*(M,\cF)$ is strongly Morita equivalent to
the $C^*$-algebra $C_0(B)$.
\end{ex}

\begin{ex} \label{ex:morita}
Consider a compact foliated manifold $(M,\cF)$. As usual, let $G$
denote the holonomy groupoid of $\cF$. For any subsets $A,
B\subset M$, denote
\[
G^A_B=\{\gamma\in G : r(\gamma)\in A, s(\gamma)\in B\}.
\]
In particular,
\[
G^M_T=\{\gamma\in G : s(\gamma)\in T\}.
\]
If $T$ is a transversal, then $G^T_T$ is a submanifold and a
subgroupoid in $G$. Let $C^*_r(G^T_T)$ be the reduced
$C^*$-algebra of this groupoid. As shown in \cite{Hil-Skan83}, if
$T$ is a complete transversal, then the algebras $C^*_r(G)$ and
$C^*_r(G^T_T)$ are strongly Morita equivalent. In particular, this
implies that
\[
C^*_r(G) \cong \cK \otimes C^*_r(G^T_T).
\]
\end{ex}

\begin{ex}\label{ex:Atheta}
Consider the linear foliation $\cF_\theta$ on the two-dimensional
torus $T^2$, where $\theta\in \RR$ is a fixed irrational number.
If we choose the transversal $T$ given by the equation $y=0$, then
the leaf space of the foliation $\cF_\theta$ is identified with
the orbit space of the $\ZZ$-action on the circle $S^1=\RR/\ZZ$
generated by the rotation
\[
R_\theta(x)=x+\theta \mod 1,\quad x\in S^1.
\]
Elements of the algebra $C^\infty_c(G^T_T)$ are determined by
matrices $a(i,j)$, where the indices $(i,j)$ are arbitrary pairs
of elements $i$ and $j$ of $T$, lying on the same leaf of $\cF$,
that is, on the same orbit of the $\ZZ$-action $R_\theta$. Since
in this case the leafwise equivalence relation on the transversal
is given by a free group action, the algebra $C^*(G^T_T)=A_\theta$
coincides with the crossed product $C(S^1)\rtimes \ZZ$ of the
algebra $C(S^1)$ by the group $\ZZ$ with respect to the
$\ZZ$-action $R_\theta$ on $C(S^1)$. Therefore (cf.
Example~\ref{ex:actions}), every element of $A_\theta$ is given by
a power series
\[
a=\sum_{n\in\ZZ}a_nU^n, \quad a_n\in C(S^1),
\]
the multiplication is given by
\[
(aU^n)(bU^m)=a(b\circ R^{-1}_{n\theta})U^{n+m}
\]
and the involution by
\[
(aU^n)^*=\bar{a}U^{-n}.
\]
The algebra $C(S^1)$ is generated by the function $V$ on $S^1$
defined as
\[
V(x)=e^{2\pi i x}, \quad x\in S^1.
\]
Hence, the algebra $A_\theta$ is generated by two elements $U$ and
$V$, satisfying the relation
\[
VU=\lambda U V, \quad \lambda = e^{2\pi i \theta}.
\]
Thus, for example, a general element of $C^\infty_c(G^T_T)$ can be
represented as a power series
\[
a=\sum_{(n,m)\in\ZZ^2}a_{nm}U^nV^m,
\]
where $a_{nm}\in \cS(\ZZ^2)$ is a rapidly decreasing sequence
(that is, for any natural $k$ we have
$\sup_{(n,m)\in\ZZ^2}(|n|+|m|)^k|a_{nm}|<\infty$).

Since, in the commutative case ($\theta=0$), the above description
defines the algebra of smooth functions on the two-dimensional
torus, the algebra $A_\theta$ is called the algebra of continuous
functions on a noncommutative torus $T^2_\theta$,
$C^\infty_c(G^T_T)= C^\infty(T^2_\theta)$,
$A_\theta=C(T^2_\theta)$. The algebra $A_\theta$ was introduced in
the paper \cite{RieffelRot} (see also \cite{CoCRAS}) and has found
many applications in mathematics and physics (cf., for instance, a
survey \cite{konechny:schw}).

The $C^*$-algebra $C^*_r(G)$ of the linear foliation on $T^2$ is
strongly Morita equivalent to $A_\theta$ (cf.
Example~\ref{ex:morita}).
\end{ex}

An important property of the groupoid $G^T_T$ associated with a
complete transversal $T$ is that it is an etale groupoid (cf.
\cite{CrainicM01}).

\begin{defn}
A smooth groupoid $G$ is called etale, if its source map $s: G\to
G^{(0)}$ is a local diffeomorphism.
\end{defn}

One can introduce an equivalence relation for groupoids, similar
to the strong Morita equivalence relation for $C^*$-algebras
\cite{Haefliger84,MRW87} (cf. also \cite{CrainicM01}). Informally
speaking, two groupoids are equivalent, if they have the same
orbit spaces and, therefore, the same transverse geometry. It is
proved in \cite{MRW87} that, if groupoids $G$ and $H$ are
equivalent, then its reduced $C^*$-algebras are strongly Morita
equivalent.

By \cite{CrainicM01}, a smooth groupoid is equivalent to an etale
one if and only if all its isotropy groups $G^x_x$ are discrete.
In the latter case, such a groupoid is called a foliation
groupoid. Examples of groupoids with these property are given by
the holonomy groupoid $G(M,\cF)$ and the fundamental groupoid
$\Pi(M,\cF)$ of a foliation $(M,\cF)$. The holonomy map defines a
groupoid morphism
\[
\operatorname{hol} : \Pi(M,\cF) \to G(M,\cF),
\]
which is identical on $G^{(0)}=M$ (a morphism over $M$).

Any foliation groupoid $G$ defines a foliation $\cF$ on
$G^{(0)}=M$. $G$ is said to be an integration of $\cF$. The next
theorem provides a precise formulation of a well-known principle,
which says that the holonomy groupoid $G(M,\cF)$ and the
fundamental groupoid $\Pi(M,\cF)$ of a foliation $(M,\cF)$ are
extreme examples of groupoids, integrating $\cF$.

A groupoid $G$ is said to be $s$-connected, if its $s$-fibers
$G_x$ are connected.

\begin{thm} \cite{CrainicM01}
Let $(M,\cF)$ be a foliated manifold. For any $s$-connected
groupoid $G$, integrating $\cF$, there is a natural representation
of the holonomy morphism $\operatorname{hol}$ as the composition
of morphisms $h_G$, $\operatorname{hol}_G$ over $M$
\[
\operatorname{hol} : \Pi(M,\cF) \stackrel{h_G}{\longrightarrow} G
\stackrel{\operatorname{hol}_G}{\longrightarrow} G(M,\cF).
\]
The maps $h_G$ and $\operatorname{hol}_G$ are surjective local
diffeomorphisms. Finally, $G$ is $s$-simply connected (i.e. has
simply connected $s$-fibers $G_x$) if and only if the morphism
$h_G$ is an isomorphism.
\end{thm}

\section{Noncommutative differential topology}
\subsection{Cyclic cohomology}
In this Section, we give the definition of the cyclic cohomology,
playing the role of a noncommutative analogue of the de Rham
homology of topological spaces (concerning to cyclic cohomology,
see the books \cite{brodzki,Co,Karoubi87,Loday} and the
bibliography therein). In the commutative case, the definition of
the de Rham cohomology requires an additional structure on a
topological space in question, for instance, the structure of a
smooth manifold. In the noncommutative case, this shows up in the
fact that cyclic cocycles are usually defined not on the
$C^*$-algebra, an analogue of the algebra of continuous functions
on the corresponding geometrical object, but on some subalgebra,
consisting of ``smooth'' functions. We postpone the consideration
of a noncommutative analogue of the notion of the algebra of
smooth functions on a smooth manifold, the notion of smooth
algebra, till the next Section, and now we turn to the definition
of the cyclic cohomology for an arbitrary algebra.

Let $\cA$ be an algebra over $\CC$. Consider the complex
$(C^*(\cA,\cA^*), b)$, where:
\begin{itemize}
  \item $C^k(\cA,\cA^*)$ is the space of $(k+1)$-linear forms on $\cA$, $k\in
  \NN$;
  \item For $\psi\in C^k(\cA,\cA^*)$, its coboundary $b\psi\in
  C^{k+1}(\cA,\cA^*)$ is given by
\begin{align*}
b\psi(a^0, \cdots, a^{k+1}) = &  \sum (-1)^j \psi(a^0, \cdots, a^j
a^{j+1}, \cdots, a^{k+1})  \\ &  + (-1)^{k+1}\psi(a^{k+1}a^0,
\cdots, a^k), \quad  a^0,a^1,\ldots, a^{k+1} \in \cA.
  \end{align*}
\end{itemize}
The cohomology of this complex is called the Hochschild cohomology
of the algebra $\cA$ with coefficients in the bimodule $\cA^*$ and
are denoted by $HH(\cA)$.

Let $C^k_\lambda(\cA)$ be the subspace of $C^k(\cA,\cA^*)$, which
consists of all $\psi\in C^k(\cA,\cA^*)$, satisfying the cyclicity
condition
\begin{equation}
    \psi(a^1, \cdots, a^k,a^0) =(-1)^k \psi(a^0, a^{1},\cdots, a^k),
    \quad a^0,a^1,\ldots, a^{k} \in
    \cA.
\end{equation}
The differential $b$ takes $C^k_\lambda(\cA)$ to
$C^{k+1}_\lambda(\cA)$, and the cyclic cohomology $HC^{*}(\cA)$ of
the algebra $\cA$ is defined as the cohomology of the complex
$(C^*_\lambda(\cA),b)$.

\begin{ex}
For $k=0$, $HC^{0}(\cA)$ is the linear space of all trace
functionals on $\cA$. By this reason, cyclic $k$-cocycles on $\cA$
are called $k$-traces on $\cA$ (for $k>0$, higher traces).
\end{ex}

\begin{ex}
For $\cA=\CC$, $HC^{n}(\CC)=0$, if $n$ is odd, and
$HC^{n}(\CC)=\CC$, if $n$ is even. A nontrivial cocycle $\phi\in
C^n_\lambda(\CC)$ is given by
\[
\phi(a^0, a^{1},\cdots, a^n)=a^0a^{1}\cdots a^n, \quad a^0,
a^{1},\cdots, a^n \in \CC.
\]
\end{ex}

Equivalently, the cyclic cohomology can be described, using a
$(b,B)$-bicomplex. Define an operator $B : C^k(\cA,\cA^*)
\rightarrow C^{k-1}(\cA,\cA^*)$ as
\[
B=AB_{0},
\]
where, for any $a^0,a^1,\ldots, a^{k-1}\in \cA$,
\begin{gather*}
     A\psi (a^0, \cdots, a^{k-1})=\sum_{j=0}^{k-1} (-1)^{(k-1)j}
    \psi(a^j,a^{j+1}\ldots, a^{k-1},a^0, a^1, \ldots,
    a^{j-1}),  \\
    B_{0}\psi (a^0, \ldots, a^{k-1}) =\psi(1, a^0, \ldots, a^{k-1}) -
    (-1)^k\psi(a^0, \ldots, a^{k-1},1).
\end{gather*}
One has that $B^2=0$, $b B=-B b$.

Consider the following double complex:
\[
C^{n,m}=C^{n-m}(\cA,\cA^*),\quad n,m\in\ZZ,
\]
with the differentials $d_1:C^{n,m}\to C^{n+1,m}$ and
$d_2:C^{n,m}\to C^{n,m+1}$ given by
\[
d_1\psi=(n-m+1)b\psi,\quad d_2\psi=\frac{1}{n-m}B\psi,\quad \psi
\in C^{n,m}.
\]
For any $q\in \NN$, consider the complex $(F^qC,d)$, where
\[
(F^qC)^p=\bigoplus_{\substack{m\geq q,\\ n+m=p}}C^{n,m}, \quad
p\in\NN,\quad d=d_1+d_2.
\]
Then one has an isomorphism
\[
HC^n(\cA)\cong H^p(F^qC), \quad n=p-2q.
\]
More precisely, this isomorphism associates to any $\psi \in
HC^n(\cA)$ a cocycle $\phi\in H^p(F^qC)$ for some $p$ and $q$ with
$n=p-2q$, which has the only nonzero component
\[
\phi_{p,q}=(-1)^{[n/2]}\psi.
\]
In particular, any cocycle in the complex $F^qC$ is cohomologic to
a cocycle of the form described above.

One can also the periodic cyclic cohomology $HP^\ev (\cA)$ and
$HP^\odd (\cA)$ by taking the inductive limit of the groups
$HC^k(\cA)$, $k\geq 0$, with respect to the periodicity operator
$S : HC^k(\cA)\to HC^{k+2}(\cA)$ given as the cup product with the
generator in $HC^{2}(\CC)$. In terms of the $(b,B)$-bicomplex, the
periodic cyclic cohomology is described as the cohomology of the
complex
\[
C^{\ev} (\cA) \stackrel{b+B}{\longrightarrow} C^{\odd}(\cA)
\stackrel{b+B}{\longrightarrow} C^{\ev}(\cA),
\]
where
\[
C^{\ev/\odd} (\cA) = \bigoplus_{k \ \text{even}/\text{odd}}
C^{k}(\cA).
\]

\begin{ex}
Let $\cA$ be the algebra $C^\infty(M)$ of smooth functions on an
$n$-dimensional compact manifold $M$ and let $\cD_{k}(M)$ denote
the space of $k$-dimensional de Rham currents on $M$. Any $C\in
\cD_{k}(M)$ determines a Hochschild cochain on $C^\infty(M)$ as
    \begin{equation}
  \psi_{C}(f^{0},f^{1}, \ldots, f^k)= \langle C, f^{0}df^{1}\wedge \ldots
  \wedge df^k\rangle,
        \quad f^{0},f^{1}, \ldots, f^k \in C^\infty(M).
    \end{equation}
This cochain satisfies the condition $B\psi_{C}=k\psi_{d^t C}$,
where $d^t$ is the de Rham boundary for currents. Thus, the map
    \begin{equation}
        \cD_{\ev/\odd}(M)\ni C=(C_{k}) \longrightarrow \varphi_{C}
        =(\frac1{k!}\psi_{C_{k}}) \in C^{\ev/\odd}(C^\infty(M))
        \label{eq:LIFNCG.morphism-homology-cyclic}
    \end{equation}
induces a homomorphism from the de Rham homology groups
$H^{\ev/\odd}(M)$ to the cyclic periodic cohomology groups
$HP^{\ev/\odd}(C^{\infty}(M))$. This homomorphism is an
isomorphism, if we restrict ourselves by the cohomology of
continuous cyclic cochains~\cite{Co}.
\end{ex}

\begin{ex}\label{ex:phi}
Another important example of cyclic cocycles is given by
normalized cocycles on a discrete group $\Gamma$. Recall (cf, for
instance, \cite{Gui}) that a (homogeneous) $k$-cocycle on $\Gamma$
is a map $h:\Gamma^{k+1} \to \CC$, satisfying the identities
\begin{gather*}
h(\gamma \gamma_0, \ldots ,\gamma \gamma_k) = h(\gamma_0, \ldots
,\gamma_k),\quad \gamma, \gamma_0, \ldots , \gamma_k \in \Gamma;\\
\sum_{i=0}^{k+1} (-1)^i h(\gamma_0, \ldots,
\gamma_{i-1},\gamma_{i+1},\ldots, \gamma_{k+1}) = 0, \quad
\gamma_0, \ldots, \gamma_{k+1}\in \Gamma.
\end{gather*}
One can associate to any homogeneous $k$-cocycle $h$ a
(nonhomogeneous) $k$-cocycle $c\in Z^k(\Gamma, \mathbb C)$ by the
formula $$ c(\gamma_1, \ldots, \gamma_k) = h (e, \gamma_1,
\gamma_1 \gamma_2, \ldots, \gamma_1\ldots \gamma_k). $$ It can be
easily checked that $c$ satisfies the following condition
\begin{multline*}
c(\gamma_1, \gamma_2,\ldots, \gamma_k)
+\sum_{i=0}^{k-1}(-1)^{i+1}c(\gamma_0,\ldots, \gamma_{i-1},
\gamma_i\gamma_{i+1}, \gamma_{i+2},\ldots,\gamma_k)\\ + (-1)^{k+1}
c(\gamma_0,\gamma_1,\ldots,\gamma_{k-1})=0.
\end{multline*}

A cocycle $c\in Z^k(\Gamma, \mathbb C)$ is said to be normalized
(in the sense of Connes), if $c(\gamma_1,\gamma_2,\ldots,
\gamma_k)$ equals zero in the case when, either $\gamma_i = e$ for
some $i$, or $\gamma_1\ldots \gamma_k = e$.

The group ring $\CC\Gamma $ consists of all functions
$f:\Gamma\to\CC$ with finite support. The multiplication in
$\CC\Gamma $ is given by the convolution operation
\[
f_1\ast f_2(\gamma)=\sum_{\gamma_1\gamma_2=\gamma}f_1(\gamma_1)
f_2(\gamma_2), \quad \gamma \in \Gamma.
\]

A normalized cocycle $c\in Z^k(\Gamma, \mathbb C)$ determines a
cyclic $k$-cocycle $\tau_c$ on $\CC\Gamma $ by the formula
\[
\tau_c(f_0, \ldots, f_k) = \sum_{\gamma_0 \ldots \gamma_k
=e}f_0(\gamma_0)\ldots f_k(\gamma_k) c(\gamma_1,\ldots, \gamma_k),
\]
where $f_j \in \CC\Gamma$ for $j=0,1,\ldots k$.
\end{ex}

In the commutative case, the $K$-theory and the cohomology of a
compact topological space $M$ are related by the Chern character
\[
\operatorname{ch} : K^*(M)\to H^*(M,\RR),
\]
which provides an isomorphism $K^*(M)\otimes \QQ \to H^*(M,\QQ)$
(cf., for instance, \cite{Mishch:book}).

If $M$ is a smooth manifold, then there is an explicit
construction of the Chern character (the Chern-Weil construction,
cf., for instance, \cite{Milnor:St,Ko-No2}), that makes use of
differential forms, currents, connections and curvature. More
precisely, if $E$ is a smooth vector bundle on $M$, then the Chern
character $\operatorname{ch}(E)\in H^\ev(M,\RR)$ of the
corresponding class $[E]$ in $K^0(M)$ is represented by the closed
differential form
\[
\operatorname{ch}(E)=\Tr \exp \left(\frac{\nabla^2}{2\pi i}\right)
\]
for any connection $\nabla : C^\infty(M,E)\to C^\infty(M,E\otimes
T^*X)$ in the bundle $E$. In the odd case (see
\cite{Baum-D82,Getzler93}), if $U \in C^{\infty}(M,U_{N}(\CC))$,
then the Chern character $\operatorname{ch}(U)\in H^\odd(M,\RR)$
of the corresponding class $[U]$ in $K^1(M)$ is given by the
cohomology class of the differential form
\[
\operatorname{ch}(U)=\sum_{k=0}^{+\infty} (-1)^{k}
\frac{k!}{(2k+1)!} \Tr (U^{-1}dU)^{2k+1}.
\]

Any closed de Rham current $C \in \cD_{k}(M)$ defines, for an even
$k$, a map $\phi_C : K^0(M)\to \CC$ by the formula
\begin{equation}\label{e:even}
  \langle\phi_C, E\rangle=\langle C, \operatorname{ch}(E) \rangle,
\quad E\in K^0(M),
\end{equation}
and, for an odd $k$. a map $\phi_C : K^1(M)\to \CC$ by the formula
\begin{equation}\label{e:odd}
  \langle [\phi_{C}], [U]\rangle = \frac1{\sqrt{2i\pi}}
            \langle C,\operatorname{ch}(U)\rangle, \quad
         U \in C^{\infty}(M,U_{N}(\CC)).
\end{equation}

The noncommutative analogue of the Chern-Weil construction
consists in the construction of a pairing between $HC^*(\cA)$ and
$K_*(\cA)$ for an arbitrary algebra $\cA$, generalizing the maps
(\ref{e:even}) and (\ref{e:odd}).

The pairing between $HC^\ev(\cA)$ and $K_{0}(\cA)$ is defined as
follows. For any cocycle $\varphi =(\varphi_{2k})$ in $C^\ev(\cA)$
and for any projection $e$ in $M_{q}(\cA)$ put
\begin{equation}
\langle[\varphi],[e]\rangle = \sum_{k \geq 0}(-1)^k
\frac{(2k)!}{k!} \varphi_{2k}\# \Tr (e-\frac{1}{2},e,\cdots,e),
\end{equation}
where $\varphi_{2k}\# \Tr$ is the $(2k+1)$-linear map on
$M_{q}(\cA)=M_{q}(\CC)\otimes \cA$, given by
\begin{equation}
  \varphi_{2k}\# \Tr(\mu^{0}\otimes a^0, \cdots, \mu^{2k}\otimes a^{2k}) =
    \Tr(\mu^{0}\ldots \mu^{2k}) \varphi_{2k}(a^0, \cdots, a^{2k}),
\end{equation}
for any $\mu^j \in M_{q}(\CC)$ and $a^j \in \cA$.

The pairing between $HC^\odd(\cA)$ and $K_{1}(\cA)$ is given by
 \begin{equation}
    \langle [\varphi], [U]\rangle = \frac1{\sqrt{2i\pi}} \sum_{k \geq 0}
    (-1)^k k! \varphi_{2k+1}\# \Tr
    (U^{-1}-1,U-1,\cdots,U^{-1}-1,U-1),
\end{equation}
where $\varphi =(\varphi_{2k+1})\in C^\odd(\cA)$ and $U$ is a
unitary matrix in $GL_{q}(\cA)$.

Now consider a $p$-summable Fredholm module $(H,F)$ over an
algebra $\cA$. The index maps~(\ref{eq:LIFNCG.index-even})
and~(\ref{eq:LIFNCG.index-odd}) are computed in terms of the
pairing of $K^{*}(\cA)$ with some cyclic cohomology class
$\tau=\operatorname{ch}_\ast (H,F) \in HC^n(\cA)$, called the
Chern character of the Fredholm module $(H,F)$. First of all, any
Fredholm module can be replaced by an equivalent one, where one
has $F^2=1$. Under this condition, in the odd case, $\tau$ is
given by
\begin{equation}\label{e:odd1}
\tau(a^0,a^1,\ldots,a^n)=\lambda_n\tr (a^0[F,a^1]\ldots [F,a^n]),
\quad a^0, a^1, \ldots, a^n\in\cA,
\end{equation}
where $n$ is odd, $n>p$, $\lambda_n$ are some constants, depending
only on $n$.

In the even case, $\tau$ is given by
\begin{equation}\label{e:even1}
\tau(a^0,a^1,\ldots,a^n)=\lambda_n\tr (\gamma a^0[F,a^1]\ldots
[F,a^n]), \quad a^0, a^1, \ldots, a^n\in\cA,
\end{equation}
where $n$ is even, $n>p$, $\gamma $ is the $\ZZ_{2}$-grading in
the space $H$, $\lambda_n$ are some constants, depending only on
$n$.

\subsection{Smooth algebras}
In this Section, we give general facts about smooth subalgebras of
$C^*$-algebras, being noncommutative analogues of the algebra of
smooth functions on a smooth manifold. Suppose that $A$ is a
$C^*$-algebra and $A^+$ is the algebra obtained by adjoining the
unit to $A$. Suppose that $\cA$ is a $*$-subalgebra of the algebra
$A$ and $\cA^+$ is the algebra obtained by adjoining the unit to
$\cA$

\begin{defn}
We say that $\cA$ is a smooth subalgebra of the algebra $A$, if
the following conditions hold:

\begin{enumerate}
\item $\cA$ is a dense $*$-subalgebra of the $C^*$-algebra $A$;

\item $\cA$ is stable under the holomorphic functional calculus, that is,
for any $a\in \cA^+ $ and for any function $f$, holomorphic in a
neighborhood of the spectrum of $a$ (considered as an element of
the algebra $A^+$) $f(a) \in \cA^+$.
\end{enumerate}
\end{defn}

Suppose that ${\cA}$ is a dense $*$-subalgebra of the algebra $A$,
endowed with the structure of a Fr\'echet algebra whose topology
is finer than the topology induced by the topology of ${A}$. A
necessary and sufficient condition for ${\mathfrak A}_0$ to be a
smooth subalgebra is given by the spectral invariance (cf.
\cite[Lemma 1.2]{Schw}):

\begin{itemize}
\item $\cA^+\cap GL(A^+) = GL(\cA^+)$, where $GL(\cA^+)$ and
$GL(A^+)$ denote the group of invertibles in $\cA^+$ and $A^+$
respectively.
\end{itemize}

This fact remains true in the case when ${\cA}$ is a locally
multiplicatively convex (i.e its topology is given by a countable
family of submultiplicative seminorms) Fr\'echet algebra such that
the group $GL(\widetilde {\cA})$ of invertibles is open
\cite[Lemma 1.2]{Schw}.

One of the most important properties of smooth subalgebras
consists in the following fact --- an analogue of smoothing in the
operator $K$-theory (cf. \cite[Sect. VI.3]{Co81}, \cite{Bost}).

\begin{thm}
If $\cA$ is a smooth subalgebra in a  $C^*$-algebra $A$, then
inclusion $\cA\hookrightarrow A$ induces an isomorphism in
$K$-theory $K(\cA) \cong K(A)$.
\end{thm}

For any $p$-summable Fredholm module $(H,F)$ over an algebra
$\cA$, the algebra, which consists of all operators $T$ in the
uniform closure $\bar{\cA}$ of $\cA$ in $\cL(\cH)$, satisfying the
condition $[F,T]\in \cL^p(\cH)$, is a smooth subalgebra of the
$C^*$-algebra $\bar{\cA}$. Similarly, every spectral triple
$(\cA,\cH,D)$ --- a noncommutative analogue of the Riemannian
structure (cf. Section~\ref{s:riem}) defines a smooth subalgebra
of the $C^*$-algebra $\bar{\cA}$ (see Section~\ref{s:sdim}).
However, in many cases the construction of an appropriate smooth
algebra is nontrivial and makes use of specific properties of the
problem in question.

\subsection{Transverse differential}\label{s:transdif}
Let $(M,\cF)$ be a foliated manifold, $H$ an arbitrary
distribution on $M$, transverse to $F=T\cF$. In this Section,
following the exposition of \cite{Co,Sau}, we define the
transverse differentiation, which is a linear map
\[
d_H: \Omega_\infty^0= C^\infty_c(G,|T\cG|^{1/2}) \to
\Omega_\infty^1=C^\infty_c(G, r^*N^*\cF\otimes |T\cG|^{1/2}),
\]
satisfying the condition
\[
d_H(k_1\ast k_2) = d_Hk_1\ast k_2+ k_1\ast d_Hk_2, \quad k_1, k_2
\in C^\infty_c(G,|T\cG|^{1/2}).
\]
The differentiation $d_H$ depends on the choice of a horizontal
distribution $H$.

For any $f\in C^\infty_c(G)$, define $d_Hf\in C^\infty_c(G,
r^*N^*\cF)$ by
\[
d_Hf(X)=d f(\widehat{X}), \quad  X\in (r^*\tau)_\gamma \cong H_y,
\quad \gamma : x\to y,
\]
where $\widehat{X}\in H_\gamma G \subset T_{\gamma}G$ is a unique
vector such that $d s(\widehat{X})= dh^{-1}_{\gamma}(X)$ and $d
r(\widehat{X})=X$.

For an arbitrary smooth leafwise density $\lambda \in
C^\infty_c(M,|T\cF|)$ define a 1-form $k(\lambda)\in
C^\infty(M,H^*)\cong C^\infty(M, N^*\cF)$ as follows. Take an
arbitrary point $m \in M$. In a foliated chart $\phi : U \to
I^p\times I^q$ defined in a neighborhood of $m$
($\phi(m)=(x^0,y^0)$), $\lambda $ can be written as
$
\lambda=f(x,y)|dx|, (x,y)\in I^p\times I^q.
$
Then
\[
k(\lambda)=d_Hf+\sum_{i=1}^p{\cL}_{\frac{\partial}{\partial
x_i}}d_Hx_i,
\]
where, for any $X=\sum_{i=1}^pX^i\frac{\partial}{\partial x_i} \in
\cX(\cF)$ and for any $\omega =\sum_{j=1}^q \omega_j\, dy_j \in
C^\infty(M,H^*)$, the Lie derivative ${\cL}_X\omega \in
C^\infty(M,H^*)$ is given by
\[
{\cL}_X\omega=\sum_{i=1}^p\sum_{j=1}^q X^i \frac{\partial \omega_j
}{\partial x_i}\,dy_j.
\]
One can give a slightly different description of $k(\lambda)$. For
any $X\in H_{m}$, let $\tilde{X}$ be an arbitrary projectable
vector field, that coincides with $X$ at $m$:
\[
\tilde{X}(x,y)=\sum_{i=1}^pX^i(x,y)\frac{\partial}{\partial
x_i}+\sum_{j=1}^qY^j(y)\frac{\partial}{\partial y_j}.
\]
Put
\begin{multline*}
k(\lambda)(X)=\sum_{i=1}^pX^i(x^0,y^0)\frac{\partial f}{\partial
x_i}(x^0,y^0)+\sum_{j=1}^qY^j(y^0)\frac{\partial f}{\partial
x_j}(x^0,y^0)\\ +\sum_{i=1}^p \frac{\partial X^i}{\partial
x_i}(x^0,y^0) f(x^0,y^0).
\end{multline*}
It can be checked that this definition is independent of the
choice of a foliated chart $\phi : U \to I^p\times I^q$ and an
extension $\tilde{X}$.

If $M$ is Riemannian, $\lambda $ is given by the induced leafwise
Riemannian volume form and $H=F^\bot$, then $k(\lambda)$ coincides
with the mean curvature 1-form of $\cF$ (cf., for instance,
\cite{Tondeur}).

In Section~\ref{s:do}, we have defined the transversal de Rham
differential $d_H$, acting from $C^\infty_c(M,\bigwedge^j N^*\cF)$
to $C^\infty_c(M,\bigwedge^{j+1} N^*\cF)$. The transverse
distribution $H$ naturally defines a transverse distribution
$HG\cong r^*H$ on the foliated manifold $(G,\cG)$ (cf.
Section~\ref{s:groupoid}) and the corresponding transversal de
Rham differential $d_H: C^\infty_c(G) \to C^\infty_c(G,
r^*N^*\cF)$ (see (\ref{e:d})).

An arbitrary leafwise half-density $\rho \in
C^\infty(M,|T\cF|^{1/2})$ can be written as
$\rho=f|\lambda|^{1/2}$ with $f\in C^\infty(M)$ and $\lambda\in
C^\infty(M,|T\cF|)$. Then $d_H\rho \in C^\infty_c(M, N^*\cF\otimes
|T\cF|^{1/2})$ is defined as
\[
d_H\rho=(d_Hf)|\lambda|^{1/2}+\frac{1}{2} f|\lambda|^{1/2}
k(\lambda).
\]
Any $f\in C^\infty_c(G,|T\cG|^{1/2})$ can be written as
$f=us^*(\rho)r^*(\rho)$, where $u\in C^\infty_c(G)$ and $\rho\in
C^\infty_c(M,|T\cF|^{1/2})$. The element $d_Hf\in C^\infty_c(G,
r^*N^*\cF\otimes |T\cG|^{1/2})$ is defined as
\[
d_Hf = d_Hus^*(\rho)r^*(\rho) + us^*(d_H\rho)r^*(\rho) +
us^*(\rho)r^*(d_H\rho).
\]
Finally, the operator $d_H$ has a unique extension to a
differentiation of the differential graded algebra
$\Omega_\infty=C^\infty_c(G, r^*\bigwedge N^*\cF\otimes
|T\cG|^{1/2})$. By definition, for any $f\in
C^\infty_c(G,|T\cG|^{1/2})$ and $\omega\in C^\infty_c(M,\bigwedge
N^*\cF)$, one has that
\[
d_H(fr^*\omega )=(d_Hf)r^*\omega + fr^*(d_H\omega ).
\]

\subsection{Transverse fundamental class of a foliation}\label{s:tfclass}
In this Section, we describe, following \cite{Co86} (cf. also
\cite{Co}), the simplest construction of a cyclic cocycle on the
algebra $C^\infty_c(G,|T\cG|^{1/2})$ associated with a foliated
manifold, and, namely, the construction of the transverse
fundamental class.

There is a general construction of cyclic cocycles on an arbitrary
algebra due to Connes \cite{Co:nc}. Let us call by a cycle of
dimension $n$ a triple $(\Omega, d, \int)$, where
$\Omega=\oplus_{j=0}^n\Omega^j$ is a graded algebra over $\CC$, $d
: \Omega \to \Omega$ is a graded differentiation of degree $1$
such that $d^2=0$, and $\int:\Omega^n\to \CC$ is a graded trace.
Thus, the maps $d$ and $\int$ satisfy the following conditions:
\begin{enumerate}
  \item
$d(\omega\omega')=d\omega\cdot\omega'+(-1)^{|\omega|}\omega\cdot
d\omega'$;
  \item $d^2=0$;
  \item $\int \omega_2\omega_1=(-1)^{|\omega_1||\omega_2|}\int
  \omega_1\omega_2$;
  \item $\int d\omega =0$ for $\omega\in\Omega^{n-1}$.
\end{enumerate}

Let $A$ be an algebra over $\CC$. A cycle over the algebra $A$ is
a cycle $(\Omega, d, \int)$ along with a homomorphism $\rho: A
\to\Omega^0$.

For any cycle $(A\stackrel{\rho}{\to}\Omega, d, \int)$ over $A$,
one can define its character as
\[
\tau(a^0, a^1, \ldots, a^n)=\int \rho(a^0)d(\rho(a^1))\ldots
d(\rho(a^n)).
\]
One can easily check that $\tau $ is a cyclic cocycle on $A$.
Moreover, any cyclic cocycle is the character of some cycle.

Now let $(M,\cF)$ be a manifold equipped with a transversally
oriented foliation, $H$ an arbitrary distribution on $M$,
transverse to $F=T\cF$. Consider the differential graded algebra
$(\Omega_\infty=C^\infty_c(G, r^*\bigwedge N^*\cF\otimes
|T\cG|^{1/2}), d_H)$. There is a closed graded trace $\tau$ on
$(\Omega_\infty, d_H)$ defined by
\[
\tau(\omega)=\int_M\left.\omega\right|_M,
\]
where $\omega\in \Omega_\infty^q=C^\infty_c(G, r^*\bigwedge^q
N^*\cF\otimes |T\cG|^{1/2})$. Here $\left.\omega\right|_M$ denotes
the restriction of $\omega$ to $M$. This is a section of the
bundle $\bigwedge^qN^*\cF\otimes |T\cF|$ on $M$, and its integral
over $M$ is well-defined, using the transverse orientation of the
foliation.

A problem is that, since $H$ is not integrable, it is not true, in
general, that $d^2_H=0$. Using (\ref{e:d}) and calculating the
$(0,2)$-component in the operator $d^2$, we get
\[
d^2_H=-(d_F\theta+\theta d_F).
\]
The operator $-(d_F\theta+\theta d_F)$ is a tangential
differential operator. Therefore, it can be written as
$R_{\bigwedge T^*M}(k)$ (cf. (\ref{rev:tang})), where $k\in
C^{\infty}_c(G,{\mathcal L}(\bigwedge T^*M))$ is given by the
exterior product with a vector-valued distribution $\Theta \in
{\cD}'(G, r^*\bigwedge^2 N^*\cF\otimes |T\cG|^{1/2})$ supported in
$G^{(0)}$. One can show that
\[
d^2_H\omega=\Theta \wedge \omega - \omega \wedge \Theta, \quad
\omega\in \Omega_\infty.
\]
Moreover, $d_H\Theta=0$.

Since the action of $d^2_H$ is given by the commutator with some
2-form, one can canonically construct a new differential graded
algebra $(\tilde{\Omega}_\infty, \tilde{d}_H)$ with
$\tilde{d}^2_H=0$ (see \cite{Co:nc,Co}). The algebra
$\tilde{\Omega}_\infty$ consists of $2\times 2$-matrices
$\omega=\{\omega_{ij}\}$ with elements from $\Omega_\infty$. An
element $\omega\in \tilde{\Omega}_\infty$ has degree $k$, if
$\omega_{11}\in \Omega_\infty^k$, $\omega_{12}, \omega_{21}\in
\Omega_\infty^{k-1}$ and $\omega_{22}\in \Omega_\infty^{k-2}$. The
product in $\tilde{\Omega}_\infty$ is given by
\[
\omega\cdot\omega'=
\begin{bmatrix}
\omega_{11} & \omega_{12} \\ \omega_{21} & \omega_{22}
\end{bmatrix}
\begin{bmatrix}
1 & 0 \\ 0 & \Theta
\end{bmatrix}
\begin{bmatrix}
\omega'_{11} & \omega'_{12} \\ \omega'_{21} & \omega'_{22}
\end{bmatrix},
\]
the differential by
\[
\tilde{d}_H\omega =\begin{bmatrix} d_H\omega_{11} & d_H\omega_{12}
\\ -d_H\omega_{21} & -d_H\omega_{22}
\end{bmatrix}
+
\begin{bmatrix}
0 & -\Theta \\ 1 & 0
\end{bmatrix}
\omega + (-1)^{|\omega|}\omega
\begin{bmatrix}
0 & 1 \\ -\Theta & 0
\end{bmatrix},
\]
and the closed graded trace $\tilde{\tau}:
\tilde{\Omega}_\infty^q\to \CC$ is defined as
\[
\tilde{\tau}\left(
\begin{bmatrix}
\omega_{11} & \omega_{12} \\ \omega_{21} & \omega_{22}
\end{bmatrix}
\right)=\tau(\omega_{11})-(-1)^q\tau(\omega_{22}\Theta).
\]
Finally, the homomorphism $\tilde{\rho}: C^\infty_c(G,
|T\cG|^{1/2}) \to \tilde{\Omega}^0_\infty$ is given by
\[
\tilde{\rho}(k)=\begin{bmatrix} k & 0 \\ 0 & 0
\end{bmatrix}.
\]
One can check that the triple $(\tilde{\Omega}_\infty,
\tilde{d}_H, \tilde{\tau})$ is a cycle over $C^\infty_c(G,
|T\cG|^{1/2})$. This cycle is called the fundamental cycle of the
transversally oriented foliation $(M,\cF)$. The character of this
cycle defines a cyclic cocycle $\phi_H$ on $C^\infty_c(G,
|T\cG|^{1/2})$. The cocycle $\phi_H$ depends on the auxiliary
choice of $H$, but the corresponding class in the cyclic
cohomology is independent of this choice. The class $[M/\cF]\in
HC^q(C^\infty_c (G, |T\cG|^{1/2}))$ defined by the cocycle
$\phi_H$ is called the transverse fundamental class of the
foliation $(M,\cF)$.

Let $C^\infty_c(G, |T\cG|^{1/2})^+$ denote the algebra
$C^\infty_c(G, |T\cG|^{1/2})$ with the adjoined unit. For an even
$q$, let us extend the cycle $(\tilde{\Omega}_\infty, d_H,
\tilde{\tau})$ to a cycle over $C^\infty_c(G, |T\cG|^{1/2})^+$,
putting $\tilde{\tau}(\Theta^{q/2})=0$. In \cite{Gor-CMP}, a
formula for a  cocycle in the  $(b,B)$-bicomplex of the algebra
$C^\infty_c(G, |T\cG|^{1/2})^+$, which gives the transverse
fundamental class of the foliation $(M,\cF)$, is obtained: the
formula
\[
\chi^r(k^0,k^1,\ldots,k^r)=\frac{(-1)^{\frac{q-r}{2}}}{\left(\frac{q+r}{2}\right)!}
\sum_{i_0+\ldots+i_r=\frac{q-r}{2}}\int_M k^0\,\Theta^{i_0}d_Hk^1
\ldots d_Hk^r\Theta^{i_r},
\]
where $r=q,q-2,\ldots$, $k^0 \in C^\infty_c(G, |T\cG|^{1/2})^+$,
$k^1,\ldots, k^r\in C^\infty_c(G, |T\cG|^{1/2})$, defines a
cocycle in the $(b,B)$-bicomplex of the algebra $C^\infty_c(G,
|T\cG|^{1/2})^+$. The class defined by $\chi$ in $HC^q
(C^\infty_c(G, |T\cG|^{1/2}))$ coincides with $[M/\cF]$.

The pairing with $[M/\cF]\in HC^q(C^\infty_c(G, |T\cG|^{1/2}))$
defines an additive map $K(C^\infty_c(G, |T\cG|^{1/2})) \to \CC$
--- ``integration over the transverse fundamental class''. An important
problem is the question of topological invariance of this map,
that is, the question of the possibility to extend it to an
additive map from $K(C^*_r(G))$ to $\CC$. A standard method of
solving this problem consists in constructing a smooth subalgebra
$\cB$ in $C^*_r(G)$, which contains the algebra $C^\infty_c(G,
|T\cG|^{1/2})$, such that the cyclic cocycle $\phi_H$ on
$C^\infty_c(G, |T\cG|^{1/2})$, which defines the transverse
fundamental class $[M/\cF]$, extends by continuity to a cyclic
cocycle on $\cB$ and defines, in that way, a map $K(\cB)\cong
K(C^*_r(G))\to\CC$. One managed to do this for so called
para-Riemannian foliations (cf. Section~\ref{para}). An example of
a para-Riemannian foliation is given by the lift of $\cF$ to a
foliation $\cV$ in the total space $P$ of the fibration $\pi : P
\to M$, whose fiber $P_x$ at $x\in M$ is the space of all
Euclidean metrics on the vector space $\tau_x=T_xM/T_x\cF$. Thus,
one can define a map $K(C^*(P,\cV))\to\CC$. Since the leaves of
the fibration $P$ are connected spin manifolds of nonpositive
curvature, one can construct an injective map $K(C^*_r(G)) \to
K(C^*(P,\cV))$, that allows to prove the existence of the map
$K(C^*_r(G))\to\CC$, given by the transverse fundamental class,
for an arbitrary foliation.

For geometrical consequences of this construction, see
\cite{Co86}.

\subsection{Cyclic cocycle defined by the Godbillon-Vey class}
Consider a smooth compact manifold $M$ equipped with a
transversally oriented codimension one foliation $\cF$. The
Godbillon-Vey class is a 3-dimensional cohomology class $GV\in
H^3(M,\RR)$. It is the simplest example of secondary
characteristic classes of the foliation. Recall its definition.
Since $\cF$ is transversally oriented, it is globally defined by a
non-vanishing smooth 1-form $\omega$ (that is, for any $x\in M$,
we have $\ker \omega_x=T_x\cF$). It follows from the Frobenius
theorem that there exists a 1-form $\alpha$ on $M$ such that
$d\omega=\alpha\wedge \omega$. One can check that the 3-form
$\alpha\wedge d\alpha$ is closed, and its cohomology class does
not depend on the choice $\omega$ and $\alpha$. This class $GV\in
H^3(M,\RR)$ is called the Godbillon-Vey class of $\cF$.

Let $T$ be a complete smooth transversal given by a good cover of
$M$ by regular foliated charts. Thus, $T$ is an oriented
one-dimensional manifold. In this Section, we construct a cyclic
cocycle on the algebra $C^\infty_c(G^T_T)$, corresponding to the
Godbillon-Vey class. This is done in several steps. To start with,
we recall the construction of the Godbillon-Vey class as a
secondary characteristic class associated with the cohomology
$H^*(W_1,\RR)$ of the Lie algebra $W_1=\RR[[x]]\partial_x$ of
formal vector field on $\RR$.

Recall (cf., for instance, \cite{Gui,Fuks}), that the cohomology
of a Lie algebra $\mathfrak{g}$ are defined as the cohomology of
the complex $(C^*(\mathfrak{g}), d)$, where:
\smallskip\par
$\bullet$ $C^q(\mathfrak{g})$ is the space of (continuous)
skew-symmetric
  $q$-linear functionals on $\mathfrak{g}$;
\smallskip\par
$\bullet $ for any $c\in C^q(\mathfrak{g})$, the cochain $dc\in
  C^{q+1}(\mathfrak{g})$ is defined by
 \begin{multline*}
dc(g_1,\ldots,g_{q+1})\\ =\sum_{1\leq i < j\leq q+1} (-1)^{i+j}
c([g_i,g_j],g_1,\ldots,\hat{g}_i,\ldots,\hat{g}_j, \ldots,
g_{q+1}),\\ g_1,\ldots,g_{q+1} \in \mathfrak{g}.
  \end{multline*}

For cochains on $W_1$, the continuity means the dependence only on
finite order jets of its arguments. The cohomology $H^*(W_1,\RR)$
were computed by Gelfand and Fuchs (see, for instance, the book
\cite{Fuks} and references therein). They are finite-dimensional
and the only nontrivial groups are $H^0(W_1,\RR)=\RR\cdot 1$ and
$H^3(W_1,\RR)=\RR\cdot gv$, where
\[
gv(p_1\partial_x, p_2\partial_x, p_3\partial_x)=
\begin{vmatrix}
p_1(0) & p_2(0) & p_3(0) \\ p'_1(0) & p'_2(0) & p'_3(0)\\ p''_1(0)
& p''_2(0) & p''_3(0)
\end{vmatrix}
, \quad p_1, p_2, p_3 \in \RR[[x]].
\]

Consider the bundle $F^k_+T \rightarrow T$ of positively oriented
frames of order $k$ on $T$ and the bundle
$F^\infty_+T=\underset{\leftarrow}{\lim} F^k_+T$ of positively
oriented frames of infinite order on $T$. By definition, a
positively oriented frame $r$ of order $k$ at $x\in T$ is the
$k$-jet at $0\in \RR$ of an orientation preserving diffeomorphism
$f$, which maps a neighborhood of $0$ in $\RR$ on some
neighborhood of $x=f(0)$ in $T$. If $y:U\to\RR$ is a local
coordinate on $T$, defined in a neighborhood $U$ of $x$, then the
numbers
\[
y_0=y(x),\quad y_1=\frac{d (y\circ f)}{dt}(0), \ldots,
y_k=\frac{d^k(y\circ f)}{dt^k}(0),
\]
are coordinates of the frame $r$, and, moreover, $y_1>0$.

There is defined a natural action of the pseudogroup $\Gamma^+(T)$
of orientation preserving local diffeomorphisms of $T$ on
$F^\infty_+T$. Let $\Omega^*(F^\infty_+T)^{\Gamma^+(T)}$ denote
the space of differential forms on $F^\infty_+T$, invariant under
the action of $\Gamma^+(T)$. There is a natural isomorphism of
differential algebras $J: C^*(W_1)\to \Omega^*
(F^\infty_+T)^{\Gamma^+(T)}$ defined in the following way. First
of all, let $v\in W_1$, and let $h_t$ be any one-parameter group
of local diffeomorphisms of $\RR$ such that $v$ is the
$\infty$-jet of the vector field
$\left.\frac{dh_t}{dt}\right|_{t=0}$. Then define a
$\Gamma^+(T)$-invariant vector field on $F^\infty_+T$, whose value
at $r=j_0f\in J^\infty_+(T)$ is given by $v(r)=j_0
\left(\left.\frac{d(f\circ h_t)}{dt}\right|_{t=0}\right)$.

For any $c\in C^q(W_1)$, put
\[
J(c)(v_1(r),\ldots,v_q(r))=c(v_1,\ldots,v_q).
\]
One can check that this isomorphism takes the cocycle $gv \in
C^3(W_1,\RR)$ to the three-form
\[
gv =\frac{1}{y_1^3}dy\wedge dy_1\wedge dy_2\in
\Omega^3(F^2_+T)^{\Gamma^+(T)}.
\]
Consider the bundle $F^\infty(M/\cF)$ on $M$, which consists of
infinite order jets of distinguished maps. There is a natural map
$F^\infty(M/\cF) \to F^\infty_+T/\Gamma^+(T)$. Using this map and
$\Gamma^+(T)$-invariance of $gv \in \Omega^3(F^2_+T)$, one can
pull back $gv$ to a closed form $gv(\cF)\in
\Omega^3(F^\infty(M/\cF))$. Since the fibers of the fibration
$F^2(M/\cF)\to M$ are contractible, $H^3(F^\infty(M/\cF),\RR)\cong
H^3(M,\RR)$, and the cohomology class in $H^3(M,\RR)$,
corresponding under this isomorphism to the cohomology class of
$gv(\cF)$ in $H^3(F^\infty(M/\cF),\RR)$, coincides with the
God\-bill\-ion-Vey class of $\cF$.

To construct the cyclic cocycle on $C^\infty_c(G^T_T)$,
corresponding to the God\-bill\-i\-on-Vey class, one uses a Van
Est type theorem (cf. \cite{Haefliger:CIME}). This theorem states
the existence of an embedding
\begin{equation}\label{e:vanest}
 \Omega^*(F^\infty_+T)^{\Gamma^+(T)}\to C^*(G^T_T,
\Omega^*(G^T_T)),
\end{equation}
where $C^*(G^T_T, \Omega^*(G^T_T))$ denotes the space of cochains
on the groupoid $G^T_T$ with values in the space $\Omega^*(G^T_T)$
of differential forms on $G^T_T$, which is a homomorphism of
complexes, inducing an isomorphism in the cohomology.

Let us describe the cocycle from $C^*(G^T_T, \Omega^*(G^T_T))$,
corresponding to the form $gv \in \Omega^3(F^2_+T)^{\Gamma^+(T)}$
(cf. \cite{Haefliger:CIME}). Let $\rho$ be an arbitrary smooth
positive density on $T$. Define a homomorphism $\delta : G^T_T\to
\RR^*_+$, called the modular homomorphism, setting
\[
\delta(\gamma)=\frac{\rho_{y}}{dh_\gamma(\rho_{x})},\quad
\gamma\in G^T_T, \quad \gamma : x \to y,
\]
where the map $dh_\gamma : |\tau_x|\to |\tau_y|$ is induced by the
linear holonomy map, and also a homomorphism $\ell=\log \delta :
G^T_T\to \RR$.

The formula
\[
c(\gamma_1, \gamma_2) =\ell(\gamma_2)\, d\ell (\gamma_1)
-\ell(\gamma_1)\, d\ell(\gamma_2), \quad \gamma_1,\gamma_2\in
G^T_T,
\]
defines a 2-cocycle on $G^T_T$ with values in the space of 1-forms
on $G^T_T$. This cocycle, called the Bott-Thurston cocycle,
corresponds to $gv \in \Omega^3(F^2_+T)^{\Gamma^+(T)}$ under the
isomorphism given by the embedding (\ref{e:vanest}).

Finally, the Bott-Thurston cocycle $c$ defines a cyclic 2-cocycle
$\psi$ on $C^\infty_c(G^T_T)$ by the formula (cf.
Example~\ref{ex:phi}):
\[
\psi(k^0,k^1,k^2)=\int_{\gamma_0\gamma_1\gamma_2\in T}
k^0(\gamma_0)k^1(\gamma_1)k^2(\gamma_2)c(\gamma_1,\gamma_2), \quad
k^0,k^1,k^2 \in C^\infty_c(G^T_T).
\]
This is the cyclic cocycle, corresponding to the Godbillon-Vey
class of $\cF$.

There is another description of cyclic cocycle, corresponding to
the God\-bill\-i\-on-Vey class of $\cF$, which relates it with
invariants of the von Neumann algebra of this foliation
\cite{Co86}. The following formula defines a cyclic $1$-cocycle on
$C^\infty_c(G^T_T)$:
\[
\tau(k^0,k^1)=\int_{G^T_T}k^0(\gamma^{-1})dk^1(\gamma),\quad k^0,
k^1 \in C^\infty_c(G^T_T).
\]
The isomorphism $HC^1(C^\infty_c(G^T_T)) \cong
HC^1(C^\infty_c(G))$, defined by the strong Morita equivalence
(cf. Example~\ref{ex:morita}), associates the class of $\tau$ in
the group $HC^1(C^\infty_c(G^T_T))$ to the transverse fundamental
class of the foliation $(M,\cF)$ in $HC^1(C^\infty_c(G))$.

The fixed smooth positive density $\rho$ on $T$ defines a faithful
normal semi-finite weight $\phi_\rho$ on the von Neumann algebra
$W^*(G^T_T)$ of the groupoid $G^T_T$. For any $k\in
C^\infty_c(G^T_T)$, the value of the weight $\phi_\rho$ is given
by
\[
\phi_\rho(k)=\int_Tk(x)\,\rho(x).
\]
Consider the one-parameter group $\sigma_t$ of automorphisms of
the von Neumann algebra $W^*(G^T_T)$ given by
\[
\sigma_t(k)(\gamma)=\delta(\gamma)^{it}k(\gamma), \quad k\in
C^\infty_c(G^T_T), \quad t\in \RR.
\]
This group is the group of modular automorphisms, associated with
the weight, by the Tomita-Takesaki theory \cite{TakesakiLNM}.

The importance of the group of modular automorphisms is explained,
in particular, by the following its characterization
\cite{TakesakiLNM}: a one-parameter group of $\ast$-automorphisms
$\sigma_t$ of a von Neumann algebra $M$ is the group of modular
automorphisms, associated with a weight $\omega$, if and only if
$\omega$ satisfies the Kubo-Martin-Schwinger condition with
respect to $\sigma_t$, that is, there is an analytic in the strip
$\operatorname{Im} z\in (0,1)$ and continuous in its closure
function $f$ such that, for any $a,b\in M, t\in\RR$,
\[
f(t)=\omega(\sigma_t(a)b), \quad f(t+i)=\omega(b\sigma_t(a)).
\]

Following Connes \cite{Co86}, we call by a $1$-trace on a Banach
algebra $B$ a bilinear functional $\phi$, defined on a dense
subalgebra $\cA\subset B$, such that
\begin{enumerate}
  \item $\phi$ is a cyclic cocycle on $\cA$;
  \item for any $a^1\in\cA$, there is a constant $C>0$ such that
\[
|\phi(a^0,a^1)|\leq C\|a^0\|, \quad a^0\in \cA,
\]
\end{enumerate}
and a $2$-trace on $B$ a trilinear functional $\phi$, defined on a
dense subalgebra $\cA\subset B$, such that
\begin{enumerate}
  \item $\phi$ is a cyclic cocycle on $\cA$;
  \item for any $a^1,a^2\in\cA$, there is a constant $C>0$ such
  that
\[
|\phi(x^0,a^1x^1,a^2)-\phi(x^0a^1,x^1,a^2)|\leq C\|a^1\|\,\|a^2\|,
\quad x^0, x^1\in \cA.
\]
\end{enumerate}
The formula
\[
\dot{\tau}(k^0,k^1)=\lim_{t\to 0}\frac{1}{t}
\left(\tau(\sigma_t(k^0),\sigma_t(k^1)) -\tau(k^0,k^1) \right),
\quad k^0, k^1\in C^\infty_c(G^T_T)
\]
defines a 1-trace on $C^*_r(G^T_T)$ with the domain
$C^\infty_c(G^T_T)$, invariant under the action of the
automorphism group $\sigma_t$.

For any $1$-trace $\phi$ on $C^*$-algebra $A$, invariant under an
action of an one-parameter automorphism group $\alpha_t$ with the
generator $D$, such that the space $\dom \phi\cap\dom D$ is dense
in $A$, a $2$-trace $\chi=i_D\phi$ on $C^*_r(G^T_T)$ (an analogue
of the contraction) is defined by
\[
\chi(a^0,a^1,a^2)=\phi(D(a^2)a^0,a^1)-\phi(a^0D(a^1),a^2), \quad
a^0,a^1,a^2 \in \dom \phi\cap\dom D.
\]

\begin{thm}
Let $(M,\cF)$ be a manifold with a transversally oriented
codimension one foliation, $T$ a complete smooth transversal,
$\rho$ a smooth positive density on $T$. Then the cyclic cocycle
$\psi\in HC^2(C^\infty_c(G^T_T))$, corresponding to the
Godbillon-Vey class of $\cF$, coincides with $i_D\dot{\tau}$.
\end{thm}

As a consequence of this statement, Connes obtained the following
geometrical fact. First, recall that one can naturally associate
to any von Neumann algebra $M$ an action of the multiplicative
group $\RR^*_+$ (called the flow of weights \cite{Connes-Tak}) on
some commutative von Neumann algebra, the center of the crossed
product $M\rtimes \RR$ of $M$ by $\RR$ relative to the action of
$\RR$ on $M$ given by the modular automorphism group $\sigma_t$.

\begin{thm}\cite{Co86}
Let $(M,\cF)$ be a manifold with a transversally oriented
codimension one foliation. If the Godbillon-Vey class $GV\in
H^3(M,\RR)$ differs from zero, then the flow of weights of the von
Neumann algebra of the foliation $\cF$ has a finite invariant
measure.
\end{thm}

In particular, this implies the following, previous result.

\begin{thm}\cite{Hurder-KatokBAMS}
Let $(M,\cF)$ be a manifold with a transversally oriented
codimension one foliation. If the Godbillon-Vey class $GV\in
H^3(M,\RR)$ differs from zero, then the von Neumann algebra of the
foliation $\cF$ has the non trivial type $III$ component.
\end{thm}

In \cite{Mor,Mor-Natsume}, one gives another construction of the
cyclic cocycle, associated with the Godbillon-Vey class, as a
cyclic cocycle on the $C^*$-algebra of foliation $C^*_r(G)$, in
the case of a foliated $S^1$-bundle, that is, when the foliation
$\cF$ is constructed, using the suspension construction for a
homomorphism $\phi : \Gamma= \pi_1(B) \to
\operatorname{Diff}(S^1)$.

\subsection{General constructions of cyclic cocycles} \label{s:second}
The secondary characteristic classes of foliations are given by
the characteristic homomorphism
\cite{Bernshtein-R72,Bott-H,Haefliger:73}
\[
\chi_{\cF}: H^*(W_q; O(q)) \to H^*(M,\RR).
\]
defined for any codimension $q$ foliation $\cF$ on a smooth
manifold $M$, where $H^*(W_q; O(q))$ denotes the relative
cohomology of the Lie algebra $W_q$ of formal vector fields in
$\RR^q$. A basic property of the secondary characteristic classes
consists in their functoriality: if a smooth map $f: N\to M$ is
transverse to a transversely oriented foliation $\cF$ and $f^*\cF$
is the foliation on $N$ induced by $f$ (by definition, the leaves
of $f^*\cF$ are the connected components of the pre-images of the
leaves of $\cF$ under the map $f$), then
\[
f^*(\chi_{\cF}(\alpha))=\chi_{f^*\cF}(\alpha), \quad \alpha\in
H^*(W_q; O(q)).
\]
The classifying space $B\Gamma_q$ of the groupoid $\Gamma_q$
classifies codimension $q$ foliations on a given manifold $M$ in
the sense that any foliation $\cF$ on $M$ defines a map $M\to
B\Gamma_q$, and, moreover, in the case when $M$ is compact, a
homotopy class of maps $M\to B\Gamma_q$ corresponds to a
concordance class of foliations on $M$ \cite{Haefliger72}.

By the functoriality of the characteristic homomorphism, it
suffices to know it for the universal foliation on $B\Gamma_q$,
that gives the universal characteristic homomorphism
\[
\chi: H^*(W_q; O(q)) \to H^*(B\Gamma_q,\RR).
\]
For any complete transversal $T$, the universal characteristic
homomorphism $\chi$ is represented as a composition
\[
\chi: H^*(W_q; O(q)) \to H^*(BG^T_T,\RR) \to H^*(B\Gamma_q,\RR).
\]
Since the groupoids $G^T_T$ and $G$ are equivalent,
$H^*(BG^T_T,\RR)\cong H^*(BG,\RR)$, that defines a map
\[
H^*(W_q; O(q)) \to H^*(BG,\RR).
\]

The relative cohomology $H^*(W_q; O(q))$ of the Lie algebra $W_q$
are calculated as follows (cf., for instance, \cite{Godbillon73}).
Consider the complex
\[
WO_q=\Lambda(u_1,u_3,\ldots, u_m)\otimes P(c_1,c_2,\ldots,c_q)_q
\]
where $m$ is the largest odd number, less than $q+1$,
$\Lambda(u_1,u_3,\ldots, u_m)$ denotes the exterior algebra with
generators $u_1,u_3,\ldots, u_m$, $\operatorname{deg} u_i=2i-1$,
$P(c_1,c_2,\ldots,c_q)$ denotes the polynomial algebra with
generators $c_1,c_2,\ldots,c_q$, $\operatorname{deg} c_i=2i$ and
$P(c_1,c_2,\ldots,c_q)_q$ denotes the quotient of
$P(c_1,c_2,\ldots,c_q)$ by the ideal generated by the monomials of
degree greater than $2q$. The differential in $WO_q$ is defined by
$d(u_i\otimes 1)=1\otimes c_i$, $d(1\otimes c_i)=0$. There is a
homomorphism of the complex $WO_q$ to the complex $C^*(W_q; O(q))$
of relative cochains on the algebra $W_q$, which induces an
isomorphism in the cohomology $H^*(WO_q)\cong H^*(W_q; O(q))$.

For any $i=1,2,\ldots, q$, the cohomology class $[1\otimes c_i]\in
H^{2i}(WO_q)$, defined by $1\otimes c_i$, is mapped by the
characteristic homomorphism $\chi_{\cF}$ to the $i$-th Pontryagin
class $p_i$ of the normal bundle $\tau$. The Godbillon-Vey class
$GV\in H^3(M,\RR)$ of a codimension one foliation on a compact
manifold $M$ is the image of $[u_1\otimes c_1]\in H^3(WO_1)$ under
the characteristic homomorphism $\chi_{\cF}$.

Connes \cite[Section 2 $\delta$, Theorem 14 and Remark b)]{Co}
constructed a natural map
\[
\Phi_*: H^*(BG,\RR)\to HP^*(C^\infty_c(G,|T\cG|^{1/2}))
\]
for any oriented etale groupoid $G$. The constructions of the
transverse fundamental class of a foliation and of the cyclic
cocycle associated with the Godbillon-Vey class are special cases
of this general construction.

In \cite{Gor:topology}, Gorokhovsky generalized the above
construction of the cyclic cocycle associated with the
Godbillon-Vey class, which uses the group of modular
automorphisms, to the case of arbitrary secondary characteristic
classes. This construction makes an essential use of the cyclic
cohomology theory for Hopf algebras developed in the paper
\cite{CoM:Hopf} (cf. Section~\ref{s:locind}).

For the calculations of the cyclic cohomology of $C^*$-algebras of
etale groupoids, cf.
\cite{Brylinsky:N,Co,crainic:cyclic,CrainicM00}.

\subsection{Index of tangentially elliptic operators}
In this Section, we briefly describe applications of the above
methods to the index theory for tangentially elliptic operators.

Let $D$ be a tangentially elliptic operator on a compact foliated
manifold $(M,\cF)$. The restrictions of $D$ to the leaves of $\cF$
define a family $(D_l)_{l\in M/\cF}$, where $D_l$ is an elliptic
differential operator on a leaf $l$ of $\cF$. It turns out that
the families $(P_{\Ker D_l})_{l\in M/\cF}$ and $(P_{\Ker
D^*_l})_{l\in M/\cF}$, which consist of the orthogonal projections
to $\Ker D_l$ and $\Ker D^*_l$ in $L^2(l)$ accordingly, define
elements of the von Neumann algebra of the foliation $W^*(M,\cF)$.
Suppose that $\cF$ has a holonomy invariant transverse measure
$\Lambda$. Then a normal semi-finite trace $\trLambda$ on
$W^*(M,\cF)$ is defined. It is proved in \cite{Co79} that the
dimensions
\begin{gather*}
\dim_\Lambda (({\Ker D_l})_{l\in M/\cF}) = \trLambda ((P_{\Ker
D_l})_{l\in M/\cF}), \\ \dim_\Lambda (({\Ker D^*_l})_{l\in M/\cF})
= \trLambda ((P_{\Ker D^*_l})_{l\in M/\cF}),
\end{gather*}
are finite, and, therefore, the index of $D$ is defined by
\[
\operatorname{ind}_\Lambda (D)=\dim_\Lambda (({\Ker D_l})_{l\in
M/\cF}) - \dim_\Lambda (({\Ker D^*_l})_{l\in M/\cF}).
\]
Like in \cite{ASIV}, the tangential principal symbol $\sigma_D$ of
$D$ defines an element of $K^1(T^*\cF)$. Using the Thom
isomorphism (for simplicity, assume that the bundle $T\cF$ is
orientable), one can define the Chern character $\operatorname{ch}
(\sigma_D)$ as an element of the rational cohomology group
$H^*(M,\QQ)$ of the manifold $M$. Let
$\operatorname{Td}(T^*\cF)\in H^*(M,\QQ)$ denote the Todd class of
the cotangent bundle $T^*\cF$ to $\cF$.

\begin{thm}\cite{Co79}\label{t:Co79}
One has the formula
\[
\operatorname{ind}_\Lambda (D)=(-1)^{p(p+1)/2}\langle C,
\operatorname{ch} (\sigma_D) \operatorname{Td}(T^*\cF) \rangle,
\]
where $C$ is the Ruelle-Sullivan current, corresponding to the
transverse measure $\Lambda$.
\end{thm}

Theorem~\ref{t:Co79} is completely similar to the Atiyah-Singer
index theorem in cohomological form \cite{ASI} with the only
difference that one uses here the pairing with the Ruelle-Sullivan
current $C$ instead of the integration over the compact manifold
in the right-hand side of the Atiyah-Singer formula. An odd
version of Theorem~\ref{t:Co79} is proved in the paper~\cite{DHK2}
(see also. \cite{DHK}). It is related with the index theory for
Toeplitz operators.

In the paper \cite{Co-skandal} (see also \cite{Co:survey}), a
$K$-theoretic version of the index theorem for tangentially
elliptic operators on an arbitrary compact foliated manifold
$(M,\cF)$ is proved. Let $D$ be a tangentially elliptic operator
on a compact manifold $M$. One can construct an analytical index
${\rm Ind}_a(D)\in K_0(C^*(G))$ of the operator $(D_l)_{l\in
M/\cF}$, using operator constructions (\cite{Co79,Co:survey}),
and, starting with the class $[\sigma_D]\in K^1(T^*\cF)$ defined
by the tangential principal symbol $\sigma_D$ of the operator $D$,
one can construct its topological index ${\rm Ind}_t(D)\in
K_0(C^*(G))$.

\begin{thm}\cite{Co-skandal,Co:survey}\label{t:CS}
For any tangentially elliptic operator $D$ on a compact foliated
manifold $(M,\cF)$, one has the formula
\[
{\rm Ind}_a(D)={\rm Ind}_t(D).
\]
\end{thm}

If $\cF$ is given by the fibers of a fibration $M\to B$, then
$K_0(C^*(G))\cong K^0(B)$, and Theorem~\ref{t:CS} is reduced to
the Atiyah-Singer theorem for families of elliptic operators
\cite{ASIV}.

In \cite{H-L} (see also \cite{H-L91}), Heitsch and Lazarov proved
an analogue of the Atiyah-Bott-Lefschetz fixed point formula
\cite{Atiyah-Bott} in the setting of Theorem~\ref{t:Co79}. More
precisely, they consider a compact foliated manifold $(M,\cF)$,
equipped with a holonomy invariant transverse measure $\Lambda$,
and a diffeomorphism $f\colon M\to M$, which takes each leaf of
the foliation to itself. They assume additionally that the fixed
point sets of $f$ are submanifolds, transverse to the foliation.
Under these assumptions, the Lefschetz theorem proved in
\cite{H-L} states the equality of the alternating sum of the
traces of the action of $f$ on the $L^2$ spaces of harmonic forms
along the leaves of the foliation (an appropriate trace is given
by $\trLambda$) with the average with respect to $\Lambda$ of
local contributions of the fixed point sets.

In \cite{Benameur97}, an equivariant generalization of
Theorem~\ref{t:CS} to the case when there is a compact Lie group
action, taking each leaf of the foliation $\cF$ to itself, is
obtained. As a consequence, the author extended the Lefschetz
theorem proved in \cite{H-L} to the case of an arbitrary
tangentially elliptic complex when the diffeomorphism $f\colon
M\to M$ is included to a compact Lie group action, taking each
leaf of $\cF$ to itself.

Finally, in \cite[Chapter III, Section 7.$\gamma$, Corollary
13]{Co}, a cyclic version of the index theorem for tangentially
elliptic operators is given.

\begin{thm}
Let $D$ be a tangentially elliptic operator on a compact manifold
with a transversely oriented foliation $(M,\cF)$, ${\rm Ind}_a(D)$
its analytic index. Then, for any $\omega\in H^*(BG)$, we have
\[
\langle \Phi_*(\omega),{\rm Ind}_a(D)
\rangle=\langle\omega,\rm{ch}_\tau(\sigma_D)\rangle.
\]
\end{thm}

Here $\Phi_* : H^{*}(BG) \to HP^*(C^\infty_c(G,|T\cG|^{1/2}))$ is
the map introduced in Section~\ref{s:second},
$\rm{ch}_\tau(\sigma_D)$ denotes the twisted Chern character
\[
\rm{ch}_\tau(\sigma_D)=Td(\tau)^{-1}\rm{ch}(\sigma_D),
\]
where ${\rm ch}: K_{\ast}(BG)\to H_*(BG,\QQ)$ is the Chern
character.

Let us also mention the papers
\cite{Mor-Natsume,Heitsch95,Heitsch-L99,Gor-Lott03,Gor-Lott04},
concerning to generalizations of the local index theorem for
families of elliptic operators \cite{BismutInv85} to the case of
tangentially elliptic operators on a foliated manifold, and the
papers \cite{Benameur95,Benameur02,Benameur03}, concerning to
cyclic versions of the Lefschetz formula for diffeomorphisms,
which take each leaf of a foliation to itself.

\section{Noncommutative differential geometry}
In this Section, we will describe noncommutative analogues of two
of the most important transverse geometric structures for foliated
manifolds: symplectic and Riemannian ones.

\subsection{Noncommutative symplectic geometry}
Based on the ideas of the deformation theory of Gerstenhaber
\cite{Gersten63}, Ping Xu \cite{Xu} and Block and Getzler
\cite{Block-Ge} introduced a noncommutative analogue of the
Poisson bracket. Namely, they defined a Poisson structure on an
algebra $A$ as a Hochschild $2$-cocycle $P\in Z^2(A,A)$ such that
$P\circ P$ as a Hochschild $3$-coboundary, $P\circ P\in B^3(A,A)$.
In other words, a Poisson structure on $A$ is given by a linear
map $P : A\otimes A\to A$ such that
\begin{equation}\label{e:p1}
(\delta P)(a_1,a_2,a_3)\equiv
a_1P(a_2,a_3)-P(a_1a_2,a_3)+P(a_1,a_2a_3)-P(a_1,a_2)a_3=0,
\end{equation}
and there is a $2$-cochain $P_1:A\otimes A\to A$ such that
\begin{multline}\label{e:p2}
P\circ P(a_1,a_2,a_3)\equiv P(a_1,P(a_2,a_3))-P(P(a_1,a_2),a_3)\\
=a_1P_1(a_2,a_3)- P_1(a_1a_2,a_3)+P_1(a_1,a_2a_3)-P_1(a_1,a_2)a_3.
\end{multline}
The identity (\ref{e:p1}) is an analogue of the Jacobi identity
for a Poisson bracket, and the identity (\ref{e:p2}) is an
analogue of the Leibniz rule.

A connection of this definition with the deformation theory
consists in the following fact. Let us call by a formal
deformation of an algebra $A$ an associative multiplication on the
vector space $A[[\nu]]=\CC[[\nu]]\otimes A$ over the field
$\CC[[\nu]]$ of formal complex series such that the induced
multiplication on $A=A[[\nu]]/\nu A[[\nu]]$ coincides with the
multiplication in $A$. Such a deformation is described by a
cochain
\[
m(a_1,a_2)=\sum_{i=0}^\infty \nu^im_i(a_1,a_2), \quad a_1, a_2\in
A[[\nu]],
\]
and, moreover,
\begin{enumerate}
  \item $m\circ m(a_1,a_2,a_3)=m(a_1,m(a_2,a_3)) -m(m(a_1,a_2),a_3)=0;$
  \item $m_0(a_1,a_2)=a_1a_2.$
\end{enumerate}
It easily follows from here that $m_1$ defines a Poisson structure
on $A$.

Block and Getzler \cite{Block-Ge} defined a Poisson structure on
the operator algebra $C^\infty_c(G, |T\cG|^{1/2})$ of a symplectic
foliation $\cF$ in the case when the normal bundle $\tau$ to $\cF$
has a basic connection $\nabla $ (recall that a basic connection
on $\tau$ is a holonomy invariant adapted connection).

Consider an arbitrary symplectic foliation $\cF$. Let $\omega $ be
the corresponding closed $2$-form of constant rank. It defines a
nondegenerate holonomy invariant $2$-form on the normal bundle
$\tau$, which gives a bundle isomorphism of $\tau$ and
$\tau^*\cong N^*\cF$ and, therefore, defines a $2$-form $\Lambda$
on $N^*\cF$.

Choose an arbitrary distribution $H$ on $M$, transverse to $\cF$.
There is defined the transverse differential
\[
d_H: \Omega_\infty^0=C^\infty_c(G, |T\cG|^{1/2}) \to
\Omega_\infty^1 =C^\infty_c(G, r^*N^*\cF\otimes |T\cG|^{1/2}).
\]
The exterior product of differential forms gives a map
$\Omega_\infty^1 \times \Omega_\infty^1 \stackrel{\wedge}{\to}
\Omega_\infty^{2}$ (cf. (\ref{e:product})). Finally, the pairing
with the $2$-form $\Lambda$ defines a map
\[
\Omega_\infty^{2}\to \Omega_\infty^0 : s\to \langle\Lambda, s
\rangle,
\]
where
\[
\langle\Lambda, s \rangle (\gamma) = \langle \Lambda_{s(\gamma)},
s(\gamma) \rangle.
\]
A Poisson bracket $P(k_1,k_2)$ of $k_1, k_2\in C^\infty_c(G,
|T\cG|^{1/2})$ is defined as
\begin{equation}\label{e:bracket}
P(k_1,k_2)=\langle\Lambda, d_Hk_1\wedge d_Hk_2\rangle.
\end{equation}
It is not difficult to check that $P$ satisfies (\ref{e:p1}).

A construction of a $2$-cochain $P_1:A\otimes A\to A$, which
satisfies (\ref{e:p2}), is given in \cite{Block-Ge} only in the
case when $\cF$ has a basic connection.

\subsection{Noncommutative Riemannian spaces}\label{s:riem}
According to \cite{Co-M,Sp-view}, the initial datum of the
noncommutative Riemannian geometry is a spectral triple.

\begin{defn}
A spectral triple is a set $({\mathcal A}, {\mathcal H}, D)$,
where:
\begin{enumerate}
\item ${\mathcal A}$ is an involutive algebra;
\item ${\mathcal H}$ is a Hilbert space, equipped with a
$\ast$-representation of ${\mathcal A}$;
\item $D$ is an (unbounded) self-adjoint operator in ${\mathcal H}$
such that
\medskip
\par
\noindent 1. for any $a\in {\mathcal A}$, the operator
$a(D-i)^{-1}$ is a compact operator in ${\mathcal H}$;
\medskip\par
\noindent 2. $D$ almost commutes with elements of ${\mathcal A}$
in the sense that $[D,a]$ is bounded for any $a\in {\mathcal A}$.
\end{enumerate}

A spectral triple is called even, if $\cH$ is endowed with a
$\ZZ_{2}$-grading $\gamma \in \cL(\cH)$, $\gamma=\gamma^{*}$,
$\gamma^{2}=1$, and, moreover,  $ \gamma D=-D\gamma$ and $\gamma a
=a \gamma$ for any $a\in \cA$. In the opposite case, a spectral
triple is called odd.
\end{defn}

Spectral triples were considered for the first time in the paper
\cite{Baaj-Julg}, where they were called unbounded Fredholm
modules. A spectral triple $(\cA, \cH, D)$ defines a Fredholm
module $(\cH,F)$ over $\cA$, where $F=D(I+D^2)^{-1/2}$
\cite{Baaj-Julg}. On the contrary, any Fredholm module $(\cH,F)$
over $\cA$ is homotopic to a Fredholm module determined by a
spectral triple. In a sense, the operator $F$ is connected with
measurement of angles and is responsible for the conformal
structure, whereas $|D|$ is connected with measurement of lengths.

\begin{defn}\label{d:dim}
A spectral triple $({\mathcal A}, {\mathcal H}, D)$ is called
$p$-summable (or $p$-dimensional), if, for any $a\in {\mathcal
A}$, the operator $a(D-i)^{-1}$ is an element of the Schatten
class $\cL^p({\mathcal H})$.

A spectral triple $({\mathcal A}, {\mathcal H}, D)$ is called
finite-dimensional, if it is $p$-summable for some $p$.

The greatest lower bound of all $p$'s, for which a
finite-dimensional spectral triple is $p$-summable, is called the
dimension of the spectral triple.
\end{defn}

As shown in \cite{Baaj-Julg}, the dimension of the spectral triple
$({\mathcal A}, {\mathcal H}, D)$ coincides with the dimension of
the corresponding Fredholm module $(\cH,F)$, $F=D(I+D^2)^{-1/2}$,
over $\cA$.

The classical Riemannian geometry is described by the spectral
triple $(\cA,{\mathcal H},D)$, associated with a compact
Riemannian manifold $(M,g)$:
\begin{enumerate}
\item An involutive algebra ${\mathcal A}$ is the algebra $C^{\infty}(M)$
of smooth functions on $M$;
\item A Hilbert space ${\mathcal H}$ is the space
$L^2(M,\Lambda^*T^*M)$, on which the algebra ${\mathcal A}$ acts
by multiplication;
\item An operator $D$ is the signature operator $d+d^*$.
\end{enumerate}

Let us show that this spectral triple contains a basic geometric
information on the Riemannian manifold $(M,g)$. First of all, it
is finite-dimensional, and the dimension of this spectral triple
coincides with the dimension of $M$. This fact is an immediate
consequence of the Weyl asymptotic formula for eigenvalues of
self-adjoint elliptic operators on a compact manifold.

For any $T\in {\mathcal L}({\mathcal H})$, denote by $\delta^i(T)$
the $i$-th iterated commutator with $|D|$. Consider the space
${\rm OP}^0$, consisting of all $T\in {\mathcal L}({\mathcal H})$
such that $\delta^i(T) \in {\mathcal L}({\mathcal H})$ for any
$i\in\NN$. Then ${\rm OP}^0\cap C(M)$ coincides with
$C^\infty(M)$. This allows to reconstruct the smooth structure of
$M$, based on its topological structure and the spectral triple
$(\cA,{\mathcal H},D)$ (observe that here one can take as $\cA$
any involutive algebra, which consists of Lipschitz functions and
is dense in $C(M)$).

Based on the spectral triple, one can compute the geodesic
distance $d(x,y)$ between any two points $x, y\in M$ by the
formula (cf. \cite{CoErgTh})
\[
d(x,y)=\sup \{|f(x)-f(y)|: f\in\cA, \|[D,f]\|\leq 1\}.
\]
Then, starting from the triple $(A,{\mathcal H},D)$, one can
reconstruct the Riemannian volume form $d\nu$, given in local
coordinates by $d\nu=\sqrt{\det g}\,dx$. For this, one uses a
trace $\operatorname{Tr}_\omega$, introduced by Dixmier in
\cite{Dixmier66} as an example of a nonstandard trace on the
algebra $\cL(\cH)$. Before we describe the Dixmier trace, we
introduce some auxiliary notions. Consider the ideal
$\cL^{1+}(\cH)$ in the algebra of compact operators $\cK(\cH)$,
which consists of all $T\in \cK(\cH)$ such that
\[
\sup_{N\in\NN}\frac{1}{\ln N}\sum_{n=1}^N\mu_n(T)<\infty,
\]
where $\mu_1(T)\geq \mu_2(T)\geq \ldots$ are the singular numbers
of $T$. For any invariant mean $\omega$ on the amenable group of
upper triangular $2\times 2$-matrices, Dixmier constructed a
linear form $\lim_\omega$ on the space ${\ell}^\infty(\NN)$ of
bounded sequences, which coincides with the limit functional
$\lim$ on the subspace of convergent sequences. The trace
$\operatorname{Tr}_\omega$ on $\cL^{1+}(\cH)$ is defined for a
positive operator $T\in \cL^{1+}(\cH)$ as
\[
\operatorname{Tr}_\omega(T)=\lim_\omega \frac{1}{\ln N}
\sum_{n=1}^N\mu_n(T).
\]

Let $M$ be a compact manifold, $E$ a vector bundle on $M$ and
$P\in\Psi^m(M,E)$ a classical pseudodifferential operator. Thus,
in any local coordinate system, its complete symbol $p$ can be
represented as an asymptotic sum $p\sim p_m+p_{m-1}+\ldots$, where
$p_l(x,\xi)$ is a homogeneous function of degree $l$ in $\xi$. As
shown in \cite{Co-action}, the Dixmier trace
$\operatorname{Tr}_\omega(P)$ does not depend on the choice of
$\omega$ and coincides with the value of the residue trace
$\tau(P)$, introduced by Wodzicki \cite{Wo} and Guillemin
\cite{Gu85}. The residue trace $\tau$ is defined as follows. For
$P\in\Psi^{*}(M,E)$, its residue form $\rho_P$ is defined in local
coordinates as $$ \rho_P = \left(\int_{|\xi|=1}\Tr
p_{-n}(x,\xi)\,d\xi\right) |dx|. $$ The density $\rho_P$ turns out
to be independent of the choice of a local coordinate system and,
therefore, gives a well-defined density on $M$. The integral of
the density $\rho_P$ over $M$ is the residue trace $\tau(P)$ of
$P$:
\begin{equation}\label{e:wodz}
\tau(P)=(2\pi)^{-n}\int_M\rho_P =(2\pi)^{-n}\int_{S^*M}\Tr
p_{-n}(x,\xi)\, dx d\xi.
\end{equation}
Wodzicki \cite{Wo} showed that $\tau$ is a unique trace on the
algebra $\Psi^*(M,E)$ of classical pseudodifferential operators of
arbitrary order.

According to \cite{Co-action} (cf. also \cite{Gracia:book}), any
$P\in\Psi^{-n}(M,E)$ ($n=\dim M$) belongs to the ideal
$\cL^{1+}(L^2(M,E))$ and, for any invariant mean $\omega$,
\[
\operatorname{Tr}_\omega(P)=\tau(P).
\]
The above results imply the formula
\[
\int_Mf\,d\nu=c(n)\operatorname{Tr}_\omega(f|D|^{-n}), \quad f\in
\cA,
\]
where $c(n)=2^{(n-[n/2])} \pi^{n/2} \Gamma(\frac{n}{2}+1)$. Thus,
the Dixmier trace $\operatorname{Tr}_\omega$ can be considered as
a proper noncommutative generalization of the integral.

Finally, the Egorov theorem for pseudodifferential operators
allows to describe the geodesic flow on the cotangent bundle
$T^*M$ in terms of the given spectral triple (cf. Section
\ref{s:geo}).

\begin{ex}
Let us give examples of spectral triples $({\mathcal A},{\mathcal
H},D)$ associated with the noncommutative torus $T^2_\theta$ (see
\cite{Connes-gravity}). These triples are parametrized by a
complex number $\tau$ with ${\rm Im}\,\tau
>0$. Put
\[
\cA=\cA_\theta=\{a = \sum_{(n,m)\in\ZZ^2}a_{nm}U^nV^m : a_{nm}\in
\cS(\ZZ^2)\}.
\]
Define a canonical normalized trace $\tau_0$ on $\cA_\theta$ as
\[
\tau_0(a)=a_{00},\quad a\in \cA_\theta.
\]
Let $L^2(\cA_\theta,\tau_0)$ be the Hilbert space, which is the
completion of $\cA_\theta$ in the inner product
$(a,b)=\tau_0(b^*a), a,b\in \cA_\theta$. The Hilbert space $\cH$
is defined as the sum of two copies of $L^2(\cA_\theta,\tau_0)$,
equipped with the grading, given by  $\gamma=\begin{pmatrix} 1 &
0\\ 0 & -1
\end{pmatrix}$.

The representation $\rho$ of $\cA_\theta$ in $\cH$ is given by the
left multiplication, that is, for any $a\in \cA_\theta$
\[
\rho(a)=\begin{pmatrix} \lambda(a) & 0\\ 0 & \lambda (a)
\end{pmatrix},
\]
where the operator $\lambda(a)$ is given on $\cA_\theta\subset
L^2(\cA_\theta,\tau_0)$ by
\[
\lambda(a)b=ab, \quad b\in \cA_\theta.
\]
Introduce the differentiations $\delta_1$ and $\delta_2$ on the
algebra $\cA_\theta$ by
\[
\delta_1(U)=2\pi iU,\quad \delta_1(V)=0;\quad \delta_2(U)=0,\quad
\delta_2(V)=2\pi iV.
\]
The operator $D$ depends explicitly on $\tau$ and has the form
\[
D=\begin{pmatrix} 0 & \delta_1+\tau\delta_2\\
-\delta_1-\bar{\tau}\delta_2 & 0
\end{pmatrix}.
\]
The triples constructed above are two-dimensional smooth spectral
triples.
\end{ex}

Finally, we describe the simplest example of a spectral triple
associated with a closed manifold $M$, equipped with a Riemannian
foliation $\cF$. Fix a bundle-like metric $g_M$ on $M$. Let
$H=F^{\bot}$ be the orthogonal complement of $F=T{\mathcal F}$.

Define a triple $({\mathcal A}, {\mathcal H}, D)$ as follows:
\begin{enumerate}
  \item ${\mathcal A}=C^{\infty}_c(G)$;
  \item ${\mathcal H}=L^2(M,\Lambda^{*}H^{*})$ is the space
of transverse differential forms, endowed with the natural action
$R_{\Lambda^{*}H^{*}}$ of $\cA$;
  \item $D=d_H+d^*_H$ is the transverse signature operator (cf.
Section~\ref{s:do}).
\end{enumerate}

\begin{theorem}\cite{noncom}
\label{triple} The spectral triple $(A,{\mathcal H},D)$ is a
finite-dimensional spectral triple of dimension
$q=\operatorname{codim} \cF$.
\end{theorem}

More general examples of spectral triples given by transversally
elliptic operators on foliated manifolds are defined later in
Section~\ref{noncom}. Before we turn to the study of properties of
spectral triples associated with Riemannian foliations, we need to
have an appropriate pseudodifferential calculus, that is a subject
of Section~\ref{transpdo}.

\subsection{Geometry of para-Riemannian foliations}\label{para}
\begin{defn}
A foliation $\cF$ on a manifold $M$ is called para-Riemannian, if
there is an integrable distribution $V$ on $M$, which contains the
tangent bundle $T\cF$ to $\cF$ such that there are holonomy
invariant Euclidean structures in the fibers of the bundles $TM/V$
and $V/T\cF$.
\end{defn}

If $V$ is an integrable distribution on $M$, which defines a
para-Riemannian structure, and $\cV$ is the corresponding
foliation on $M$, $\cF\subset \cV$, then $\cV$ is Riemannian, and
the restriction of $\cF$ to every leaf of $\cV$ is a Riemannian
foliation.

An interest in studying of para-Riemannian foliations is related
with the fact that, in some cases, the study of arbitrary
foliations can be reduced to the case of para-Riemannian
foliations. A basic observation consists in the following. Let
$(M,\cF)$ be an arbitrary foliated manifold. Consider the bundle
$P$ over $M$, whose fiber $P_x$ at $x\in M$ is the space of all
Euclidean metrics on the vector space $\tau_x=T_xM/T_x\cF$. There
is a natural lift of $\cF$ to a foliation $\cV$ on $P$, moreover,
this foliation is para-Riemannian. For the first time, this
construction was used by Connes \cite{Co86} in order to extend the
map $K_0(C^\infty_c(G, |T\cG|^{1/2}))\to\CC$ given by the
transverse fundamental class of an arbitrary foliation $(M,\cF)$
to a map $K_0(C^*_r(G))\to\CC$ (cf. Section~\ref{s:tfclass}).
Connes constructed such an extension for $(P,\cV)$, using the
para-Riemannian condition, and then, using the Thom isomorphism,
he showed that the map $K_0(C^*_r(P,\cV))\to\CC$ given by the
transverse fundamental class of $(P,\cV)$ defines the desired
extension $K_0(C^*_r(G))\to\CC$ for the initial foliation
$(M,\cF)$.

In \cite{Hil-Skan}, Hilsum and Skandalis constructed a Fredholm
module associated with an arbitrary para-Riemannian foliation. To
do this, they made use of transversally hypoelliptic operators and
pseudodifferential operators of type $(\rho, \delta)$.

In \cite{Co-M}, Connes and Moscovici described a spectral triple
associated with a para-Riemannian foliations. More precisely, they
considered the closely related setting (strongly Morita
equivalent), when there is a manifold equipped with an action of a
discrete group, preserving a triangular structure.

Let us describe the construction of Connes and Moscovici. Let $W$
be an oriented smooth manifold, and $\Gamma $ a group of
diffeomorphisms of $W$. Consider the bundle $\pi:P(W)\to W$, whose
fiber $P_x(W)$ at $x\in W$ is the space of all Euclidean metrics
on the vector space $T_xW$. Thus, a point $p\in P(W)$ is given by
a point $x\in W$ and a nondegenerate quadratic form on $T_xW$. Let
$F_+(W)$ be the bundle of positive frames in the bundle $W$, whose
fiber $F_x(W)$ at $x\in W$ is the space of orientation preserving,
linear isomorphisms $\RR^n\to T_xW$. Equivalently, the bundle
$P(W)$ can be described as the quotient of the bundle $F_+(W)$ by
the fiberwise action of the subgroup $SO(n)\subset GL(n,\RR)$,
$P(W)=F_+(W)/SO(n)$. We will use the natural invariant Riemannian
metric on the symmetric space $GL_+(n,\RR)/O(n)$, given by the
matrix Hilbert-Schmidt norm on the tangent space to
$GL(n,\RR)/SO(n)$, which is identified with the space of symmetric
$n\times n$-matrices. If we transfer this metric to the fibers
$P_x$ of the fibration $P(W)=P$, we  get a Euclidean structure on
the vertical distribution $V\subset TP$. The normal space
$N_p=T_pP/V_p$ is naturally identified with the space $T_xW,
x=\pi(p)$. Thus, the quadratic form on $T_xW$, corresponding to
$p$, defines a natural Euclidean structure on $N_p$.

There are natural actions of the group $\Gamma$ on $F(W)$ and on
$P$. The action takes the fibers of the fibration $\pi:P\to W$ to
itself. Moreover, the Euclidean structures on the distributions
$V$ and $N$, introduced above, are invariant under the action of
$\Gamma$. In this case, one say that there is a triangular
structure on $P$, invariant under the group action.

Consider the Hermitian vector bundle $E=\wedge^*V^*_\CC\otimes
\wedge^*N_\CC^*$ over $P$. The Hermitian structure in the fibers
of $E$ is given by the Euclidean structures on $V$ and $N$. The
grading operators $\gamma_V$ and $\gamma_N$ in $\wedge^*V^*_\CC $
and $\wedge^*N_\CC^*$ are given by the Hodge operators of the
Euclidean structures and the orientations of the bundles $V$ and
$N$. The Euclidean structures on $V$ and $N$ also define a natural
volume form $v\in \wedge^*V^*\otimes \wedge^*N^*=\wedge T^*P$.

Let $\cA$ be the crossed product $C^\infty_c(P)\rtimes \Gamma$.
Recall that this algebra is generated as a linear space by the
expressions of the form $fU_g$, where $f\in C^\infty_c(P)$ and
$g\in \Gamma$ (cf. Example~\ref{ex:actions}).

Let $\cH$ be the space $L^2(P,E)$ of $L^2$ sections of the bundle
$E$, which is equipped by the Hilbert structure, given by the
volume form $v$ and the Hermitian structure  on $E$. The action of
$\cA$ in $\cH$ is given in a following way. A function $f\in
C^\infty_c(P)$ acts as the corresponding multiplication operator
in $\cH$. For any $g\in\Gamma$, the unitary operator $U_g$ in
$\cH$ is given by the natural actions of $\Gamma$ in the section
of the bundles $V$ and $N$.

Consider the foliation $\cV$, given by the fibers of the fibration
$P$. Then $V=T\cV$, $N$ is the normal bundle to $\cV$. Denote by
$d_L:C^\infty(P,E)\to C^\infty(P,E)$ the tangential de Rham
differential associated with the foliation $\cV$ (see
(\ref{e:d})). Let $Q_L$ be the second order tangential
differential operator in $C^\infty(P,E)$ given by
\[
Q_L=d_Ld^*_L-d^*_Ld_L.
\]
As shown in \cite{Co-M}, the principal symbol of $Q_L$ is
homotopic to the principal symbol of the signature operator
$d_L+d^*_L$.

Choose an arbitrary distribution $H$ on $P$, transverse to $V$.
Consider the corresponding transverse de Rham differential $d_H$
(see (\ref{e:d})) and the transverse signature operator
\[
Q_H=d_H+d^*_H.
\]
This operator depends on the choice of $H$, but its transverse
principal symbol is independent of $H$. Define a mixed signature
operator $Q$ in $C^\infty(P,E)$ as
\[
Q=Q_L(-1)^{\partial_N}+Q_N,
\]
where $(-1)^{\partial_N}$ denotes the parity operator in the
transverse direction, that is, it coincides with $1$ on
$\wedge^{\rm ev}N^*$ and with $-1$ on $\wedge^{\rm odd}N^*$. One
can show that $Q$ is essentially self-adjoint in $\cH$. Using the
functional calculus, define the operator $D$ as
\[
Q=D |D|.
\]

\begin{thm}\cite{Co-M}\label{t:para}
The triple $(\cA,\cH,D)$ is a spectral triple of dimension $\dim
V+2\dim N$.
\end{thm}

The proof of this theorem makes an essential use of the
pseudodifferential calculus on Heisenberg manifolds, constructed
by Beals and Greiner \cite{BG}. Note that the above construction
can be applied to any manifold endowed with a triangular
structure, invariant under a group action. In this case, a
complicated analytic problem is the essential self-adjointness of
$Q$ in the case of a noncompact manifold.

\section{Transverse pseudodifferential calculus}\label{transpdo}
Throughout in this Section, we will consider a closed, connected,
oriented, foliated manifold $(M,{\mathcal F})$, $\dim M = n$,
$\dim {\mathcal F} = p$, $p + q = n$, and a Hermitian vector
bundle $E$ on $M$ of rank $r$. We will consider operators, acting
on half-densities. We will denote by $C^{\infty}(M,E)$ the space
of smooth sections, by $L^2(M,E)$ the Hilbert space of square
integrable sections, by ${\cD}'(M,E)$ the space of distributional
sections, ${\cD}'(M,E)=C^{\infty}(M,E)'$, and by $H^s(M,E)$ the
Sobolev space of order $s$ for the vector bundle $E\otimes
|TM|^{1/2}$.

\subsection{Classes $\Psi^{m,-\infty}(M,{\mathcal F}, E)$} \cite{noncom}
Consider the $n$-dimensional cube $I^{n} = I^{p} \times I^{q}$
equipped with the trivial foliation, whose leaves are $I^{p}\times
\{y\}$, $y\in I^{q}$. As usual, the coordinates in $I^{n}$ are
denoted by $(x,y)$, $x\in I^{p}$, $y \in I^{q}$, and the dual
coordinates by $(\xi ,\eta )$, $\xi \in {\RR}^{p}$, $\eta \in
{\RR}^{q}$.

\begin{defn}
A function $k\in C^{\infty}(I^{p} \times I^{p} \times I^{q} \times
{\mathbb R}^{q},{\mathcal L}({\mathbb C}^r))$ belongs to the class
$S^{m}(I^{p}\times I^{p}\times I^{q} \times {\mathbb R}^{q},
{\mathcal L}({\mathbb C}^r))$, if, for any multiindices
$\alpha\in\NN^q $ and $\beta\in\NN^{2p+q}$, there is a constant
$C_{\alpha \beta} > 0$ such that
\begin{multline*}
\|
\partial^{\alpha}_{\eta}
\partial^{\beta}_{(x,x',y)}k(x,x',y,\eta )\| \leq C_{\alpha \beta
}(1 +\vert \eta \vert )^{ m-\vert \alpha \vert },\\ (x,x',y)\in
I^{p}\times I^p\times I^q,\quad \eta \in {\mathbb R}^{q}.
\end{multline*}
\end{defn}

In the following, we will only consider classical symbols.

\begin{defn}
A function $k\in C^{\infty}(I^{p}\times I^{p}\times I^{q} \times
{\mathbb R}^{q}, {\mathcal L}({\mathbb C}^r))$ is called a
classical symbol of order $z\in {\mathbb C}$, if it can be
represented as an asymptotic sum  $$ k(x,x',y,\eta)\sim
\sum_{j=0}^{\infty} \theta(\eta) k_{z-j}(x,x',y,\eta), $$ where
$k_{z-j}\in C^{\infty}(I^{p}\times I^{p}\times I^{q} \times
({\mathbb R}^{q}\backslash \{0\}), {\mathcal L}({\mathbb C}^r))$
are homogeneous in $\eta$ of degree $z-j$, that is, $$
k_{z-j}(x,x',y,t\eta)=t^{z-j}k_{z-j}(x,x',y,\eta), \quad t>0, $$
and $\theta$ is a smooth function in ${\mathbb R}^{q}$ such that
$\theta(\eta)=0$ for $|\eta|\leq 1$, $\theta(\eta)=1$ for
$|\eta|\geq 2$.
\end{defn}

In this definition, the sign of asymptotic summation $\sim$ means
that $k - \sum_{j=0}^{N} \theta k_{z-j} \in S^{{\rm Re}\,z-N-1}$
for any natural $N$.

A symbol $k \in S ^{m} (I ^{p} \times I^p\times I^q\times {\mathbb
R}^{q}, {\mathcal L}({\mathbb C}^r))$ defines an operator
\[
A: C^{\infty}_c(I^n, {\mathbb C}^r) \rightarrow C^{\infty}(I^n,
{\mathbb C}^r)
\]
as
\begin{equation}
\label{pdo:loc} Au(x,y)=(2\pi)^{-q} \int
e^{i(y-y')\eta}k(x,x',y,\eta) u(x',y') \,dx'\,dy'\,d\eta,
\end{equation}
where \(u \in C^{\infty}_{c}(I^{n}, {\mathbb C}^r), x \in I^{p}, y
\in I^{q}\). Denote by $\Psi^{m,-\infty}(I^n,I^p,{\mathbb C}^r)$
the class of operators of the form (\ref{pdo:loc}) with $k \in S
^{m} (I ^{p} \times I^p\times I^q\times {\mathbb R}^{q}, {\mathcal
L}({\mathbb C}^r))$ such that its Schwartz kernel is compactly
supported in $I^n\times I^n$.

If $\phi: U \subset M\to I^p\times I^q, \phi': U' \subset M \to
I^p\times I^q$ are compatible foliated charts on $M$ endowed with
trivializations of $E$, then an operator
$A\in\Psi^{m,-\infty}(I^n,I^p,{\mathbb C}^r)$ defines an operator
$A': C^{\infty}_c(U,E)\rightarrow C^{\infty}_c(U',E)$, which can
extended in a trivial way to an operator in $C^{\infty}(M,E)$. The
operator obtained in such a way  will be also denoted by $A'$ and
called an elementary operator of class $\Psi^{m,-\infty}
(M,{\mathcal F},E)$.

\begin{defn}
The class $\Psi ^{m,-\infty}(M,{\mathcal F},E)$ consists of
operators $A$, acting in $C^{\infty}(M,E)$, which can be
represented in the form $A=\sum_{i=1}^k A_i + K$, where $A_i$ are
elementary operators of class $\Psi ^{m,-\infty}(M,{\mathcal
F},E)$, corresponding to pairs $\phi_i,\phi'_i$ of compatible
foliated charts, $K\in \Psi ^{-\infty}(M,E)$.
\end{defn}

Operators from $\Psi ^{m,-\infty}(M,{\mathcal F},E)$ are called
transversal pseudodifferential operators. These operators can be
represented as Fourier integral operators associated with a
natural canonical relation in the punctured cotangent bundle
$\widetilde{T}^*M =T^{*}M \backslash \{0\}$. More precisely, let
${\mathcal F}_N$ be the linearized foliation in
$\widetilde{N}^*{\mathcal F}$, $G_{{\mathcal F}_N}$ the holonomy
groupoid of the foliation ${\mathcal F}_N$. The map
\begin{equation}
\label{can} (r,s):G_{{\mathcal F}_N}\rightarrow
\widetilde{T}^*M\times \widetilde{T}^*M
\end{equation}
defines an injectively immersed relation in $\widetilde{T}^*M$.
The algebra of Fourier integral operators associated with the
canonical relation~(\ref{can}) coincides with
$\Psi^{*,-\infty}(M,{\mathcal F},E)$. Observe also that the class
$\Psi^{*,-\infty}(M,{\mathcal F}, E)$ coincides with the algebra
${\mathcal R}_{\Sigma}$, associated with the coisotropic conic
submanifold $\Sigma=\widetilde{N}^*\cF$ in $\widetilde{T}^*M$,
which was introduced by Guillemin and Sternberg in \cite{GS79}.

Recall that a Fourier integral operator on $M$ is a linear
operator $F:C^{\infty}(M)\to {\cD}'(M)$, microlocally
representable in the form
\begin{equation}\label{e:FIO}
Fu(x)=\int e^{\phi(x,y,\theta)}a(x,y,\theta)\,u(y)\,dy\,d\theta,
\end{equation}
where $x\in X\subset \RR^n, y\in Y\subset\RR^n, \theta\in
\RR^N\setminus 0$. Here $a(x,y,\theta)\in S^m(X\times
Y\times\RR^N)$ is an amplitude, $\phi$ is a nondegenerate phase
function. (Concerning to Fourier integral operators cf., for
instance, \cite{H4,Taylor,Treves2})

Consider the smooth map from $X\times Y\times\RR^N$ to $T^*X\times
T^*Y$ given by
\[
(x,y,\theta)\mapsto (x,\phi_x(x,y,\theta),y,-\phi_y(x,y,\theta)).
\]
The image of the set
\[
\Sigma_\phi=\{(x,y,\theta)\in X\times Y\times\RR^N
:\phi_\theta(x,y,\theta)=0\}
\]
under this map turns out to be a homogeneous canonical relation
$\Lambda_\phi$ in $T^*X\times T^*Y$. (Recall that a closed conic
submanifold $C\subset T^*(X\times Y)\setminus 0$ is called a
homogeneous canonical relation, if it is contained in
$\tilde{T}^*X\times \tilde{T}^*Y$ and is Lagrangian with respect
to the 2-form $\omega_X-\omega_Y$, where $\omega_X, \omega_Y$ are
the canonical symplectic forms on $T^*X, T^*Y$ respectively.)

A Fourier integral operator $F$ is said to be associated with
$\Lambda_\phi$. Let us write $F\in I^m(X\times Y, \Lambda_\phi)$,
if $a\in S^{m+n/2-N/2}(X\times Y\times\RR^N)$.

Let $\phi: U \subset M \to I^p\times I^q, \phi': U' \subset M \to
I^p\times I^q$ be compatible foliated charts on $M$, and
$A:C^{\infty}_c(U,\left. E\right|_U)\to C^{\infty}_c(U',\left.
E\right|_{U'})$ an elementary operator given by (\ref{pdo:loc})
with $k \in S ^{m} (I ^{p} \times I^p\times I^q\times {\RR}^{q},
{\cL}({\CC}^r))$. It can be represented in the form (\ref{e:FIO}),
if one takes $X=U$ with the coordinates $(x,y)$, $Y=U'$ with the
coordinates $(x',y')$, $\theta=\eta $, $N=q$, the phase function
$\phi(x,y,x',y')=(y-y')\eta$ and the amplitude $a=k(x,x',y,\eta)$.
The associated homogeneous canonical relation $\Lambda_\phi$ can
be described as
\[
\Lambda_\phi=\{(x,y,\xi,\eta,x',y',\xi',\eta')\in T^*U\times T^*U'
: y=y',\xi=\xi'=0,\eta=-\eta'\},
\]
that coincides with the intersection of $G'_{{\mathcal F}_N}$ with
$T^*U\times T^*U'$. (For any relation $C\subset T^*X\times T^*Y$,
$C'$ denotes the image of $C$ under the map $T^*X\times T^*Y \to
T^*X\times T^*Y : (x,\xi,y,\eta)\mapsto (x,\xi,y,-\eta)$).
Furthermore, note that the class $\Psi^{m,-\infty}(M, {\mathcal
F},E)$ consists of all operators in $C^{\infty}(M,E)$ with the
Schwartz kernels from the class $I^{m-p/2} (M\times M,
G'_{{\mathcal F}_N}; \cL(E)\otimes |T(M\times M)|^{1/2})$. Since
$G_{{\mathcal F}_N}$ is, in general, a (one-to-one) immersed
canonical relation, it is necessary to be more precise with the
definition of the classes $I^{m}(M\times M,G'_{{\mathcal
F}_N};\cL(E)\otimes |T(M\times M)|^{1/2})$. This can be done by
analogy with the definition of the classes of leafwise
pseudodifferential operators on a foliated manifold given in
\cite{Co79} (see also the definition of the classes
$\Psi^{*,-\infty}(M,\cF,E)$ given above). Namely, the space
$I^m(M\times M, G'_{{\mathcal F}_N})$ of compactly supported
Lagrangian distributions is defined as the set of finite sums of
elementary Lagrangian distributions, corresponding to pairs of
compatible foliated charts on $M$.

\subsection{Symbolic calculus in $\Psi^{m,-\infty}$}
\label{symb:psi} Symbolic properties of the classes $\Psi
^{m,-\infty}(M,{\mathcal F},E)$ can be obtained as a consequence
of the corresponding facts for the Guillemin-Sternberg algebras
${\mathcal R}_{\Sigma}$ \cite{GS79}, but in many cases it is
simpler to give a direct proof (cf. \cite{noncom}).

The principal symbol $\sigma_A$ of an elementary operator
$A\in\Psi^{m,-\infty}(I^n,I^p,{\mathbb C}^r)$ given by
(\ref{pdo:loc}) is defined to be the matrix-valued half-density
$\sigma_A$ on $I^p\times I^p\times I^q\times ({\mathbb R}^q
\backslash \{0\})$ given by
\begin{multline}
\label{princ} \sigma_A(x,x',y,\eta) = k_m(x,x',y,\eta)
|dx\,dx'|^{1/2},\\ (x,x',y,\eta) \in I^p\times I^p\times I^q\times
({\mathbb R}^q \backslash \{0\}),
\end{multline}
where $k_m$ is the homogeneous of degree $m$ component of $k$.

Before we give the global description of the principal symbol for
operators of class $\Psi^{m,-\infty} (M, {\mathcal F},E)$, we
introduce some auxiliary notions. Denote by $\pi^*E$ the lift of
the vector bundle $E$ to $\tilde{N}^*\cF$ under the map
$\pi:\tilde{N}^*\cF\to M$. Since $\tilde{N}^*\cF$ is noncompact,
it is impossible to define the structure of an involutive algebra
on the whole space $C^{\infty}(G_{{\mathcal F}_N}, {\mathcal
L}(\pi^*E)\otimes |T{\mathcal G}_N|^{1/2})$. Introduce the space
$C^{\infty}_{prop}(G_{{\mathcal F}_N}, {\mathcal L}(\pi^*E)\otimes
|T{\mathcal G}_N|^{1/2})$, which consists of all properly
supported elements $k\in C^{\infty}(G_{{\mathcal F}_N}, {\mathcal
L}(\pi^*E)\otimes |T{\mathcal G}_N|^{1/2})$ (this means that the
restriction of $r:G_{\cF_N}\to \tilde{N}^*\cF$ to $\supp k$ is a
proper map). The structure of an involutive algebra on
$C^{\infty}_{prop}(G_{{\mathcal F}_N}, {\mathcal L}(\pi^*E)\otimes
|T{\mathcal G}_N|^{1/2})$ is defined, using the standards formulas
(cf. Section~\ref{s:calgebra}).

The space of all sections $s\in C^{\infty}_{prop}(G_{{\mathcal
F}_N},{\mathcal L}(\pi^*E)\otimes |T{\mathcal G}_N|^{1/2})$,
homogeneous of degree $m$ with respect to the fiberwise
multiplication in the fibers of the bundle $\pi:\tilde{N}^*\cF\to
M$, is denoted by $S^{m}(G_{{\mathcal F}_N},{\mathcal
L}(\pi^*E)\otimes |T{\mathcal G}_N|^{1/2})$. The space $$
S^*(G_{{\mathcal F}_N},{\mathcal L}(\pi^*E)\otimes |T{\mathcal
G}_N|^{1/2})=\bigcup_m S^m(G_{{\mathcal F}_N},{\mathcal
L}(\pi^*E)\otimes |T{\mathcal G}_N|^{1/2}) $$ is a subalgebra of
$C^{\infty}_{prop}(G_{{\mathcal F}_N}, {\mathcal L}(\pi^*E)\otimes
|T{\mathcal G}_N|^{1/2})$.

Let $\phi: U\subset M\to I^p\times I^q, \phi': U' \subset M \to
I^p\times I^q$ be two compatible foliated charts on $M$ endowed
with trivializations of $E$. Then the corresponding coordinate
charts $\phi_n: U_1\subset N^*\cF\to I^p\times I^q \times \RR^q,
\phi_n': U_1^\prime \subset N^*\cF \to I^p\times I^q \times \RR^q$
(see Section~\ref{s:fol}) are compatible foliated charts on the
foliated manifold $(N^*\cF,\cF_N)$ endowed with obvious
trivializations of $\pi^*E$. Thus, there is a foliated chart
$\Gamma_N: W(\phi_n,\phi'_n) \subset G_{{\mathcal F}_N}\to
I^p\times I^p\times I^q\times {\mathbb R}^q$ on the foliated
manifold $(G_{\cF_N},{\cG}_N)$. The local half-density defined by
(\ref{princ}) in any foliated chart $W(\phi_n,\phi'_n)$ gives a
well-defined element $\sigma_A$ of $S^m(G_{{\mathcal
F}_N},{\mathcal L}(\pi^*E)\otimes |T{\mathcal G}_N|^{1/2})$ ---
the principal symbol of $A\in \Psi^{m,-\infty} (M, {\mathcal
F},E)$.

If we fix a smooth positive leafwise density, then $\sigma_A$ is
identified with an element of $C^{\infty}(G_{{\mathcal
F}_N},{\mathcal L}(\pi^*E))$. The notion of the principal symbol
of an operator of class $\Psi^{m,-\infty}(M,{\mathcal F},E)$ as a
Fourier integral operator (cf., for instance, \cite{H4}) agrees
with this definition of the principal symbol.

\begin{proposition}\cite{noncom}
\label{prop} The principal symbol map
\[
\sigma: \Psi^{*,-\infty}(M,{\mathcal F},E)\rightarrow
S^*(G_{{\mathcal F}_N},{\mathcal L}(\pi^*E)\otimes |T{\mathcal
G}_N|^{1/2})
\]
is a $\ast$-homomorphism of involutive algebras. Thus, if $A\in
\Psi ^{m_{1},-\infty}(M,{\mathcal F},E)$ and $B\in \Psi
^{m_{2},-\infty}(M,{\mathcal F},E)$, then  \(C = AB\) belongs
\(\Psi ^{m_{1}+m_{2},-\infty}(M,{\mathcal F},E)\) and
$\sigma_{AB}=\sigma_A\sigma_B$. If $A\in \Psi ^{m,-\infty}
(M,{\mathcal F},E)$, then \(A^{*} \in \Psi
^{m,-\infty}(M,{\mathcal F},E)\) and $\sigma_{A^*}= (\sigma_A)^*$.
\end{proposition}

Recall that the transversal principal symbol $\sigma _{P}$ of an
operator $P\in \Psi ^{m}(M,E)$ is the restriction of its principal
symbol $p_m$ to $\widetilde{N}^{\ast}{\mathcal F}$.

\begin{proposition}
\label{module} If $A\in  \Psi ^{\mu}(M,E)$ and $B\in  \Psi
^{m,-\infty}(M,{\mathcal F},E)$, then $AB$ and $BA$ belong to
$\Psi ^{\mu+m,-\infty}(M,{\mathcal F},E)$ and
\begin{align*}
\sigma_{AB}(\gamma,\eta)&=\sigma_A(\eta)\sigma_B(\gamma,\eta),
\quad (\gamma,\eta)\in G_{{\mathcal F}_N},\\[6pt]
\sigma_{BA}(\gamma,\eta)&=
\sigma_B(\gamma,\eta)\sigma_A(dh_{\gamma}^*(\eta)), \quad
(\gamma,\eta)\in G_{{\mathcal F}_N}.
\end{align*}
\end{proposition}

\begin{proposition}\label{as:p:vanishes}
If the transversal principal symbol of an operator $P\in
\Psi^{\mu}(M,E)$ vanishes, then, for any $K\in
\Psi^{m,-\infty}(M,\cF,E)$, the operators $K P$ and $P K$ belong
to $\Psi^{m+\mu-1,-\infty} (M,{\mathcal F},E)$.
\end{proposition}

Suppose that $E$ is holonomy equivariant, that is, there is an
action $T(\gamma):E_x\rightarrow E_y, \gamma \in G,
\gamma:x\rightarrow y$ of the holonomy groupoid $G$ in the fibers
of $E$. Then the bundle ${\mathcal L}(\pi^*E)$ on $N^*\cF$ is
holonomy equivariant. Denote by  $\operatorname{ad} T$ the
corresponding action of the holonomy groupoid $G_{\cF_N}$ in the
fibers of ${\mathcal L}(\pi^*E)$.

\begin{defn}\label{d:hinv}
The transversal principal symbol $\sigma_P$ of an operator
$P\in\Psi^m(M,E)$ is holonomy invariant, if, for any leafwise path
$\gamma $ from $x$ to $y$ and for any $\nu \in N^{
\ast}_{y}{\mathcal F}$, the following identity holds: $$ {\rm
ad}\,T(\gamma,\nu)[\sigma _{P}(dh_{\gamma}^{\ast}(\nu ))] = \sigma
_{P}(\nu ). $$
\end{defn}

Let $p_{m} \in S^{m}(I^{n} \times ({\RR}^{n}\setminus\{0\}))$ be
the principal symbol of $P\in \Psi ^{m}(M)$ in some foliated
chart, then its transversal principal symbol $\sigma _{P}$ is
given by $$ \sigma _{P}(x,y,\eta ) = p_{m}(x,y,0,\eta ), \quad
(x,y)\in I^{n}, \quad \eta  \in {\RR}^{q}\setminus \{0\}, $$ and
the holonomy invariance of $\sigma _{P}$ means that $\sigma_P$ is
independent of $x$.

The assumption of the existence of a positive order
pseudodifferential operator with the holonomy invariant
transversal principal symbol on a foliated manifold imposes
sufficiently strong restrictions on the geometry of the foliation.
An example of an operator with the holonomy invariant transverse
principal symbols is the transverse signature operator on a
Riemannian foliation.

There is a canonical embedding $$ i:C^{\infty}_{prop}(G_{{\mathcal
F}_N},|T{\mathcal G}_N|^{1/2})\hookrightarrow
C^{\infty}_{prop}(G_{{\mathcal F}_N},{\mathcal L}(\pi^*E)\otimes
|T{\mathcal G}_N|^{1/2}),$$ which takes any $k\in
C^{\infty}_{prop}(G_{{\mathcal F}_N},|T{\mathcal G}_N|^{1/2})$ to
$i(k)=k \operatorname{ad} T$. We will identify
$C^{\infty}_{prop}(G_{{\mathcal F}_N},|T{\mathcal G}_N|^{1/2})$
with its image in $C^{\infty}_{prop}(G_{{\mathcal F}_N},{\mathcal
L}(\pi^*E)\otimes |T{\mathcal G}_N|^{1/2})$ under the map $i$.

\begin{defn}
An operator $P\in \Psi^{m,-\infty}(M,{\mathcal F},E)$ is said to
have the scalar principal symbol, if its principal symbol belongs
to $C^{\infty}_{prop}(G_{{\mathcal F}_N},|T{\mathcal
G}_N|^{1/2})$.
\end{defn}

Denote by $\Psi_{sc}^{m,-\infty}(M,{\mathcal F},E)$ the set of all
operators $K\in \Psi^{m,-\infty}(M,{\mathcal F},E)$ with the
scalar principal symbol. Observe that, for any $k\in
C^{\infty}_c(G,|T{\cG}|^{1/2})$, the operator $R_E(k)$ belongs to
$\Psi^{0,-\infty}_{sc}(M,\cF,E)$, and
\[
\sigma(R_E(k))=\pi_G^*k\in C^{\infty}_{prop}(G_{{\mathcal
F}_N},|T{\mathcal G}_N|^{1/2}),
\]
where $\pi_G :G_{\cF_N}\to G$ is the natural map.

\begin{proposition}
Let $(M,\cF)$ be a compact foliated manifold and $E$ a holonomy
equivariant vector bundle. If
$A\in\Psi_{sc}^{m,-\infty}(M,{\mathcal F},E)$, and $P\in\Psi
^{\mu}(M,E)$ has the holonomy invariant transversal principal
symbol, then $[A,P]\in\Psi^{m+\mu-1,-\infty}(M,{\mathcal F})$.
\end{proposition}

Finally, one has the following statement on the continuity of the
symbol map. Any $A\in \Psi^{0,-\infty}(M,{\mathcal F},E)$ defines
a bounded operator in the Hilbert space $L^2(M,E)$. Denote by
$\bar{\Psi}^{0,-\infty}(M,{\mathcal F},E)$ the closure of
$\Psi^{0,-\infty}(M,{\mathcal F},E)$ in the uniform topology of
${\mathcal L}(L^2(M,E))$.

Recall that, for any $\nu\in \tilde{N}^*\cF$, there is a natural
$\ast$-representation $R_\nu$ of the algebra $S^{0}(G_{{\mathcal
F}_N}, {\mathcal L}(\pi^*E)\otimes |T{\mathcal G}_N|^{1/2})$ in
$L^2(G^\nu_{{\cF}_N}, s_N^*(\pi^*E))$ (cf. (\ref{e:Rx})). Thus,
for any $k \in S^{0}(G_{{\mathcal F}_N}, {\mathcal
L}(\pi^*E)\otimes |T{\mathcal G}_N|^{1/2})$, the continuous
operator family $\{R_\nu(k)\in \cL(L^2(G^\nu_{{\cF}_N},
s_N^*(\pi^*E))):\nu\in \tilde{N}^*\cF\}$ defines a bounded
operator in $L^2(G_{{\cF}_N},s_N^*(\pi^*E))$. We will identify an
element $k\in S^{0}(G_{{\mathcal F}_N}, {\mathcal
L}(\pi^*E)\otimes |T{\mathcal G}_N|^{1/2})$ with the corresponding
bounded operator in $L^2(G_{{\cF}_N},s_N^*(\pi^*E))$ and denote by
$\bar{S}^{0}(G_{{\mathcal F}_N}, {\mathcal L}(\pi^*E)\otimes
|T{\mathcal G}_N|^{1/2})$ the closure of $S^{0}(G_{{\mathcal
F}_N}, {\mathcal L}(\pi^*E)\otimes |T{\mathcal G}_N|^{1/2})$ in
the uniform topology of $\cL(L^2(G_{{\cF}_N},s_N^*(\pi^*E)))$.

\begin{prop}[\cite{egorgeo}]
\label{seq} $($1$)$ The symbol map $$ \sigma :
\Psi^{0,-\infty}(M,{\mathcal F},E)\rightarrow S^{0}(G_{{\mathcal
F}_N},{\mathcal L}(\pi^*E)\otimes |T{\mathcal G}_N|^{1/2}) $$
extends by continuity to a homomorphism $$
\bar{\sigma}:\bar{\Psi}^{0,-\infty}(M,{\mathcal F},E) \rightarrow
\bar{S}^{0}(G_{{\mathcal F}_N},{\mathcal L}(\pi^*E)\otimes
|T{\mathcal G}_N|^{1/2}). $$

$($2$)$ The ideal $\Ker \bar{\sigma}$ contains the ideal of
compact operators in $L^2(M,E)$.
\end{prop}

\subsection{Zeta-function of transversally elliptic operators}
\label{s:zeta} Let $(M,{\mathcal F})$ be a compact foliated
manifold and $E$ a Hermitian vector bundle on $M$. Suppose that an
operator $A\in\Psi^{m}(M,E)$ satisfies the following conditions:
\begin{description}
  \item[(T1)] $A$ is a transversally elliptic operator with the positive transversal
principal symbol;
  \item[(T2)] the operator $A$, considered as an unbounded operator
in the Hilbert space $L^2(M,E)$, is essentially self-adjoint on
the initial domain $C^{\infty}(M,E)$, and its closure is an
invertible and positive operator.
\end{description}

The condition (T2) can be considered as an invariance type
condition, which is usually assumed for transversally elliptic
operators.

Before we formulate a general result on the meromorphic
continuation of the zeta-function of an operator $A$, satisfying
the conditions (T1) and (T2), we define a Wodzicki-Guillemin type
residue trace for operators from $\Psi^{m,-\infty}(M,{\mathcal
F},E)$. It suffices to do this for elementary operators of class
$\Psi^{m,-\infty}(M,{\mathcal F},E)$. For an operator
$P\in\Psi^{m,-\infty}(I^n,I^p,{\mathbb C}^r)$, define the residue
form $\rho_P$ as $$ \rho_P = \left(\int_{|\eta|=1} \Tr
k_{-q}(x,x,y,\eta)\,d\eta\right) |dx dy|, $$ and the residue trace
$\tau(P)$ as $$ \tau(P)=(2\pi)^{-q}\int_{|\eta|=1}\Tr
k_{-q}(x,x,y,\eta)\, dx dy d\eta, $$ where $k_{-q}$ is the
homogeneous of degree $-q$ component of the complete symbol $k$ of
 $P$.

For any $P\in \Psi^{m,-\infty}(M,{\mathcal F},E)$, its residue
form $\rho_P$ is a well-defined density on $M$, and the residue
trace $\tau(P)$ is obtained by the integration of $\rho_P$ over
$M$: $$ \tau(P)=(2\pi)^{-q} \int_{M}\rho_P. $$

It should be noted that these results are particular cases of the
results of \cite{Gu_1,Gu_2}, concerning to Fourier integral
operators, but one can get them directly (cf. \cite{noncom}).

\begin{theorem}
\label{zeta1} Let an operator $A\in\Psi^{m}(M,E)$ satisfy the
conditions (T1) and (T2), and $Q\in \Psi^{l,-\infty}(M,{\mathcal
F},E)$, $l\in {\mathbb Z}$. Then the function $z\mapsto
\tr(QA^{-z})$ is holomorphic for  ${\rm Re}\, z> l+q/m$ and admits
a (unique) meromorphic extension to ${\mathbb C}$ with at most
simple poles at $z_k=k/m$ with integer $k\leq l+q$. Its residue at
the point $z=z_k$ equals
\[
\underset{z=z_k}{\res} \tr (QA^{-z})=q \tau(QA^{-k/m}).
\]
\end{theorem}

Now suppose that $E$ is holonomy equivariant. For any $z\in \CC,
{\rm Re}\, z>q/m$, and for any $k\in C^\infty_c(G)$, the operator
$R_E(k)A^{-z}$ is a trace class operator, and one can define the
distributional zeta-function of $A$ as
\begin{equation}
\label{zeta} \langle\zeta_{A}(z), k\rangle =\tr R_E(k)A^{-z},\quad
\text{Re}\, z> q/m, \quad k\in C^\infty_c(G),
\end{equation}
which is an analytic function in the half-plane ${\rm Re}\,
z>q/m$.

\begin{theorem}
\label{zeta0} For any $k\in C^{\infty}_c(G)$, the function $
\langle\zeta_{A}(z), k\rangle $ extends to a meromorphic function
on the complex plane with at most simple poles at the points
$z=q/m,(q-1)/m,\ldots.$
\end{theorem}

\subsection{Egorov theorem}
\label{s:Egorov} The classical Egorov theorem \cite{egorov} is one
of fundamental results of microlocal analysis, which relates the
quantum evolution of pseudodifferential operators with the
classical dynamics of principal symbols. Recall the formulation of
this theorem. Let $M$ be a compact manifold, $E$ a vector bundle
on $M$ and $P\in \Psi^1(M,E)$ a positive self-adjoint
pseudodifferential operator with the positive principal symbol
$p$. The Egorov theorem states that, if $A\in \Psi^0(M,E)$, then
$A(t)=e^{itP}Ae^{-itP}\in \Psi^0(M,E)$. Moreover, if $E$ is the
trivial line bundle and $a\in S^0(\tilde{T}^*M)$ is the principal
symbol of $A$, then the principal symbol  $a_t\in
S^0(\tilde{T}^*M)$ of $A(t)$ is given by
\[
a_t(x,\xi)=a(f_t(x,\xi)), \quad (x,\xi)\in \tilde{T}^*M,
\]
where  $f_t$ is the bicharacteristic flow of $P$, that is, the
Hamiltonian flow on $T^*M$, defined by its principal symbol as the
Hamiltonian.

In \cite{egorgeo} a version of the Egorov theorem for
transversally elliptic operators on compact foliated manifolds is
proved. Suppose that $(M,{\mathcal F})$ is a compact foliated
manifold. Let $A$ be a linear operator in $C^{\infty}(M)$,
satisfying the following conditions:
\begin{description}
  \item[(A1)] $A\in \Psi^2(M,|TM|^{1/2})$ is a transversally elliptic
operator with the scalar principal symbol and the positive,
holonomy invariant transversal principal symbol;
\item[(A2)] $A$ is an essentially self-adjoint positive operator
in $L^2(M)$ (with the initial domain $C^{\infty}(M)$).
\end{description}

We start with the definition of the transverse bicharacteristic
flow of the operator $\sqrt{A}$. Let $a_2\in S^2(\tilde{T}^*M)$ be
the principal symbol of $A$. Take any function ${\tilde p}\in
S^1(\tilde{T}^*M)$, which coincides with $\sqrt{a_2}$ in some
conic neighborhood of $\tilde{N}^*{\mathcal F}$. Denote by
$X_{\tilde p}$ the Hamiltonian vector field on $T^*M$ with the
Hamiltonian $\tilde p$. For any $\nu\in \tilde{N}^*{\mathcal F}$,
the vector $X_{\tilde p}(\nu)$ is tangent to $\tilde{N}^*{\mathcal
F}$. Therefore, the Hamiltonian flow $\tilde{f}_t$ with the
Hamiltonian $\tilde{p}$ preserves $\tilde{N}^*{\mathcal F}$.
Denote by $f_t$ its restriction to $N^*{\mathcal F}$. One can show
that the vector field $X_{\tilde p}$ on $\tilde{N}^*{\mathcal F}$
is an infinitesimal transformation of the foliation $\cF_N$, and,
therefore, the flow $f_t$ preserves the foliation $\cF_N$.

Using the fact that $X_{\tilde p}$ is an infinitesimal
transformation of $\cF_N$, it is not difficult to show the
existence of a unique vector field $\cH_p$ on $G_{{\mathcal F}_N}$
such that $ds_N(\cH_p)=X_{\tilde p}$ and $dr_N(\cH_p)=X_{\tilde
p}$. Let $F_t$ be the flow on $G_{{\mathcal F}_N}$ defined by
$\cH_p$. It is easy to see that $s_N\circ F_t=f_t\circ s_N$,
$r_N\circ F_t=f_t\circ r_N$ and the flow $F_t$ preserves $\cG_N$.

\begin{defn}\label{d:trans}
Let $A$ be a linear operator in $C^{\infty}(M)$ satisfying the
conditions (A1) and (A2). The transverse bicharacteristic flow of
the operator $\sqrt{A}$ is the one-parameter group $F^*_t$ of
automorphisms of the involutive algebra
$C^{\infty}_{prop}(G_{{\mathcal F}_N},|T{\mathcal G}_N|^{1/2})$,
induced by the action of $F_t$.
\end{defn}

It is easily seen that the definition of the transverse
bicharacteristic flow is independent of the choice of $\tilde{p}$.

The construction of the transverse bicharacteristic flow is an
example of noncommutative symplectic (or, maybe, better to say,
Poisson) reduction. Here we mean by symplectic reduction the
following procedure \cite[Chapter III, Section 14]{LM87} (see also
\cite{Li75,Li77}).

Let $(X,\omega)$ be a symplectic manifold and $Y$ a submanifold of
$X$ such that the $2$-form $\omega_Y$ on $Y$, induced by the
symplectic form $\omega$, has constant rank. Let ${\cF}_Y$ be the
characteristic foliation of $Y$ with respect to the form
$\omega_Y$. If ${\cF}_Y$ is given by the fibers of a submersion
$p: Y \to B$, then there is a unique symplectic form $\omega_B$ on
$B$ such that $p^*\omega_B=\omega_Y$. The symplectic manifold $(B,
\omega_B)$ is called the reduced symplectic manifold associated
with $Y$. In the particular case when $Y$ is the pre-image of a
point under the momentum map defined by a Hamiltonian action of a
Lie group, the symplectic reduction associated with $Y$ is a
well-known procedure of a symplectic reduction developed by
Marsden and Weinstein \cite{MW}.

If $Y$ is invariant under the action of a Hamiltonian flow with a
Hamiltonian $H\in C^\infty(X)$ (this condition is equivalent to
the condition that the restriction $(dH)\left|_Y\right.$ is
constant along the leaves of the characteristic foliation
${\cF}_Y$), then there is a unique function $\hat{H}\in
C^\infty(B)$, called the reduced Hamiltonian, such that
$H\left|_Y\right. = \hat{H}\circ p$ (cf., for instance,
\cite[Chapter III, Theorem 14.6]{LM87}). Moreover, the map $p$
takes the restriction of the Hamiltonian flow $H$ on $Y$ to the
reduced Hamiltonian flow on $B$, defined by the reduced
Hamiltonian $\hat{H}$.

Now let $(M,{\cF})$ be a foliated manifold. Consider the
symplectic reduction given by the coisotropic submanifold
$Y=N^*{\cF}$ of the symplectic manifold $X=T^*M$. The
corresponding characteristic foliation ${\cF}_Y$ coincides with
the linearized foliation ${\cF}_N$. The algebra
$C^{\infty}_{prop}(G_{{\cF}_N},|T{\cG}_N|^{1/2})$ plays a role of
a noncommutative analogue of the algebra of smooth functions on
the base $B=N^*{\cF}/{\cF}_N$ - ``the cotangent bundle'' to
$M/\cF$. Under the condition (A1) and (A2), the above construction
can be considered as a noncommutative analogue of the symplectic
reduction procedure applied to the Hamiltonian flow $\tilde{f}_t$
with the Hamiltonian $\tilde p$. Its result is the transverse
bicharacteristic flow $F^*_t$, which is a noncommutative analogue
of the corresponding reduced Hamiltonian flow on
$N^*{\cF}/{\cF}_N$. See also \cite{Xu} with regard to the
reduction procedure in the noncommutative Poisson geometry.

\begin{ex}
Let $(M,{\mathcal F})$ be a compact Riemannian manifold, equipped
with a Riemannian foliation and a bundle-like metric $g_M$. Let
$H$ be the orthogonal complement to $F=T{\mathcal F}$. Consider
the transverse Laplacian $\Delta_H$. Then the geodesic flow $g_t$
of the Riemannian metric $g_M$ can be restricted to $N^*{\mathcal
F}$, that defines a flow $G_t$, and the transverse geodesic flow
of $\sqrt{\Delta_H}$ is given by the induced action of $G_t$ on
$C^{\infty}_{prop} (G_{{\mathcal F}_N}, |T{\mathcal G}_N|^{1/2})$.

If ${\mathcal F}$ is given by the fibers of a Riemannian
submersion $f:M\rightarrow B$, then there is a natural isomorphism
$N^*_m{\mathcal F}\rightarrow T^*_{f(m)}B$, and, under this
isomorphism, the transverse geodesic flow $G_t$ on $N^*{\mathcal
F}$ corresponds to the geodesic flow on $T^*B$ (cf., for instance,
\cite{ONeil,Re}).
\end{ex}

Now we turn to the Egorov theorem for transversally elliptic
operators. Let $E$ be a Hermitian vector bundle on $M$,
$D\in\Psi^1(M,E)$ a formally self-adjoint transversally elliptic
operator in $L^2(M,E)$. Suppose that $D^2$ has the scalar
principal symbol and the holonomy invariant transversal principal
symbol. $D$ is essentially self-adjoint with the initial domain
$C^\infty(M,E)$ \cite{noncom}. Therefore, the operator $\langle
D\rangle=(D^2+I)^{1/2}$ is well-defined as a positive,
self-adjoint operator in $L^2(M,E)$, in addition, its domain
contains the Sobolev space $H^1(M,E)$.

By the spectral theorem, for any $s\in\RR$, the operator $\langle
D\rangle^s=(D^2+I)^{s/2}$ is a well-defined positive self-adjoint
operator in the space $\cH=L^2(M,E)$, which is unbounded for
$s>0$. For any $s\geq 0$, define by $\cH^s$ the domain of $\langle
D\rangle^s$, and, for $s<0$, put $\cH^s=(\cH^{-s})^*$. Let also
$\cH^{\infty}=\bigcap_{s\geq 0}\cH^s, \quad
\cH^{-\infty}=(\cH^{\infty})^*$. It is clear that $C^\infty(M,E)
\subset \cH^s$ for any $s$.

\begin{defn}\label{d:ideals}
A bounded operator $A$ in $\cH^{\infty}$ belongs to
$\cL(\cH^{-\infty},\cH^{\infty})$
($\cK(\cH^{-\infty},\cH^{\infty})$), if, for any $s$ and $r$, it
extends to a bounded (compact) operator from $\cH^s$ to $\cH^r$,
or, equivalently, the operator $\langle D\rangle^rA\langle
D\rangle^{-s}$ extends to a bounded (compact) operator in
$L^2(M,E)$.

Introduce also the class $\cL^1(\cH^{-\infty},\cH^{\infty})$,
which consists of all operators of class
$\cK(\cH^{-\infty},\cH^{\infty})$ such that, for any $s$ and $r$,
the operator $\langle D\rangle^rA\langle D\rangle^{-s}$ is a trace
class operator in $L^2(M,E)$.
\end{defn}

It is easily seen that $\cL(\cH^{-\infty},\cH^{\infty})$ is an
involutive subalgebra in $\cL(\cH)$, and the classes
$\cK(\cH^{-\infty},\cH^{\infty})$ and
$\cL^1(\cH^{-\infty},\cH^{\infty})$ are its ideals. Observe that
any operator with the smooth kernel belongs to the class
$\cL^1(\cH^{-\infty},\cH^{\infty})$.

By the spectral theorem, the operator $\langle D\rangle$ defines a
strongly continuous group $e^{it\langle D\rangle}$ of bounded
operators in $L^2(M,E)$. Consider the one-parameter group $\Phi_t$
of $\ast$-automorphisms of the algebra ${\mathcal L}(L^2(M,E))$,
defined as
\begin{displaymath}
\Phi_t(T)=e^{i t\langle D\rangle}Te^{-i t\langle D\rangle}, \quad
T\in {\mathcal L}(L^2(M,E)).
\end{displaymath}

Recall that any scalar operator $P\in\Psi^m(M)$, acting on
half-densities, has the subprincipal symbol, which is a globally
defined, degree $m-1$ homogeneous, smooth function on
$T^*M\setminus 0$, given in local coordinates by
\begin{equation}\label{e:subprincipal1}
p_{sub}=p_{m-1}-\frac{1}{2i}\sum_{j=1}^n\frac{\partial^2p_m}{\partial
x_j\partial \xi_j}.
\end{equation}
Note that $p_{sub}=0$, if $P$ is a real self-adjoint differential
operator of even order.

\begin{thm}
\label{Egorov} Let $D\in\Psi^1(M,E)$ be a formally self-adjoint
transversally elliptic operator in $L^2(M,E)$ such that $D^2$ has
the scalar principal symbol and the holonomy invariant transversal
principal symbol. Then:
\medskip\par
(1) for any $K\in \Psi^{m,-\infty}(M,{\mathcal F},E)$, there is a
$K(t)\in\Psi^{m,-\infty}(M,{\mathcal F},E)$ such that the family
$\Phi_t(K)-K(t), t\in \RR,$ is a smooth family of operators of
class  $\cL^1(\cH^{-\infty},\cH^{\infty})$.
\medskip\par
(2) If, in addition, $E$ is the trivial line bundle, and the
subprincipal symbol of $D^2$ vanishes, then, for any $K\in
\Psi^{m,-\infty}(M,{\mathcal F})$ with the principal symbol $k\in
S^{m}(G_{{\mathcal F}_N},|T{\mathcal G}_N|^{1/2})$ the operator
$K(t)$ has the principal symbol $k(t)\in S^m(G_{{\mathcal
F}_N},|T{\mathcal G}_N|^{1/2})$ given by $k(t)=F^*_t(k)$.
\end{thm}

In \cite{dgtrace}, one proves an analogue of the
Duistermaat-Guillemin formula \cite{DG} for transversally elliptic
operators on a compact foliated manifold $(M,{\mathcal F})$.
Consider the operator $P=\sqrt{A}$ in $C^{\infty}(M)$, where $A$
satisfies the conditions (A1) and (A2). It is proved that, for any
$k\in C^{\infty}_c(G,|T{\cG}|^{1/2})$, the trace of
$R(k)\,e^{itP}$ is well-defined as a distribution on the real line
$\RR$, $\theta_k=\tr R(k)\,e^{itP} \in {\mathcal D}'({\RR})$.

As above, let $f_t$ denote the transverse bicharacteristic flow on
$N^*{\mathcal F}$, which preserves $\cF_N$. We say that $\nu\in
\widetilde{N}^*{\mathcal F}$ is a relatively periodic (with
respect to $\cF_N$) point of the flow $f$ with a relative period
$t$, if $f_{t}(\nu)$ and $\nu$ are on the same leaf of $\cF_N$.

As shown in \cite{dgtrace}, the distribution $\theta_k$ is smooth
outside the set of all relative periods ${\mathcal T}_k$ of the
transverse bicharacteristic flow $f_t$ such that, if $\nu\in
\widetilde{N}^*{\mathcal F}$ is a corresponding relatively
periodic point, then the projection of the point
$(f_t(\nu),\nu)\in G_{\cF_N}$ to $G$ belongs to the support of $k$
(by compactness of the support of $k$, the set ${\mathcal T}_k$ is
discrete). Moreover, if the set of relatively periodic points with
the fixed relative period $t$ is transversally clean (that is a
natural generalization of the similar notion, introduced in
\cite{DG}), then one can show that $\theta_k$ has an asymptotic
expansion in a neighborhood of $t$ and write down a formula for
the leading term of this asymptotic expansion.

\section{Noncommutative spectral geometry}  \label{noncom}
\subsection{Dimension spectrum}\label{s:sdim}
Let $(M,{\mathcal F})$ be a compact foliated manifold. We
introduce a more general, than in Theorem~\ref{triple}, class of
spectral triples associated with transversally elliptic operators
on $M$. Let:
\begin{enumerate}
\item ${\mathcal A}$ is the algebra $C^{\infty}_c(G)$;
\item ${\mathcal H}$ is the Hilbert space $L^2(M,E)$ of square
integrable sections of a holonomy equivariant Hermitian vector
bundle $E$, where the action of $k\in {\mathcal A}$ is given by
the $\ast$-representation $R_E$;
\item $D$ is a first order self-adjoint transversally elliptic operator
with the holonomy invariant transversal principal symbol such that
$D^2$ is self-adjoint and has the scalar principal symbol.
\end{enumerate}

\begin{theorem}\cite{noncom}
\label{triple1} The spectral triple $({\mathcal A},{\mathcal
H},D)$ described above is a finite-dimensional spectral triple of
dimension $q=\codim \cF$.
\end{theorem}

It is easy to see that the spectral triple introduced in
Section~\ref{s:riem} is a particular case of the spectral triples
constructed in Theorem~\ref{triple1}.

Recall (see Definition~\ref{d:dim}) that the notion of dimension
of a spectral triple $({\mathcal A},{\mathcal H},D)$ is given by
the greatest lower bound of all $p$'s such that the operator
$a(D-i)^{-1}, a\in {\mathcal A},$ is an operator from the Schatten
ideal ${\mathcal L}^p({\mathcal H})$. If we look at a geometric
space as a union of parts of different dimensions, this notion of
dimension gives only the maximum of the dimensions of the parts of
this space. To take into account lower dimensional parts of this
geometric space, Connes and Moscovici \cite{Co-M} suggested to
consider as a more correct notion of dimension not a single real
number $d$, but a subset ${\rm Sd}\subset {\mathbb C}$, which is
called the dimension spectrum of this spectral triple. Recall
briefly this definition \cite{Co-M,Sp-view}.

Let $({\mathcal A}, {\mathcal H}, D)$ be a spectral triple.
Consider the operator $\langle D\rangle=(D^2+I)^{1/2}$. Denote by
$\delta$ the (unbounded) differentiation on ${\mathcal
L}({\mathcal H})$ given by
\begin{equation}
\label{derivative} \delta(T)=[\langle D\rangle,T],\quad T\in
\Dom\delta\subset {\mathcal L}({\mathcal H}).
\end{equation}

\begin{defn}
We say that $P\in {\rm OP}^{\alpha}$ if and only if $P\langle
D\rangle^{-\alpha}\in \bigcap_n {\rm Dom}\;\delta^n$. In
particular, ${\rm OP}^{0}=\bigcap_n {\rm Dom}\;\delta^n$.
\end{defn}

Observe that the classes  $\cL(\cH^{-\infty},\cH^{\infty})$,
$\cK(\cH^{-\infty},\cH^{\infty})$ and
$\cL^1(\cH^{-\infty},\cH^{\infty})$ introduced in
Definition~\ref{d:ideals} belong to the domain of $\delta$ and are
invariant under the action of $\delta$. Moreover, they are ideals
in ${\rm OP}^{0}$.

\begin{defn}
We will say that a triple $({\mathcal A}, {\mathcal H}, D)$ is
smooth, if, for any $a\in {\mathcal A}$, we have the inclusions
$a, [D,a] \in  {\rm OP}^{0}$.
\end{defn}

The space ${\rm OP}^0$ is a smooth algebra (see, for instance,
\cite[Theorem 1.2]{Ji}). Recall that, for the spectral triple
associated with a smooth manifold $M$, ${\rm OP}^0\cap C(M)$
coincides with $C^\infty(M)$. Therefore, informally speaking, the
smoothness condition for a spectral triple $({\mathcal A},
{\mathcal H}, D)$ means that $\cA$ consists of smooth functions on
the corresponding geometric space.

Let $({\mathcal A}, {\mathcal H}, D)$ be a smooth spectral triple.
Denote by ${\mathcal B}$ the algebra generated by all elements of
the form $\delta^n(a)$, where $a\in{\mathcal A}$ and $n\in
{\mathbb N}$. Thus, $\cB$ is the smallest subalgebra in ${\rm
OP}^{0}$, which contains $\cA$ and is invariant under the action
of $\delta$.

\begin{defn}
A spectral triple $({\mathcal A},{\mathcal H},D)$ has the discrete
dimension spectrum ${\rm Sd}\subset{\mathbb C}$, if ${\rm Sd}$ is
a discrete subset in ${\mathbb C}$, the triple is smooth, and, for
any $b\in {\mathcal B}$, the distributional zeta-function
$\zeta_b(z)$ of $\langle D\rangle $ given by $$ \zeta_b(z)=\tr
b\langle D\rangle^{-z}, $$ is defined in the half-plane
$\{z\in\CC:{\rm Re}\, z>d\}$ and extends to a holomorphic function
on ${\mathbb C}\backslash {\rm Sd}$ such that the function
$\Gamma(z)\zeta_b(z)$ is rapidly decreasing on the vertical lines
$z=s+it$ for any $s$ with ${\rm Re}\,s
>0$.

The dimension spectrum is said to be simple, if the singularities
of $\zeta_b(z)$ at $z\in {\rm Sd}$ are at most simple poles.
\end{defn}

\begin{ex}\label{ex:classic}
Let $M$ be a compact manifold of dimension $n$, $\cE$ a vector
bundle on $M$ and $D\in\Psi^1(M,\cE)$ a self-adjoint elliptic
operator. Then the triple $(C^{\infty}(M),L^{2}(M,\cE),D)$ is a
smooth $n$-dimensional spectral triple. The algebra $\mathcal B$
is contained in the algebra $\Psi^0(M, \cE)$ of zero order
pseudodifferential operators. For an operator $P\in
\Psi^m(M,\cE)$, $m\in\ZZ$, the function $z \rightarrow \tr
P|D|^{-z}$ has a meromorphic extension to $\CC$ with at most
simple poles at integer points $k$, $k\leq m+n$. The residue of
this function at $z=0$ coincides with the Wodzicki-Guillemin
residue $\tau(P)$ of $P$ (see Section~\ref{s:riem}):
\begin{equation}
\underset{z=0}{\res} \tr P|D|^{-z}=\tau(P).
     \label{eq:LIFNCG-NCR}
\end{equation}

In particular, this spectral triple has the simple discrete
dimension spectrum, which is contained in $\{k \in \ZZ : \ k\leq
n\}$.
\end{ex}

An immediate consequence of the technique described in
Section~\ref{transpdo}, in particular, of Theorem~\ref{zeta0}, is
a description of the dimension spectrum for the spectral triples
associated with transversally elliptic operators.

\begin{theorem}\cite{noncom}
\label{dim} The spectral triple $({\mathcal A},{\mathcal H},D)$
introduced in Theorem~\ref{triple1} has the discrete dimension
spectrum ${\rm Sd}$, which is contained in $\{v\in {\mathbb
N}:v\leq q\}$ and is simple.
\end{theorem}

\subsection{Noncommutative pseudodifferential calculus}
Since the algebra ${\mathcal A}$ for the spectral triples
associated with transversally elliptic operators on a foliated
manifold is not unital, it is natural to consider the space
$M/{\mathcal F}$ to be noncompact. Thus, it is necessary to take
into account the behavior of geometrical objects at ``infinity''.
In \cite{Co-M,Sp-view}, the definition of the algebra
$\Psi^*({\mathcal A})$ of pseudodifferential operators is given
for a unital algebra $\cA$. In this Section, following the paper
\cite{egorgeo}, we introduce the algebra $\Psi^*_0({\mathcal A})$,
which can be considered as an analogue of the algebra of
pseudodifferential operators on a noncompact manifold, whose
symbols vanish at infinity along with all the derivatives of an
arbitrary order. First, we introduce some auxiliary notions.

Let $({\mathcal A},{\mathcal H},D)$ be a smooth spectral triple.
Denote by ${\rm OP}_0^{0}$ the space of all $P\in {\rm OP}^{0}$
such that $\langle D\rangle^{-1}P$ and $P\langle D\rangle^{-1}$
are compact operators in $\cH$. We also say that $P\in {\rm
OP}_0^{\alpha}$, if $P\langle D\rangle^{-\alpha}$ and $\langle
D\rangle^{-\alpha}P$ belong to ${\rm OP}_0^{0}$. It is easy to see
that ${\rm OP}_0^{-\infty}=\bigcap_{\alpha}{\rm OP}_0^{\alpha}$
coincides with ${\mathcal K}({\mathcal H}^{-\infty}, {\mathcal
H}^\infty)$. If the algebra $\cA$ is unital, then ${\rm
OP}_0^{\alpha}= {\rm OP}^{\alpha}$.

We will assume that the subalgebra ${\mathcal B}\subset {\mathcal
L}({\mathcal H})$, generated by the operators $\delta^n(a), a\in
{\mathcal A}, n\in {\NN}$, is contained in ${\rm OP}_0^{0}$. In
particular, this implies that $({\mathcal B},{\mathcal H},D)$ is a
spectral triple in the above sense.

By the definition of a spectral triple, $\cA\subset {\rm
OP}^{0}_0$. Informally speaking, this means that $\cA$ consists of
smooth functions on the corresponding geometric space, vanishing
at infinity. The condition $\cB\subset {\rm OP}^{0}_0$ means that
elements of $\cA$ vanish at infinity along with all its
derivatives of arbitrary order.

\begin{defn}
We say that an operator $P$ in $\cH^{-\infty}$ belongs to the
class $\Psi^*_0({\mathcal A})$, if it admits an asymptotic
expansion:
\begin{equation}
\label{psif} P\sim \sum_{j=0}^{+\infty}b_{q-j}\langle
D\rangle^{q-j}, \quad b_{q-j}\in {\mathcal B},
\end{equation}
that means that, for any $N$
\[
P - \left(b_q\langle D\rangle^q + b_{q-1}\langle
D\rangle^{q-1}+\ldots+b_{-N}\langle D\rangle^{-N}\right)\in {\rm
OP}_0^{-N-1}.
\]
\end{defn}

Using a slight modification of the proof of Theorem~B.1 in
\cite[Appendix B]{Co-M}, it can be shown that $\Psi^*_0({\mathcal
A})$ is an algebra.

\begin{ex}
For the spectral triple $(C^{\infty}(M),L^{2}(M,\cE),D)$ described
in Example~\ref{ex:classic}, the algebra $\Psi^*_0({\mathcal A})$
is contained in $\Psi^*(M, \cE)$.
\end{ex}

Let us give a description of the noncommutative pseudodifferential
calculus for the spectral triples associated with transversally
elliptic operators, following \cite{egorgeo}.

Let $({\mathcal A},{\mathcal H},D)$ be a spectral triple,
constructed in Theorem~\ref{triple1}. Then it is a smooth spectral
triple. Moreover, the algebra ${\mathcal B}$ in contained in ${\rm
OP}_0^{0}$. In this case, any element $b\in \cB$ can be written as
\[
b=B+T, \quad B\in \Psi_{sc}^{0,-\infty}(M,\cF,E), \quad T\in
\cL^1(\cH^{-\infty},\cH^{\infty}),
\]
and the algebra $\Psi^*_0({\mathcal A})$ is contained in
$\Psi_{sc}^{*,-\infty}(M,{\mathcal F},E)+{\rm OP}_0^{-N}$ for any
$N$.

\subsection{Noncommutative local index theorem}\label{s:locind}
A spectral triple $(\cA, \cH, D)$ defines a Fredholm module
$(\cH,F)$ over $\cA$ (see Section~\ref{s:riem}), and, therefore,
the index map $\operatorname{ind} : K_*(\cA)\to\CC$ (see the
formulas (\ref{eq:LIFNCG.index-even}) and
(\ref{eq:LIFNCG.index-odd})). As shown above, this map can be
expressed in terms of the pairing with the cyclic cohomology class
$\operatorname{ch}_\ast (\cH,F) \in HP^*(\cA)$, the Chern
character of the Fredholm module $(\cH,F)$ (see the formulas
(\ref{e:odd1}) and (\ref{e:even1})). Let us call
$\operatorname{ch}_\ast (\cH,F)$ the Chern character of the
spectral triple $(\cA, \cH, D)$ and denote by
$\operatorname{ch}_\ast (\cA, \cH, D)$. The formulas
(\ref{e:odd1}) and (\ref{e:even1}) have a defect, which consists
in the fact that they give expressions of the map
$\operatorname{ind}$ in terms of the operator traces, but these
traces are nonlocal functionals, and it is impossible to compute
them in coordinate charts. To correct this defect, Connes and
Moscovici proved another formula for the index map, which involves
local functionals of Wodzicki-Guillemin residue trace type.

Suppose that a spectral triple $(\cA,\cH,D)$ is smooth and has the
simple discrete dimension spectrum. Recall that $\cB$ denotes the
algebra generated by the operators of the form $\delta^k(a)$, $a
\in \cA$, $k\in \NN$. Define the noncommutative integral
determined by this spectral triple, setting
\begin{equation}
    \bint b = \underset{z=0}{\res} \Tr b|D|^{-z}, \quad b\in \cB.
\end{equation}
The functional $\bint$ is a trace on $\cB$, which is local in the
sense of noncommutative geometry, since it vanishes for any
element of $\cB$, which is a trace class operator in $\cH$.

\begin{thm}[{\cite[Thm.~II.3]{Co-M}}]\label{t:cmeven}
Suppose that $(\cA,\cH,D)$ is an even spectral triple, which is
$p$-summable and has the simple discrete dimension spectrum. Then:

\begin{enumerate}
  \item The following formulas define an even cocycle
$\varphi_{\mathop{CM}}^\ev =(\varphi_{2k})$ in the
$(b,B)$-bicomplex of $\cA$. For $k=0$,
             \begin{equation}
    \varphi_{0}(a^{0})  = \underset{z=0}{\res}z^{-1}\tr \gamma a^0|D|^{-z},
    \label{eq:chern.connes-moscovici.constant}
    \end{equation}
whereas for $k\neq0$,
    \begin{equation}
        \varphi_{2k}(a^{0}, \ldots, a^{2k}) =   \sum_{\alpha\in
        \NN^n}c_{k,\alpha}\bint \gamma a^0 [D,a^1]^{[\alpha_{1}]} \ldots
          [D,a^{2k}]^{[\alpha_{2k}]} |D|^{-(2|\alpha|+k)},
          \label{eq:chern.connes-moscovici.even}
       \end{equation}
where
\[
c_{k,\alpha}=
\frac{(-1)^{|\alpha|}2\Gamma(|\alpha|+k)}{\alpha!(\alpha_{1}+1)
\cdots (\alpha_{1}+\cdots +\alpha_{2k}+ 2k)},
\]
and the symbol $T^{[j]}$ denotes the $j$-th iterated commutator
with $D^2$.
  \item The cohomology class defined by $\varphi_{\mathop{CM}}^\ev$ in $HP^\ev
(\cA)$ coincides with $\operatorname{ch}_\ast (\cA, \cH, D)$.
\end{enumerate}
\end{thm}

\begin{theorem}[{\cite[Thm.~II.2]{Co-M}}]\label{t:cmodd}
Suppose that $(\cA,\cH,D)$ is a spectral triple, which is
$p$-summable and has the simple discrete dimension spectrum. Then:
\begin{enumerate}
  \item There is defined an odd cocycle $\varphi_{\mathop{CM}}^\odd
=(\varphi_{2k+1})$ in the $(b,B)$-bicomplex of $\cA$ given by
     \begin{multline}
             \varphi_{2k+1}(a^{0}, \ldots, a^{2k+1}) =  \\
            \sqrt{2i\pi}
               \sum_{\alpha\in \NN^n}c_{k,\alpha}
               \bint a^0 [D,a^1]^{[\alpha_{1}]} \ldots
              [D,a^{2k+1}]^{[\alpha_{2k+1}]} |D|^{-2(|\alpha|+k)-1)},
             \label{eq:chern.connes-moscovici.odd}
       \end{multline}
where
\[
c_{k,\alpha}= \frac{(-1)^{|\alpha|}\Gamma(|\alpha|+k
+\frac12)}{\alpha!(\alpha_{1}+1) \cdots (\alpha_{1}+\cdots
+\alpha_{2k+1}+2k +1)}.
\]
  \item The cohomology class defined by $\varphi_{\mathop{CM}}^\odd$
in $HP^\odd (\cA)$ coincides with $\operatorname{ch}_\ast (\cA,
\cH, D)$.
\end{enumerate}
 \end{theorem}

\begin{ex}(\cite{Co-M,Ponge:newpf})
Let $M$ be a compact manifold of dimension $n$ and $D$ a first
order, self-adjoint, elliptic, pseudodifferential operator on $M$,
acting on sections of a vector bundle $\cE$ on $M$. Then the
noncommutative integral $\bint$ defined by the spectral triple
$(C^{\infty}(M),L^{2}(M,\cE),D)$ coincides with the residue trace
$\tau$ (see~(\ref{e:wodz})).

Moreover, in the case when $D$ is the Dirac operator, acting on
sections of the spin bundle $\cS$ on a compact spin Riemannian
manifold $M$, in Theorem~\ref{t:cmeven}, for any $f^0, f^1,\ldots,
f^{m}\in C^{\infty}(M)$, when $|\alpha|\neq 0$, we have
\[
\tau (\gamma f^0 [D,f^1]^{[\alpha_{1}]} \ldots
          [D,f^{m}]^{[\alpha_{m}]} |D|^{-(2|\alpha|+m)})=0,
\]
and, when $\alpha_1=\alpha_2=\ldots=\alpha_{2k}=0$,
\[
\tau (\gamma f^0 [D,f^1] \ldots [D,f^{m}] |D|^{-m}) = c_m\int_M
f^0\,df^1\wedge \ldots df^{m}\wedge \hat{A}(R),
\]
where $c_m$ is some constant, and $\hat{A}(R)$ is the total
$\hat{A}$-form of the Riemann curvature $R$,
$\hat{A}(R)=\left[\det\left(\frac{R/2}{\operatorname{sh}
(R/2)}\right)\right]^{1/2}$.

If $\dim M$ is even, then the spectral triple is even, and the
components of the corresponding even cocycle
$\varphi_{\mathop{CM}}^\ev = (\varphi_{2k})$ are given by
\[
\varphi_{2k}(f^0,\ldots,f^{2k})=\frac{1}{(2k)!}\int_M
f^0\,df^1\wedge \ldots df^{2k}\wedge \hat{A}(R)^{(n-2k)},
\]
where $f^0, f^1,\ldots, f^{2k}\in C^{\infty}(M)$.

If $\dim M$ is odd, then the spectral triple is odd, and the
components of the corresponding odd cocycle
$\varphi_{\mathop{CM}}^\odd = (\varphi_{2k+1})$ are given by
\[
\varphi_{2k+1}(f^0,\ldots,f^{2k})=\sqrt{2i\pi}
\frac{(2i\pi)^{-[\frac{n}{2}]+1}}{(2k+1)!}\int_M f^0\,df^1\wedge
\ldots df^{2k+1}\wedge \hat{A}(R)^{(n-2k-1)},
\]
where $f^0, f^1,\ldots, f^{2k+1}\in C^{\infty}(M)$.
\end{ex}

In~\cite{CoM:Hopf}, Connes and Moscovici computed the index of
transversally hypoelliptic operators associated with triangular
structures on a smooth manifold (see Section~\ref{para}, in
particular, Theorem~\ref{t:para}), using the noncommutative local
index theorem, Theorem~\ref{t:cmodd}. As it was already mentioned
in Section~\ref{para}, the study of these spectral triples makes
an essential use of the hypoelliptic pseudodifferential calculus
on Heisenberg manifolds developed by Beals and Greiner \cite{BG}.
In particular, the noncommutative integral $\bint$ determined by
such a spectral triple coincides with the Wodzicki-Guillemin type
residue trace $\tau$ defined on the Beals-Greiner algebra of
pseudodifferential operators.

The direct computation of the index given by Theorem~\ref{t:cmodd}
for the spectral triples associated with triangular structures on
smooth manifolds is quite cumbersome even in the one-dimensional
case. One gets formulas, involving thousands terms, most of them
give zero contribution. To simplify apriori the computations,
Connes and Moscovici introduced a Hopf algebra $\cH_n$ of
transverse vector fields in $\RR^n$, which plays a role of the
quantum symmetry group. They constructed the cyclic cohomology
$HC^*(\cH)$ for an arbitrary Hopf algebra $\cH$ and a map
\[
HC^*(\cH_n)\to HC^*(C^\infty_c(P)\rtimes\Gamma).
\]
Moreover, they showed that the cyclic cohomology $HC^*(\cH_n)$ are
canonically isomorphic to the Gelfand-Fuchs cohomology
$H^*(W_n,SO(n))$. Therefore, it is defined a characteristic
homomorphism
\[
\chi_{SO(n)}^* : H^*(W_n,SO(n))\to HP^*(C^\infty(P)\rtimes
\Gamma).
\]

The following theorem is the main result of \cite{CoM:Hopf} (cf.
also \cite{CoMEssays} and the survey \cite{skandal-hopf})

\begin{thm}
Let $(\cA,\cH,D)$ be the spectral triple introduced in
Theorem~\ref{t:para}. The Chern character
$\operatorname{ch}_*(\cA,\cH,D)\in HP^*(C^\infty(PW)\rtimes
\Gamma)$ is the image of an universal class $\cL_n\in
H^*(W_n,SO(n))$ under the characteristic homomorphism
$\chi_{SO(n)}^*$:
\[
\operatorname{ch}_*(\cA,\cH,D)=\chi_{SO(n)}^*(\cL_n).
\]
\end{thm}

There is one more computation, illustrating the Connes-Moscovici
local index theorem, which is given in \cite{Perrot}. Let $\Sigma$
be a closed Riemann surface and $\Gamma$ a discrete pseudogroup of
local conformal maps of $\Sigma$ without fixed points. Using
methods of \cite{CoM:Hopf} (the bundle of Kaehler metrics $P$ on
$\Sigma$, hypoelliptic operators, Hopf algebras), in \cite{Perrot}
a spectral triple, which a generalization of the classical
Dolbeaut complex to this setting, is constructed. The Chern
character of this spectral triple as a cyclic cocycle on the
crossed product $C^\infty_c(\Sigma)\rtimes \Gamma$ is computed in
terms of the fundamental class $[\Sigma]$ and of a cyclic
2-cocycle, which is a generalization of the Poincar\'e dual class
to the Euler class. This formula can be considered as a
noncommutative version of the Riemann-Roch theorem.

A Hopf algebra of the same type as the algebra $\cH_1$ was
constructed by Kreimer \cite{Kreimer98} for the study of the
algebraic structure of the perturbative quantum field theory. This
connection was elaborated further in \cite{CoKr1,CoKr2,CoKr3}. One
should also mention the papers by Connes and Moscovici
\cite{CoM:mmj1,CoM:mmj2} on modular Hecke algebras, where, in
particular, it is shown that such an important algebraic structure
on the modular forms as the Rankin-Cohen brackets has a natural
interpretation in the language of noncommutative geometry in terms
of the Hopf algebra $\cH_1$.

\subsection{Noncommutative geodesic flow}\label{s:geo}
Let a spectral triple $({\mathcal A},{\mathcal H},D)$ be given. We
will also suppose that the subalgebra ${\mathcal B}\subset
{\mathcal B}({\mathcal H})$, generated by the operators
$\delta^n(a), a\in {\mathcal A}, n\in {\NN}$, is contained in
${\rm OP}_0^{0}$.

Put ${\mathcal C}_0={\rm OP}_0^0\bigcap \Psi^*_0({\mathcal A})$.
Let $\bar{\mathcal C}_0$ be the closure of ${\mathcal C}_0$ in
${\mathcal L}(\cH)$. For any $T\in \cL(\cH)$, define
\begin{equation}\label{e:alphat}
\alpha_t(T)=e^{it\langle D\rangle}Te^{-it\langle D\rangle}, \quad
t\in {\RR}.
\end{equation}
As usual, ${\mathcal K}$ denotes the ideal of compact operators in
${\mathcal H}$.

\begin{defn}
For the spectral triple $({\mathcal A}, {\mathcal H}, D)$, the
unitary cotangent bundle $S^*{\mathcal A}$ is defined as the
quotient of the $C^*$-algebra, generated by the union of all
spaces of the form $\alpha_t(\bar{\mathcal C}_0)$ with $t\in\RR$
and of $\mathcal K$, by its ideal ${\mathcal K}$.
\end{defn}

\begin{defn}
For the spectral triple $({\mathcal A}, {\mathcal H}, D)$, the
noncommutative geodesic flow is the one-parameter group $\alpha_t$
of automorphisms of the algebra $S^*{\mathcal A}$ defined by
(\ref{e:alphat}).
\end{defn}

\begin{ex}
As shown in \cite{Sp-view} using the classical Egorov theorem, for
the spectral triple $(\cA,{\mathcal H},D)$ associated with a
compact Riemannian manifold $(M,g)$, the $C^*$-algebra
$S^*{\mathcal A}$ is canonically isomorphic to the algebra
$C(S^*M)$ of continuous functions on the cosphere bundle to $M$,
$S^*M=\{(x,\xi)\in T^*M : \|\xi\|=1\}.$ Moreover, the
one-parameter group $\alpha_t$ is given by the natural action of
the geodesic flow on $C(S^*M)$.
\end{ex}

Some examples of the computation of the noncommutative geodesic
flow can be found in \cite{Golse:Leichtnam}.

Theorem~\ref{Egorov} allows to give a description of the
noncommutative flow defined by a spectral triple associated with a
Riemannian foliation in the case when $E$ is the trivial line
bundle (see \cite{egorgeo}). Thus, suppose that $(M,{\mathcal F})$
is a compact foliated manifold and a spectral triple $({\mathcal
A},{\mathcal H},D)$ is of the following form:
\begin{enumerate}
\item The involutive algebra ${\mathcal A}$ is the algebra $C^\infty_c(G,|T\cG|^{1/2})$;
\item The Hilbert space ${\mathcal H}$ is the space $L^2(M)$
with the action of ${\mathcal A}$ given by the
$\ast$-representation $R$;
\item The operator $D$ is a first order, self-adjoint, transversally elliptic
operator with the holonomy invariant transversal principal symbol
such that the subprincipal symbol of $D^2$ vanishes.
\end{enumerate}

\begin{thm}
\label{noncom:flow} For the given spectral triple, the transverse
bicharacteristic flow $F^*_t$ of $\langle D\rangle$ extends by
continuity to a strongly continuous one-parameter group of
automorphisms of the algebra $\bar{S}^0(G_{{\mathcal
F}_N},|T{\cG}_N|^{1/2})$, and there is a $\ast$-homomorphism $P :
S^*{\mathcal A}\rightarrow \bar{S}^0(G_{{\mathcal
F}_N},|T{\cG}_N|^{1/2})$ such that the following diagram commutes:
\begin{equation}\label{e:cd}
  \begin{CD}
S^*{\mathcal A} @>\alpha_t>> S^*{\mathcal A}\\ @VPVV @VVPV
\\ \bar{S}^0(G_{{\mathcal F}_N},|T{\cG}_N|^{1/2})@>F^*_t>>
\bar{S}^0(G_{{\mathcal F}_N},|T{\cG}_N|^{1/2})
  \end{CD}
\end{equation}
\end{thm}

\subsection{Transversal Laplacian and noncommutative diffusion}
In this Section, we describe the results of Sauvageot
\cite{SauCRAS,Sau} on the existence and properties of the
noncommutative diffusion defined by the transversal Laplacian on
the $C^*$-algebra of a Riemannian foliation.

Let $(M,\cF)$ be a compact manifold equipped with a Riemannian
foliation, $g$ a bundle-like metric, $H=F^\bot$. Recall that the
transverse differential $d_H: \Omega_\infty^0 \to \Omega_\infty^1$
has been defined in Section~\ref{s:transdif}. Denote by
$L^2\Omega^1$ the Hilbert completion
$\Omega_\infty^1=C^\infty_c(G, r^*N^*\cF\otimes |T\cG|^{1/2})$ in
the inner product
\[
\langle \omega_1,\omega_2\rangle_{L^2\Omega^1}= {\rm tr}_{\mathcal
F}\;(\langle\omega_1,\omega_2\rangle_{C^*_r(G)}) =
\int_{G}\langle\omega_1(\gamma),\omega_2(\gamma)\rangle_{N^*_{r(\gamma)}\cF}\,dv(\gamma),
\]
where ${\rm tr}_{\mathcal F}$ denotes the trace on the von Neumann
algebra $W^*(M,\cF)$ defined by the transverse Riemannian volume
form, and $dv$ is the transverse Riemannian volume form on $G$.
Consider the operator $d_H$ as a densely defined operator from the
Hilbert space $L^2\Omega^0=L^2(G)$ to the Hilbert space
$L^2\Omega^1$. Then, by Lemma 5.1.1 in \cite{Sau}, the operator
$d_H$ is closable, $\Omega_\infty^1$ belongs to the domain of the
adjoint $d^*_H$ (the transversal divergence operator) and the
transversal Laplacian
\[
\Delta_H=d^*_Hd_H
\]
is defined on $\Omega_\infty^0=C^\infty_c(G, |T\cG|^{1/2})$. The
operator $\Delta_H$ is a positive, self-adjoint operator in
$L^2(G)$. Therefore, it generates a strongly continuous
one-parameter operator semigroup $e^{-t\Delta_H}$ in $L^2(G)$.

Let $\cN$ be a two-sided ideal of the von Neumann algebra
$W^{\ast}(M,{\mathcal F})$ defined as
\[
\cN=\{k\in W^{\ast}(M,{\mathcal F}) : {\rm tr}_{\mathcal
F}\;(k^*k)<\infty\}.
\]
There is a natural isometric embedding of the ideal $\cN$ to
$L^2(G)$, denoted by $\Lambda$ and satisfying the condition
\[
\Lambda (R(k))=k, \quad k\in C^\infty_c(G).
\]

Using the results of \cite{Sau1396,Sau1442}, Sauvageot
\cite{SauCRAS,Sau} showed the existence of an one-parameter
semigroup $\{\Phi_t : t>0\}$ of normal, completely positive
contractions of the von Neumann algebra $W^{\ast}(M,{\mathcal
F})$, mapping $\cN$ to itself and satisfying the condition
\[
\Lambda (\Phi_t(k))=e^{-t\Delta_H}\Lambda (k), \quad k\in \cN.
\]
Moreover, the $C^*$-algebra of the foliation $C^*_r(G)$ is
invariant under the action of $\Phi_t$: $\Phi_t(C^*_r(G))\subset
C^*_r(G)$ for any $t\geq 0$. Finally, it is proved in
\cite{SauCRAS,Sau} that $\{\Phi_t : t>0\}$ is a Markov semigroup,
that is,
\begin{itemize}
  \item for any state $\rho$ on $C^*_r(G)$ and
  for any $t\geq 0$, $\rho\circ \Phi_t$ is a state on $C^*_r(G)$;
\end{itemize}
or, equivalently,
\begin{itemize}
  \item for any $t\geq 0$, we have
\[
\Phi_t(1)=1,
\]
where $1$ is the unit in the von Neumann algebra
$W^{\ast}(M,{\mathcal F})$.
\end{itemize}

\subsection{Adiabatic limits and semiclassical Weyl formula}
Let $(M,{\mathcal F})$ be a closed foliated manifold, $\dim M =
n$, $\dim {\mathcal F} = p$, $p+q=n$, endowed with a Riemannian
metric $g_M$. Let $F=T{\mathcal F}$ be the tangent bundle to $\cF$
and $H=F^{\bot}$ the orthogonal complement of $F$. Thus, the
tangent bundle $TM$ is represented as the direct sum:
\begin{equation}
\label{ad:decomp} TM=F\bigoplus H.
\end{equation}
The decomposition~(\ref{ad:decomp}) induces the decomposition of
the metric $ g_{M}=g_{F}+g_{H}$. Define a one-parameter family
$g_{h}$ of Riemannian metrics on $M$ by
\begin{equation}\label{e:gh}
g_{h}=g_{F} + {h}^{-2}g_{H}, \quad 0 < h \leq 1.
\end{equation}

For any $h>0$, consider the Laplace operator on differential forms
defined by $g_h$:
\begin{equation*}
\Delta_{h}=d^{*}_{g_h}d+dd^{*}_{g_h},
\end{equation*}
where $d :C^{\infty}(M,\Lambda^k T^{*}M)\rightarrow
C^{\infty}(M,\Lambda^{k+1}T^{*}M)$ is the de Rham differential,
$d^{*}_{g_h}$ is the adjoint of $d$ with respect to the inner
product in $C^{\infty}(M,\Lambda T^{*}M)$ defined by $g_{h}$.
$\Delta_{h}$ is a self-adjoint, elliptic, differential operator
with the positive, scalar principal symbol in the Hilbert space
$L^2(M,\Lambda T^{*}M,g_h)$ of square integrable differential
forms on $M$, endowed with the inner product induced by $g_h$,
which has discrete spectrum.

In \cite{adiab}, the asymptotic behavior of the trace of
$f(\Delta_h)$ when $h\to 0$ was studied for any $f\in S(\RR)$.
Such asymptotic limits are called adiabatic limits after Witten.

Recall that the decomposition~(\ref{ad:decomp}) induces a
bigrading on $\Lambda T^{*}M$ by
\begin{equation}
\label{ad:bigrad} \Lambda^k
T^{*}M=\bigoplus_{i=0}^{k}\Lambda^{i,k-i}T^{*}M, \quad
\Lambda^{i,j}T^{*}M=\Lambda^{i}F^{*}\bigotimes \Lambda^{j}H^{*}.
\end{equation}

Introduce a bounded operator $\Theta_h$ in $L^{2}(M,\Lambda
T^{*}M)$ by the formula
\begin{equation}\label{e:thetah}
\Theta_{h}u = h^{j}u, \quad u \in L^{2}(M,\Lambda^{i,j}T^{*}M ,
g_{h}).
\end{equation}
It is easy to see that $\Theta_h$ is an isomorphism of Hilbert
spaces $L^{2}(M,\Lambda T^{*}M , g_{h})$ and $L^{2}(M,\Lambda
T^{*}M ,g)$. Using $\Theta_h$, one can transfer our considerations
in the fixed Hilbert space $L^{2}(M,\Lambda T^{*}M, g)$. The
operator $\Delta_h$ considered as an operator in $L^{2}(M,\Lambda
T^{*}M , g_{h})$ corresponds by the isometry $\Theta_{h}$ to the
operator
\begin{equation*}
L_{h}= \Theta_{h}\Delta_h\Theta_{h}^{-1}
\end{equation*}
in $L^2(M,\Lambda T^{*}M)=L^{2}(M,\Lambda T^{*}M ,g)$.

We will use the notation introduced in Section~\ref{s:do}. One can
show that $d_h = \Theta_{h}d\Theta_{h}^{-1} = d_F + hd_H +
h^{2}\theta$, and the adjoint of $d_h$ in  $L^2(M,\Lambda T^{*}M)$
is $\delta_h=\Theta_{h}\delta \Theta_{h}^{-1} = \delta_F + h
\delta_H + h^{2}\theta^{*}$. Therefore, one has
\begin{equation*}
\begin{aligned}
L_h & = d_h\delta_h + \delta_h d_h \\ & =\Delta_F + h^2\Delta_H +
h^4\Delta_{-1,2}+ hK_1+h^2K_2 +h^3K_3.
\end{aligned}
\end{equation*}

Suppose that ${\mathcal F}$ is a Riemannian foliation and $g_{M}$
is a bundle-like metric. Then $K_1\in D^{0,1}(M, {\mathcal F},
\Lambda T^{*}M )$. Due to this fact, one can show that the leading
term of the asymptotic of the trace of $f(\Delta_h)$ or, that is
the same, of the trace of $f(L_h)$ as $h\to 0$ coincides with the
leading term of the asymptotic of the trace of $f(\bar{L}_h)$ as
$h\to 0$, where
\[
\bar{L}_h=\Delta_F + h^2\Delta_H.
\]
Observe that the operator $\bar{L}_h$ can be considered as a
Schr\"odinger operator on the leaf space $M/\cF$, where $\Delta_H$
plays a role of the Laplace operator, and $\Delta_F$ a role of the
operator-valued potential on $M/\cF$.

Recall that in the case of a Schr\"odinger operator $H_h$ on a
compact manifold $M$ with the operator-valued potential $V\in
C^\infty(M,{\mathcal L}(H))$, where $H$ is a finite-dimensional
Euclidean space and $V(x)^{*}=V(x)$
\[
H_h=-h^2\Delta +V(x),\quad x\in M,
\]
the corresponding asymptotic formula (the semiclassical Weyl
formula) has the following form:
\begin{equation}
\label{ad:semi1} \operatorname{tr}
f(H_h)=(2\pi)^{-n}h^{-n}\int_{T^*M} \operatorname{Tr}
f(p(x,\xi))\,dx\,d\xi+o(h^{-n}),\quad h\rightarrow 0+,
\end{equation}
where $p(x,\xi)$ is the operator-valued principal $h$-symbol of
$H_h$:
\begin{equation*}
p(x,\xi)=|\xi|^2+V(x),\quad (x,\xi)\in T^{*}M.
\end{equation*}
It turns out that the asymptotic formula for the trace of
$f(\Delta_h)$ can be written in a form, similar
to~(\ref{ad:semi1}), using the language of noncommutative
geometry. For this, it is necessary, in particular, to replace the
usual integration over $T^{*}M$ and the fiberwise trace
$\operatorname{Tr}$ by the integration in the sense of the
noncommutative integration theory given by the trace
$\operatorname{tr}_{{\mathcal F}_N}$ on the twisted von Neumann
algebra $W^{\ast}(N^*\cF, \cF_N, \pi^{*}\Lambda T^{*}M)$.

Let us call by the principal $h$-symbol of $\Delta_h$ the
tangentially elliptic operator
$\sigma_h(\Delta_h):C^{\infty}(N^*{\mathcal F},\pi^{*} \Lambda
T^{*}M) \rightarrow C^{\infty}(N^*{\mathcal F},\pi^{*} \Lambda
T^{*}M)$ on the foliated manifold $(N^*{\mathcal F},{\mathcal
F}_N)$ given by
\begin{equation*}
\sigma_h(\Delta_h) = \Delta_{{\mathcal F}_N}+g_N,
\end{equation*}
where
\begin{itemize}
  \item $\Delta_{{\mathcal F}_N}$ is the lift of
$\Delta_F$ to a tangentially elliptic operator $\Delta_{{\mathcal
F}_N}$ on the foliated manifold $(N^*{\mathcal F},{\mathcal
F}_N)$, acting in $C^{\infty}(N^*{\mathcal F},\pi^{*} \Lambda
T^{*}M)$;
  \item $g_N$ is the multiplication operator by $g_N \in
C^\infty(N^*\cF)$, where $g_N$ is the Riemannian metric on
$N^*{\mathcal F}$, induced by the Riemannian metric on $M$.
(Observe that $g_N$ coincides with the transversal principal
symbol of $\Delta_H$)
\end{itemize}
We will consider $\sigma_h(\Delta_h)$ as a family of elliptic
operators along the leaves of the foliation ${\mathcal F}_N$. For
any function $f\in C^\infty_c({\mathbb R})$, the operator
$f(\sigma_h(\Delta_h))$ belongs to the $C^*$-algebra $C^*(N^*\cF,
{\mathcal F}_N,\pi^{*}\Lambda T^{*}M)$. Moreover, the transverse
Liouville measure for the symplectic foliation ${\mathcal F}_N$
defines a trace $\operatorname{tr}_{{\mathcal F}_N}$ on the
$C^*$-algebra $C^*(N^*\cF, {\mathcal F}_N,\pi^{*}\Lambda T^{*}M)$,
and the value of this trace on $f(\sigma_h(\Delta_h))$ is finite.

\begin{theorem}[\cite{adiab}]
\label{ad:main} For any function $f\in C^\infty_c({\mathbb R})$,
the asymptotic formula holds:
\[
\operatorname{tr} f(\Delta_{h}) =(2\pi)^{-q}h^{-q}
\operatorname{tr}_{{\mathcal F}_N} f(\sigma_h(\Delta_h))
+o(h^{-q}),\quad h\rightarrow 0.
\]
\end{theorem}

\end{document}